\documentclass[11pt,a4paper]{article}

\title{MAGNET: an open-source library for mesh agglomeration by Graph Neural Networks}
\author{Paola F. Antonietti$^{a}$, Matteo Caldana$^{a,\!\!\!}$
\thanks{Corresponding author: {\tt matteo.caldana@polimi.it}}
,
Ilario Mazzieri$^a$,
Andrea Re Fraschini$^{a}$ \\[0.3cm]
\small\textit{$^a$MOX, Dipartimento di Matematica, 
Politecnico di Milano,}\\
\small\textit{
Piazza Leonardo da Vinci 32,
20133 Milano, Italy
}}

\date{\today}

\usepackage[utf8]{inputenc}
\usepackage[english]{babel}
\usepackage[T1]{fontenc} 
\usepackage{placeins}
\usepackage{graphicx}
\graphicspath{{images/}}
\usepackage{eso-pic} 
\usepackage{subfig} 
\usepackage[font=footnotesize]{caption} 
\usepackage{transparent}

\usepackage{amsmath}
\usepackage{amsthm}
\usepackage{bm}
\usepackage[overload]{empheq}  

\usepackage{tabularx}
\usepackage{longtable} 
\usepackage{colortbl}

\usepackage{algorithm}
\usepackage{algpseudocode}

\usepackage[colorlinks=true,linkcolor=black,anchorcolor=black,citecolor=black,filecolor=black,menucolor=black,runcolor=black,urlcolor=black]{hyperref} 
\usepackage{cleveref}
\usepackage[square, numbers, sort&compress]{natbib} 
\usepackage{appendix}

\usepackage{enumitem}

\usepackage{amsthm,thmtools,xcolor} 
\usepackage{comment} 
\usepackage{fancyhdr} 
\usepackage{lipsum} 
\usepackage{tcolorbox} 

\newcommand{\bea}{\begin{eqnarray}} 
  \newcommand{\eea}{\end{eqnarray}}

\usepackage{listings}
\usepackage{amssymb}
\usepackage{xspace}
\usepackage{tikz}
\usepackage{adjustbox} 
\usepackage[margin=1in]{geometry}

\newcommand{\maggnn}{\texttt{MAGNET}\xspace}
\newcommand{\py}[1]{\texttt{\lstinline!#1!}}

\definecolor{codegreen}{rgb}{0,0.6,0}
\definecolor{codegray}{rgb}{0.5,0.5,0.5}
\definecolor{codepurple}{rgb}{0.58,0,0.82}
\definecolor{backcolour}{rgb}{0.95,0.95,0.92}

\lstdefinestyle{mystyle}{
  backgroundcolor=\color{backcolour},   
  commentstyle=\color{codegreen},
  keywordstyle=\color{magenta},
  numberstyle=\tiny\color{codegray},
  stringstyle=\color{codepurple},
  basicstyle=\ttfamily\footnotesize,
  breakatwhitespace=false,         
  breaklines=true,                 
  captionpos=b,                    
  keepspaces=true,                 
  numbers=left,                    
  numbersep=5pt,                  
  showspaces=false,                
  showstringspaces=false,
  showtabs=false,                  
  tabsize=2
}

\newlength{\bibitemsep}
\newlength{\bibparskip}
\let\oldthebibliography\thebibliography
\renewcommand\thebibliography[1]{%
  \footnotesize
  \oldthebibliography{#1}%
  \setlength{\parskip}{\bibitemsep}%
  \setlength{\itemsep}{\bibparskip}%
}

\begin{document}
  
  \maketitle
  
  \begin{abstract}
We introduce \maggnn, an open-source Python library designed for mesh agglomeration in both two- and three-dimensions, based on employing Graph Neural Networks (GNN). \maggnn serves as a comprehensive solution for training a variety of GNN models, integrating deep learning and other advanced algorithms such as METIS and k-means to facilitate mesh agglomeration and quality metric computation.
The library's introduction is outlined through its code structure and primary features. The GNN framework adopts a graph bisection methodology that capitalizes on connectivity and geometric mesh information via SAGE convolutional layers, in line with the methodology proposed in \cite{antonietti2022refinementpolygridsCNN,antonietti2024polytopalmeshagglomerationgeometrical}. Additionally, the proposed \maggnn library incorporates reinforcement learning to enhance the accuracy and robustness of the model initially suggested in \cite{antonietti2022refinementpolygridsCNN,antonietti2024polytopalmeshagglomerationgeometrical} for predicting coarse partitions within a multilevel framework.
A detailed tutorial is provided to guide the user through the process of mesh agglomeration and the training of a GNN bisection model.
We present several examples of mesh agglomeration conducted by \maggnn, demonstrating the library's applicability across various scenarios. Furthermore, the performance of the newly introduced models is contrasted with that of METIS and k-means, illustrating that the proposed GNN models are competitive regarding partition quality and computational efficiency.
Finally, we exhibit the versatility of \maggnn's interface through its integration with \texttt{lymph}, an open-source library implementing discontinuous Galerkin methods on polytopal grids for the numerical discretization of multiphysics differential problems.
  \end{abstract}

\noindent\textbf{Keywords}: agglomeration, polytopal meshes, graph neural networks, open-source library.\\
\noindent\textbf{MSC subject classification}: 65N22, 65N30, 65N50, 68T07

\section{Introduction}
\label{sec:introduction}

Many problems arising in the numerical solution of Partial Differential Equations (PDEs) involve domains characterized by highly complex shapes, heterogeneous media, moving geometries, fractures, inclusions, and/or immersed boundaries that can make the mesh generation process challenging. The ever-increasing need to compute the numerical solution with enough accuracy and at reasonable computational costs has favored in the last decade the development of polytopal Finite Element Methods, i.e., Finite Element Methods that can support general polygonal and polyhedral (polytopal, for short) meshes. These methods offer a natural and flexible way to facilitate the process of mesh generation and allow to describe the computational domain in detail while employing a lower number of degrees of freedom when compared to classical FEM. Examples of polytopal FEMs include the Virtual Element Method (VEM) \cite{veiga2012VEMbasics, veiga2016VEMpolygonal, veiga2018VEM3dmagnetostatics, tushar2022VEMsmalledges}, the mimetic finite difference method \cite{brezzi2005MFEpolytopal, brezzi2005MFEdiffusion, veiga2014MFEelliptic}, the Polytopal Discontinuous Galerkin (PolyDG) method \cite{bassi2012polyDGagglomeration, antonietti2013hppolyDGcomplicated, antonietti2021PolyDGgeophysicalsimulations, cangiani2021hppolyDGpolytopal},  the Hybrid High-Order method (HHO) \cite{dipietro2015hhodiffusion, dipietro2016hhoreview, dipietro2019hhopolytopal}, and the hybridizable discontinuous Galerkin method \cite{cockburn2008hybridGDsuperconv, cockburn2009hybrifDGunified, cockburn2010hybrifDGerroranalysis}.
One of the many advantages of polytopal methods is that they can incorporate naturally mesh agglomeration. In the context of the numerical solution of PDEs, mesh agglomeration finds many essential applications. First, it can be used to construct coarser grids for the underlying differential problems that employ fewer mesh elements, being an accurate geometric description of complicated geometrical features (e.g., complicated boundaries, inclusions, microstructures, interior barriers, layers). Second, it can be employed in an adaptive mesh framework to coarsen the mesh in regions where the error is already under control. Finally, it is one of the main pillars in the construction of multilevel/multigrid iterative solvers to generate a hierarchy of grids to be employed to accelerate the solution of the resulting algebraic problems \cite{antonietti2017multigridDG, dargaville2020comparisonelementagglomerationalgorithms, feder2024r3mgrtreebasedagglomeration}.

Performing automatic and suitable quality (polytopal) mesh agglomeration is an open field of research. During mesh agglomeration, it is crucial to preserve the quality of the original mesh, as quality deterioration could negatively affect the numerical method in terms of stability and accuracy. However, there are no well-established effective strategies for this task yet, especially when the underlying computational domain involves heterogeneous materials, inclusions, microstructures, or interior barriers that should be preserved to make the coarse mesh consistent with the different problems under investigation.

Most state-of-the-art classical methods for graph partitioning use a multilevel approach \cite{karypis1998metis, chevalier2008scotch, dhillon2007graclusandmultilevel}: the graph is initially recursively coarsened by collapsing together groups of adjacent nodes, obtaining a sequence of progressively smaller graphs until a target size is reached; then, the coarsest graph is partitioned cheaply and uncoarsened by successively projecting the partition onto the finer graphs. Since we are interpolating back to a larger graph, a refining step is necessary after each projection to preserve the quality of the initial partition.
The main advantage of the multilevel approach is being able to efficiently tackle very large graphs by running the partitioning algorithm only on the smallest graph. Examples of available tools for multilevel graph partitioning are METIS \cite{karypis1998metis} and SCOTCH \cite{chevalier2008scotch}. Other graph partitioning methods include spectral methods \cite{barnard1994fast}, which are based on eigenvalue decomposition of the graph Laplacian matrix. These methods produce high-quality partitions, but can be prohibitively expensive on large graphs.\\

In the last years, much effort has been directed towards the application of Machine Learning (ML) techniques to enhance and accelerate numerical methods. ML algorithms have the advantage of automatically extracting information from large datasets, making them ideal for situations where establishing \emph{a-priori} criteria for the task is not feasible due to the number of possible configurations, which is the case for mesh agglomeration. 
Among ML techniques, Graph Neural Networks (GNNs) have gathered great interest due to their ability to handle at the same time both mesh connectivity and geometric/physical information about the underlying cells. GNNs are also flexible, easily incorporating additional properties of the mesh simply by adding new input features during the training process.

In this work, we introduce \maggnn (Mesh Agglomeration by Graph Neural Network), an open-source Python library, released under the GNU Lesser General Public License, version 3 (LGPLv3). \maggnn provides a simple framework for mesh agglomeration in both two- and three-dimensions, it allows to experiment with different neural network architectures, training them, and comparing their performance to state-of-the-art methods like METIS \cite{karypis1998metis} and k-means on standard quality metrics. \maggnn can also be easily integrated with other software, in particular, it already interfaces with \texttt{lymph} \cite{antonietti2024lymph}, an open-source Matlab library for PolyDG methods: this allows an assessment of the agglomerated meshes by solving differential problems on it. 
Our approach extends the results presented in \cite{antonietti2024polygonalmeshagglomeration, antonietti2024polytopalmeshagglomerationgeometrical} and consists of re-framing the problem of mesh agglomeration as a graph partitioning problem by exploiting the dual mesh.
A key innovation of the present work is the adoption of the reinforcement learning (RL) approach introduced in \cite{gatti2022reinforcementlearningforgraphpartitioning} to train the model, offering a novel paradigm to adaptively optimize graph-based partitioning in heterogeneous domains.
One of the main advantages of GNN-based algorithms is the possibility to incorporate the physical properties of the computational domain during mesh agglomeration and thus taking into account heterogeneities of the computational physical domain, embedded microstructures, and/or interior barriers, and preserve them during the agglomeration process. We remark that our library is independent of both the PDE and the polytopal method under investigation, as it solely relies on the underlying grid agglomeration procedure, making it adaptable in principle for integration into any polytopal code. \\

The rest of the paper is structured as follows.
In Section~\ref{sec:mesh_agglomeration}, we frame the problem of mesh agglomeration as a graph partitioning problem. In particular, we describe the GNN approaches used in \maggnn, including a novel reinforcement learning approach to mesh agglomeration. In Section~\ref{sec:library} we describe the main features of the library and its code structure. In~Section~\ref{sec:examples} we provide a brief user guide, proceeding step-by-step in the agglomeration of a mesh, and showcase some agglomerated examples, comparing the performance of the GNNs approach with state-of-the-art methods. Finally, in Section~\ref{sec:lymph_interface} we illustrate the interface with \texttt{lymph}, solving two-dimensional test cases on polygonal meshes agglomerated by \maggnn.

The library is available at the GitHub repository \texttt{\href{https://github.com/lymphlib/magnet}{lymphlib/magnet}} together with the complete documentation and tutorials.

\section{Mesh Agglomeration}
\label{sec:mesh_agglomeration}
In this section, we frame the mesh agglomeration problem as a graph partitioning problem. We also briefly review the state-of-the-art models available in the literature and included in \maggnn. In particular, we introduce the GNN agglomeration algorithm, and we show how reinforcement learning can be used to better improve and fine-tune our models.
\begin{figure}[!t]
  \centering
  \begin{tikzpicture}
    \node (image1) at (0.3,0) {\includegraphics[width=0.47\textwidth]{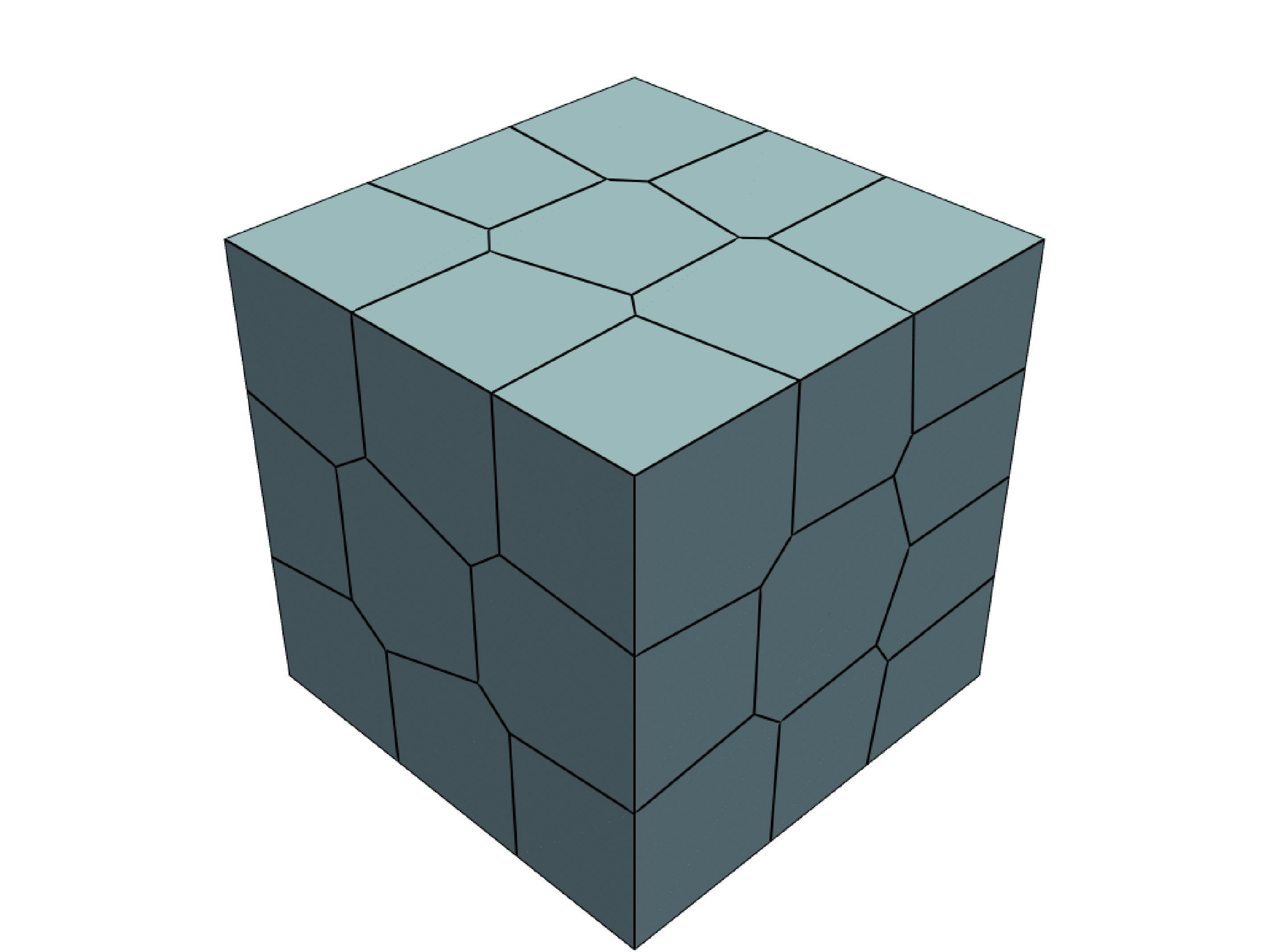}}; 
    \node (arrow1) at (0.25\textwidth,0) {$\rightarrow$};
    \node (image2) at (0.5\textwidth,0) {\includegraphics[width=0.47\textwidth]{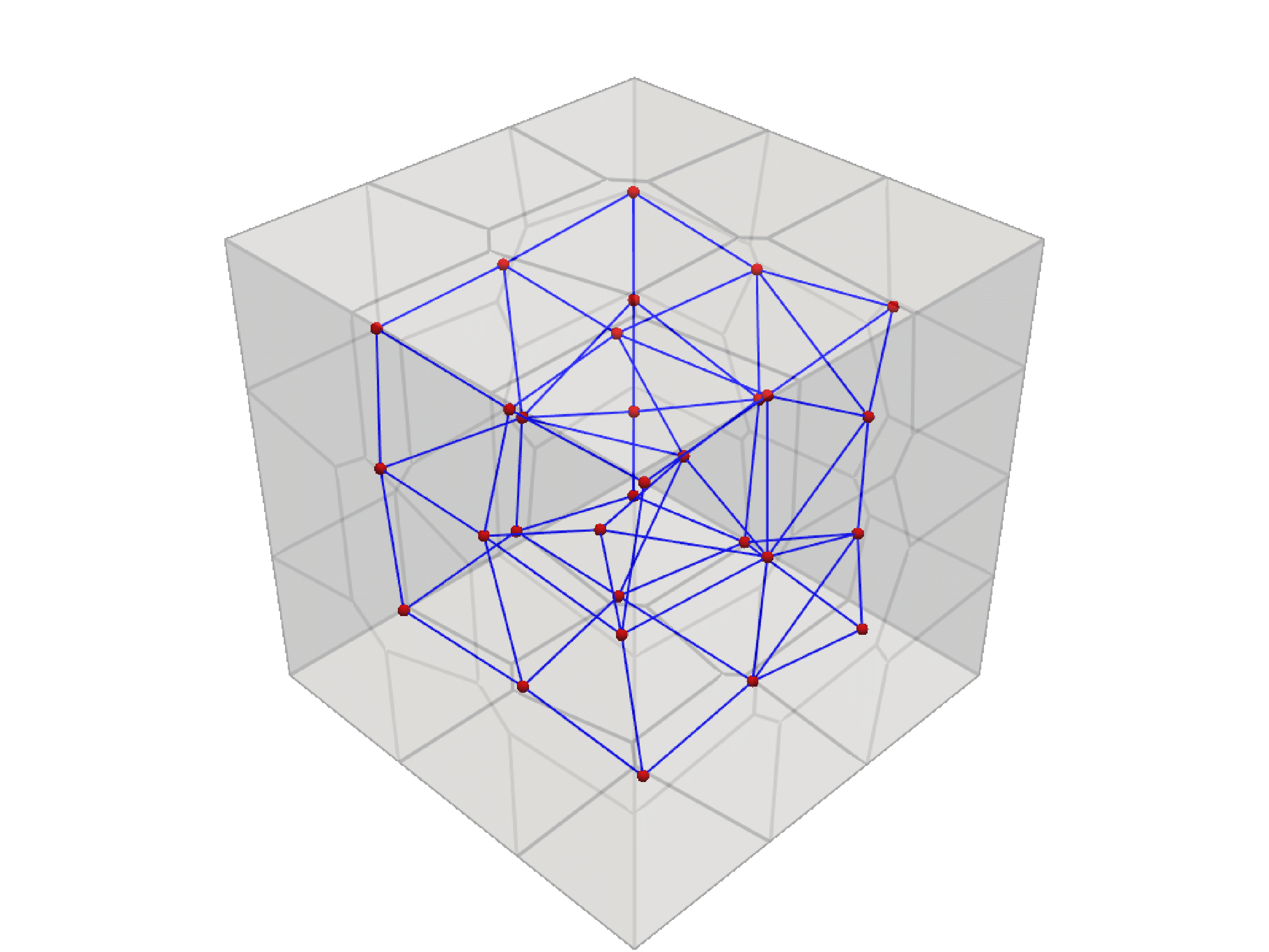}};  
  \end{tikzpicture}
  \caption{Example of a mesh obtained from the centroidal Voronoi tessellation of a cube (\textit{left}) and a representation of its computational graph $\mathcal{G}$: in red the nodes $\mathcal{V}$, in blue the edges $\mathcal{E}$ (\textit{right}).}
  \label{fig:graph_extraction}
\end{figure}

\subsection{Mesh Agglomeration and a Graph Partitioning Problem}
\label{sec:graph_partitioning}
Let us consider an open and bounded computational domain $\Omega \subset \mathbb{R}^d$, $d=2,3$. We consider a discretization $\mathcal{T}_h$ that approximates $\Omega$ composed of disjoint open polytopal elements $P$.
Mesh agglomeration consists of merging some elements $P$ of the mesh $\mathcal{T}_h$ to obtain a coarser (connected) element. Most state-of-the-art methods reframe mesh agglomeration as a graph partitioning problem. Namely, they consider the dual mesh, an undirected graph $\mathcal{G} = (\mathcal{V}, \mathcal{E})$ where each node $v \in \mathcal{V}$ corresponds to an element $P$ of the mesh, and two nodes $v_i, v_j$ are connected by an edge $e = (v_i, v_j) \in \mathcal{E} \subset \mathcal{V} \times \mathcal{V}$ if the corresponding elements are adjacent, i.e., they share an edge in two-dimensions or a face in three-dimensions, see Figure~\ref{fig:graph_extraction}. Mesh agglomeration is then equivalent to partitioning the dual graph.

Given a connected subgraph $\mathcal{S} = (\mathcal{V}_S, \mathcal{E}_S) $, that is a connected graph where $\mathcal{V}_S \subset \mathcal{V}$ and $\mathcal{E}_S = \{(v_i, v_j) \in \mathcal{E}: v_i, v_j \in \mathcal{V}_S\}$, we define its cut and volume respectively as:
  \begin{equation*}
    \text{cut}(\mathcal{S}) = |(v_i,v_j)\in \mathcal{E}: v_i\in \mathcal{V}_S, v_j \in \mathcal{V} \setminus \mathcal{V}_S|, \quad
    \text{vol}(\mathcal{S}) = \sum_{v_i \in \mathcal{V}_S} \text{deg}(v_i),
  \end{equation*}
where $\text{deg}(v)=|(v,v_i)\in \mathcal{E}, \forall \ v_i \in \mathcal{V}|$ is the degree of node $v$, and $|\cdot|$ indicates the cardinality of the set.
The definition of cut can be extended to a generic partition $(\mathcal{S}_1, \dots ,\mathcal{S}_M)$, where $M \in \mathbb{N}$ is the total number of subsets in the partition:
\begin{equation*}
  \text{cut}(\mathcal{S}_1, \dots ,\mathcal{S}_M) = \sum_{i=1}^M \text{cut}(\mathcal{S}_i).
\end{equation*}
Typically, when partitioning a graph, we want to minimize the cut while keeping each set balanced, i.e., the subgraphs should have approximately equal size. In the context of mesh agglomeration, this corresponds to requiring that the interface between different agglomerated elements is small and that they have similar sizes (areas or volumes). So, we consider the normalized cut instead, which penalizes partitions including subgraphs with very different volumes:
\begin{equation}
  \label{eq:normalized_cut}
  \text{NC}(\mathcal{S}_1, \dots, \mathcal{S}_M)=\sum_{i=1}^M \frac{\text{cut}(\mathcal{S}_i)}{\text{vol}(\mathcal{S}_i)}.
\end{equation}
Depending on the application, other notions of volume may be considered, e.g., the number of nodes. The graph partitioning problem can then be formulated as follows: Find $M > 1$ disjoint subgraphs $(\mathcal{S}_1, \dots ,\mathcal{S}_M)$, $\cup_{i=1}^M \mathcal{S}_i = \mathcal{V}$, $\mathcal{S}_i \cap \mathcal{S}_j = \emptyset \ \forall i \neq j$ minimizing the total normalized cut $\text{NC}(\mathcal{S}_1, \dots,\mathcal{S}_M)$. 
The problem of finding a balanced partition that minimizes the cut is NP-complete \cite{garey1990npcompleteness}. As such, current solution algorithms rely on heuristics to find approximate solutions; however, these algorithms tend to be sequential and not easily parallelizable (like Kernigham-Lin (KL) \cite{kernighan1970KLheuristic} and Fiduccia--Mattheyses (FM) \cite{fiduccia1982FMheuristic}), so they cannot fully exploit modern hardware.
On the other hand, deep learning techniques, like graph neural networks, can run efficiently on GPUs, are parallelizable \cite{ma2018GNNparallelism}, and very flexible in their application; this has motivated research on their employment as graph partitioners \cite{gatti2021deeplearningspectralembedding, gatti2022reinforcementlearningforgraphpartitioning}.
  
\subsection{Graph Neural Network for Mesh Agglomeration} 
\label{sec:GNN}
In this section, we introduce the GNN approach to mesh agglomeration; further details can be found in \cite{antonietti2024polygonalmeshagglomeration, antonietti2024polytopalmeshagglomerationgeometrical}.

The mesh agglomeration algorithm we employ consists of recursively applying a single GNN bisection model to the graph extracted from the mesh. Namely, the graph is partitioned by performing a node classification task: taking as input a graph $\mathcal{G}$ and a node feature tensor $\mathbf X \in \mathbb{R}^{N \times F}$, where $N$ is the number of nodes and $F$ is the number of features per node, the GNN model returns a probability tensor $\mathbf Y \in \mathbb{R}^{N \times 2}$, where $Y_{ij}$ is the probability assigned by the GNN that node $v_i$ belongs to the partition $\mathcal{S}_j$. Since the output of the GNN is a probability distribution, the loss function used during the training process is the expected normalized cut, which takes into account the uncertainty in our predictions on node classification:
\begin{equation*}
  \sum_{k=1,2} \frac{ \sum_{i= 1}^N \sum_{(v_i, v_j) \in \mathcal{E}} Y_{ik}(1-Y_{jk})}{\sum_{i=1}^N Y_{ik}\text{deg}(v_i)} = \sum_{i,j=1}^N \big((\mathbf Y \oslash \mathbf Y^T \mathbf D)(1- \mathbf Y)^T\odot \mathbf A\big)_{ij},
\end{equation*}
where $\mathbf D$ is the vector of node degrees $D_i = \text{deg}(v_i)$, $\mathbf A$ is the adjacency matrix of the graph, $\oslash$ and $\odot$ denote the element-wise division and multiplication respectively, and the sum is over all the elements of the resulting matrix (for a more detailed explanation of this framework, we refer to \cite{nazi2019GAPexpectednormcut}). To bisect a graph, a single forward pass of the GNN model is computed, and nodes are assigned to one of the two classes based on which probability is highest in the output probability matrix $\mathbf Y$.
This approach could be easily generalized to any number of partitions instead of just two, but this would require a differently trained GNN, with several outputs equal to the number of classes, for each desired one; by instead bisecting the graph recursively, we can obtain any number of partitions using only one model. Moreover, GNNs have the unique advantage of being able to process the connectivity together with geometric data that has been embedded into the graph as node features.

\subsubsection{The SAGE-Base Model}
\label{sec:sagebase}
The standard GNN model implemented in \maggnn, that we will call SAGE-Base in the following, employs four SAGE convolutional layers (sample and aggregate \cite{hamilton2018SAGE}) followed by four linear, or dense, layers; the hyperbolic tangent is used after each layer to keep the output in the range $[-1,1]$ and facilitate learning. A final softmax layer ensures that the GNN output is a probability distribution over the two classes for each node. We use as node features the coordinates of the centroids of each mesh element and the measure of each element, so that we have a total of four node features (three in the two-dimensional case). Before feeding them to the GNN, an additional normalization step, consisting of normalizing coordinates to zero mean and unit variance and rotating the mesh so that the widest direction is aligned with the x-axis, is performed to reduce variability in input.
Since graphs are dimensionless objects, the approach is the same for polygonal and polyhedral meshes: only the number of input features and the normalization procedure have to be adapted between the two. 
The two-dimensional version of this model has been trained on a set of 800 meshes comprising structured quadrilaterals, structured triangles, random Delaunay triangulations, and random Voronoi tessellations in equal measure. The training has been performed using the Adam optimizer \cite{kingma2017adam} for 300 epochs with learning rate $\gamma =10^{-5} $, weight decay parameter $\lambda = 10^{-5}$ and batch size 4. The three-dimensional version has a greater number of parameters to account for the increased variability of the three-dimensional case, doubling the size of SAGE and linear layers from 64 and 32 to 128 and 64 units, respectively.
The training dataset comprised 400 tetrahedral meshes, of which 100 were of the unit cube and 300 were of randomly generated portions of it (see Figure~\ref{subfig:cube_portion}); other changes include using a learning rate of $\lambda = 10^{-4}$ and lengthening training to 400 epochs. In both cases, meshes have been randomly rotated during training to artificially augment the dataset.
  
\subsubsection{The SAGE-Heterogeneous Model}
\label{sec:sagehetero}
The SAGE-Base model can be easily extended to include additional node features; in particular, it is possible to perform agglomeration of heterogeneous meshes by including a parameter taking values in the interval $[0,1]$ that describes the heterogeneity of the mesh, which we will refer to as the physical group. By updating the loss function to include a penalty term for agglomerating elements with very different physical parameters, the GNN can learn to partition it while respecting the heterogeneity of the physical parameters. We choose such a term in this way
\begin{equation}
  \frac{a}{|\mathcal{V}|}\sum_{i=1}^{N}\sum_{j=1}^2 \big( \mathbf  P \odot \mathbf Y \big)_{ij},
  \label{eq:sagehetero-loss}
\end{equation}
where $\mathbf P \in \mathbb{R}^{N \times 2}$ is a matrix that contains the physical parameter $p$ and its complementary $1-p$ on each row, $\odot$ denotes element-wise multiplication, $|\mathcal{V}|$ is the number of nodes in the graph and $a$ is a suitable coefficient chosen by performing a hyperparameter search; the term is normalized with respect to graph size to make it scale invariant. 
While this model works only when there are at most two different physical parameters, this approach could be extended to an arbitrary number of different physical groups by adding as node features a one-hot encoding of the various heterogeneous parts; however, this would require a different GNN for each total number of physical groups.
The obtained model is referred to as SAGE-Heterogeneous in the code. The bisection algorithm is the same since the output matrix $\mathbf Y$ has the same meaning.
The SAGE-Heterogeneous model included in \maggnn has the same architecture as the SAGE-Base model and was trained using mostly the same hyperparameters. The two-dimensional version has been trained for 200 epochs on a dataset of 600 meshes of the unit square, of which 200 are homogeneous, 200 are divided into two parts along a line, and 200 have 6-12 circular inclusions.
The three-dimensional version has been trained on a set of 400 meshes of the unit cube divided into one to four physically heterogeneous parts each for 150 epochs. 
  
\subsection{Reinforcement Learning}
\label{sec:reinforcementlearn}
In this section, we aim to introduce the main tools of reinforcement learning and show how to use them to train the GNN model, following the work of \cite{gatti2022reinforcementlearningforgraphpartitioning}. For a complete overview of RL, we refer the interested reader to \cite{sutton2018reinforcementlearning}.

Reinforcement learning considers the sequential decision process of an agent (or actor) aiming
to optimize its interaction with the environment. 
At every time step $t \in \mathbb{N}$, the agent observes a state $s_t \in \mathbb{R}^{d_S}$ and takes accordingly a feasable action $a_{t} \in \mathcal{A} \subset \mathbb{R}^{d_A}$ according to a stochastic policy $\pi_t(a | s_t)$. The environment returns to the agent a new state $s_{t+1}$ and a scalar reward $r_{t+1}$.
The decision process is assumed to be Markov, that is, the policy only depends on the current state $s_t$ and not the previous ones. The objective of the agent is to maximize the cumulative discounted reward:
\begin{equation*}
  R_t = \sum_{k = 0}^{T-t-1} \gamma^k r_{t+1+k},
\end{equation*}
where $\gamma \in (0,1]$ is the discount factor, which describes how far into the future rewards are taken into account by the agent. The value 
\begin{equation*}
  V^\pi(s) = \mathbb{E}_\pi[R_t|s_t=s]
\end{equation*}
is the expected return for starting in state $s$ and following policy $\pi$, which is a measure of how ``good'' state $s_t$ is.
In deep reinforcement learning, neural networks are employed to approximate a suitable policy that maximizes $R_t$. Among a very wide variety of available algorithms, we employ the so-called synchronous advantage actor-critic (A2C) policy gradient approach \cite{mnih2016asynchronous}. A2C is a model-free reinforcement learning algorithm that employs a neural network to approximate the value function $V$, together with $\pi$. Namely, the neural network is formed by two heads, an actor and a critic. Given a state $s_t$, the actor generates a probability distribution over the actions (that is, a discrete vector of probabilities if $\mathcal{A}$ is finite or parameters of a distribution otherwise). The critic branch of the network is then used to estimate the value $V(s_t)$ of the state $s_t$. The value is then used to compute the advantage $(R_t - v(s_t))$, which is then used to reinforce the action taken by the agent. In the A2C framework, actor and critic parameters are updated simultaneously, leading to a generally more stable training procedure. The neural network parameters $\boldsymbol \theta$ are updated in the direction of the gradient of the policy log-probabilities 
$$\mathbb{E}_{\pi_\theta}[\nabla_{\boldsymbol \theta} \log \pi (a_t|s_t) (R_t-v(s_t))],$$
where the advantage is used instead of directly $R_t$ to stabilize training.
An additional hyperparameter $\alpha \in (0,1]$ dictates how fast the critic learns with respect to the actor. Since RL algorithms heavily depend on implementation-level details \cite{henderson2018deep}, we refrain from describing how the gradient updates are computed in practice and direct the reader to the original paper \cite{mnih2016asynchronous} for a complete overview.

\subsubsection{The Reinforcement Learning  Partitioner}
\label{sec:rlpartitioner}
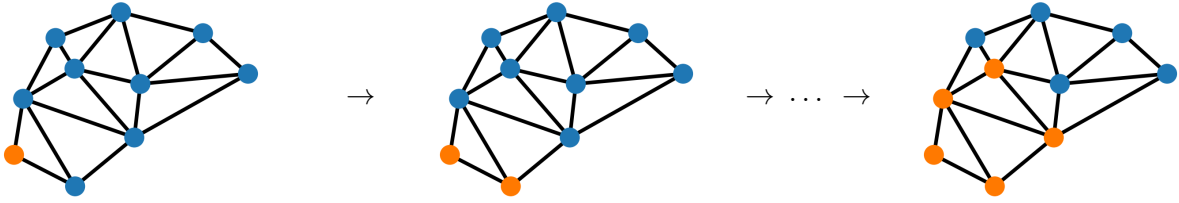
\begin{figure}[!t]
  \centering
  \begin{tikzpicture}
    \node (image1) at (0,0) {\includegraphics[width=0.27\textwidth]{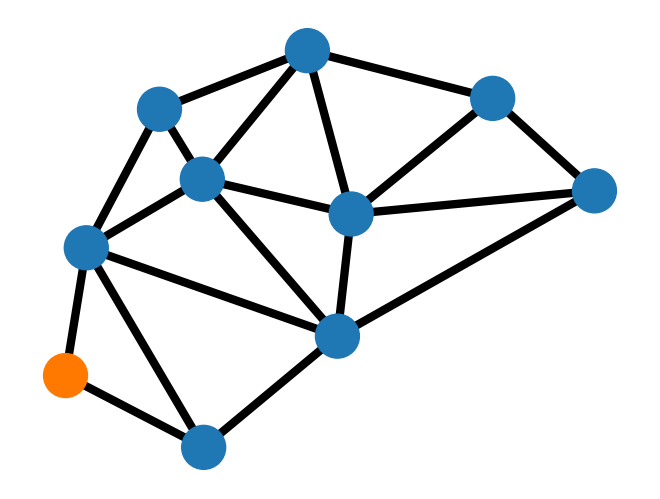}}; 
    \node (arrow1) at (0.15\textwidth,0) {$\rightarrow$};   
    \node (image2) at (0.30\textwidth,0) {\includegraphics[width=0.27\textwidth]{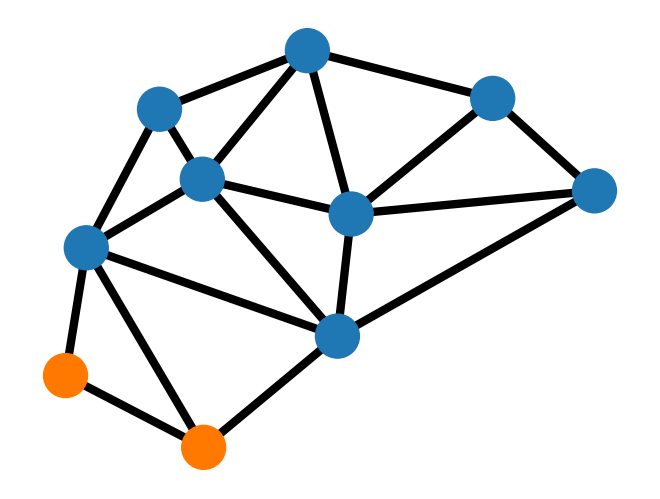}};   
    \node (arrow2) at (0.46\textwidth,0) {$\rightarrow$};    
    \node (dots) at (0.50\textwidth,0) {$\dots$};    
    \node (arrow3) at (0.53\textwidth,0) {$\rightarrow$};    
    \node (image3) at (0.68\textwidth,0) {\includegraphics[width=0.27\textwidth]{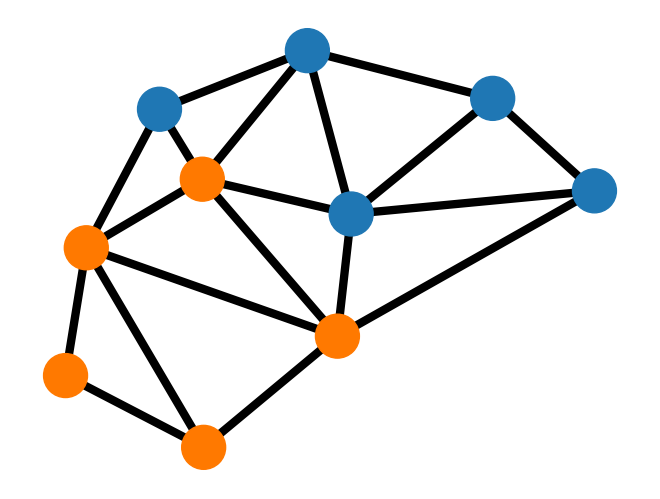}};
 \end{tikzpicture}
  \caption{Representation of how the reinforcement learning agent bisects the graph. First, all nodes are in subset $\mathcal{S}_1$ (blue) except for one node with minimum degree (the orange one). At each step, the agent picks a blue node and flips it to an orange one until the two sets have the same number of nodes.}
  \label{fig:rlpart_illustration}
\end{figure}

The reinforcement learning environment for mesh agglomeration that we propose is an extension of the one presented in \cite{gatti2022reinforcementlearningforgraphpartitioning} for graph partitioning to the specific application of mesh agglomeration.
The objective of the actor is to partition the graph extracted from the mesh into two subgraphs $\mathcal{S}_1, \mathcal{S}_2$ while minimizing the normalized cut Eq.~\eqref{eq:normalized_cut}. 
The state corresponds to the partition of the graph, which is represented in the node features tensor $\mathbf X \in \mathbb{R}^{N \times 2}$ by a one-hot encoding: $X_i = [1,0]$ if $v_i \in \mathcal{S}_1$, while $X_i =[0, 1]$ if $v_i \in \mathcal{S}_2$. Each action, chosen according to the policy of the actor, corresponds to moving one node from $\mathcal{S}_1$ to $\mathcal{S}_2$ or vice versa, modifying the underlying feature tensor. The actor's policy is thus a probability tensor over all the nodes of the graph.
We employ as a reward for each action its corresponding decrease in normalized cut, since this is the quantity we want to minimize. For simplicity, initially, all nodes are in subgraph $\mathcal{S}_1$, except one node which is chosen among the ones with minimum degree so that the cut is as low as possible, and the actor can only move nodes from $\mathcal{S}_1$ to $\mathcal{S}_2$. This procedure is illustrated in Figure~\ref{fig:rlpart_illustration}.
  
It is possible to further customize the learning process by including additional rewards or penalties, e.g., for having non-connected subgraphs.
To improve the quality of the partition, we include among the features the centroid coordinates of each cell, so that the actor may also exploit geometric data to inform its policy.
The length of the episode is chosen as $|\mathcal{V}|/2$, that is the actor moves nodes from $\mathcal{S}_1$ to $\mathcal{S}_2$ until they have the same cardinality; in this way, at the end of the episode we naturally obtain a balanced partition.
Like for the SAGE-Base model, a partition is obtained by recursively applying this bisection model to the mesh, except now a forward pass of the GNN is needed for each action in the episode instead of only once per bisection. We note that the computational cost could be reduced by selecting multiple actions from the same policy, effectively reducing the number of forward passes needed to perform the task.\\
  
The model is implemented in \maggnn using four convolutional SAGE layers and two linear layers that are common to both actor and critic; the critic branch is then formed by two further linear layers with a final softmax layer to create a probability distribution over the graph nodes; to enforce the constraint that chosen nodes can no longer be moved back, we set to minus infinity the input to the softmax layer corresponding to nodes that have already been flipped, so that the probability that they will be chosen is zero. The critic branch uses an Attentional Aggregation layer \cite{li2019AttentionalAggregation} followed by two linear layers that taper off in dimension, finally returning a singular scalar corresponding to the value estimate. All layers, except for the last ones, are interlaced with a hyperbolic tangent activation function.
The pre-trained RL Partitioner model bundled with \maggnn has been trained on a dataset of 1000 meshes (with the same composition as the one of the two-dimensional SAGE-Base model) for one epoch, updating the model parameters at the end of each episode using the Adam optimizer with a learning rate of $10^{-3}$. The discount factor $\gamma$ is chosen to be 0.9, since we are not interested in a greedy approach but rather the long-term result of the episode, while $\alpha$ is set to 0.1 to let the actor learn faster than the critic.
  
\subsubsection{The Reinforcement Learning  Refiner}
\label{sec:rlrefiner}
As we explained in Section~\ref{sec:state_of_the_art}, most state-of-the-art methods use a multilevel approach, which requires a refinement step after each successive uncoarsening of the graph. This process can be tackled by using an actor-critic agent similar to the one presented for graph partitioning, with the difference that nodes can now be moved back and forth between the two subgraphs instead of only progressively growing one of the two from zero, and that the starting state is the projected coarse partition. 
Assuming that the bisection algorithm has already achieved a good quality partition on the coarser problem, only a few nodes will need to be moved during the refinement process, and these nodes will be close to the interface between the two subgraphs. Consequently, we choose a short episode length equal to the cut itself $\text{cut}(\mathcal{S}_1, \mathcal{S}_2)$, and instead of running the model on the entire graph we consider only the $k$-hop subgraph around the cut (that is the subgraph obtained by taking all the nodes with distance at most $k$ from a given center), with $k = 2, 3, 4$. This enhancement significantly reduces the computational cost. The loss of the global context means that two node features indicating the current partition volume balance need to be added to the feature tensor as a trade-off. 
As before, the reduction in the normalized cut will be the actor's reward. 
Since the initial coarse bisection is likely to already have very balanced volumes, a decrease in cut tends to be more significant and disproportionately preferred by the agent than an increase in partition balance.  To counterbalance this effect, we add a penalty term for imbalanced volumes:
\begin{equation*}
  b \frac{\left(\text{Vol}(\mathcal{S}_1)-\text{Vol}(\mathcal{S}_2)\right)^2}{\text{Vol}(\mathcal{V})},
\end{equation*}
where $b>0$ is a hyperparameter that weights the imbalance penalty. Enforcing a stricter balance requirement is crucial since small imbalances can compound over many uncoarsening and refinement steps, leading to agglomerated elements with very different sizes.
In the RL refiner model implemented in \maggnn, actor and critic have one SAGE and one dense layer each, sharing two additional common SAGE layers; the hyperbolic tangent is used as an activation function, and the actor branch terminates with a softmax layer as usual. The refiner architecture can afford to be very light due to the limited task it needs to perform and lower variability in starting configurations.
The model has been trained on a dataset of around 5000 meshes generated starting from a kernel of 1000 meshes (of the same type as the ones used for the two-dimensional SAGE-Base model) by recursively coarsening them and adding the coarsened versions to the dataset. The SAGE-Base model of Section~\ref{sec:sagebase} has been used as a coarse partitioner to generate the initial cut. The network parameters are updated every few steps (in this case 8) instead of only at the end of each episode because, by starting at a configuration close to the optimum, almost all episodes will lead to a negative reward and the network will struggle to learn; more frequent parameters updates give the agent more granular information on the value of its decisions, alleviating the issue.
The imbalance penalty used was $b=0.35$; otherwise, the same hyperparameters of the RL Partitioner model were used.

\subsection{Other Agglomeration Strategies Available in \maggnn}
\label{sec:state_of_the_art}

In addition to the GNN strategies described in the previous Section~\ref{sec:GNN} and Section~\ref{sec:reinforcementlearn}, \maggnn also includes both METIS and k-means as possible agglomeration approaches.
METIS \cite{karypis1998metis} is a widely used library for computing k-way partitions of large irregular graphs based on a multilevel procedure. It uses simple and efficient algorithms (like greedy approaches) for partitioning the coarsest graph, while the refinement is performed based on local heuristics like KL and FM, while respecting balancing constraints. To agglomerate a mesh, METIS is run with the requirement of creating connected subgraphs; also, we use as node weights the volumes of the corresponding cells so that the balance constraint in the partitioning algorithm is imposed in terms of geometric size, which is more relevant to our application than simple node cardinality.

K-means \cite{macqueen1967kmeans} is a quantization algorithm that clusters a set of points in $\mathbb{R}^N$ into $k$ sets by minimizing the squared Euclidean distances of each point to the centroid of its group, i.e.\ by minimizing the within-group variance. While it cannot be applied directly to graphs because they are non-Euclidean data, k-means can be used with node embedding techniques to map each node onto a feature space, producing an effective clustering of the nodes. Indeed, instead of applying k-means to the mesh graph, we use it to directly cluster the cells by using as features the centroid coordinates of the mesh elements, producing fairly rounded agglomerated elements of similar size. In this way, we are neglecting the mesh connectivity, but this is usually not an issue because for meshes with similarly sized elements, centroid coordinates are strongly correlated with two cells being adjacent.
  
\section{\maggnn: Library Overview}
\label{sec:library}
In this section, we describe the tools and design principles of \maggnn.
The main design goal of \maggnn is to provide a single flexible framework for comparing the performance of different methods for mesh agglomeration, paying particular attention to emerging Machine Learning approaches.
\maggnn is able to agglomerate two-dimensional and three-dimensional meshes, including ones obtained from mesh generation of heterogeneous domains, for which agglomeration should preserve the underlying distribution of physical parameters.
Functionalities of the package include: basic input-output of meshes, GNN training, generation of meshes for testing/training, agglomeration by various methods, computation of quality metrics for agglomerated meshes, and visualization.
  
\subsection{Code Structure}
\label{sec:code_structure}
\begin{figure}
  \centering
  \includegraphics[width=1.0\textwidth]{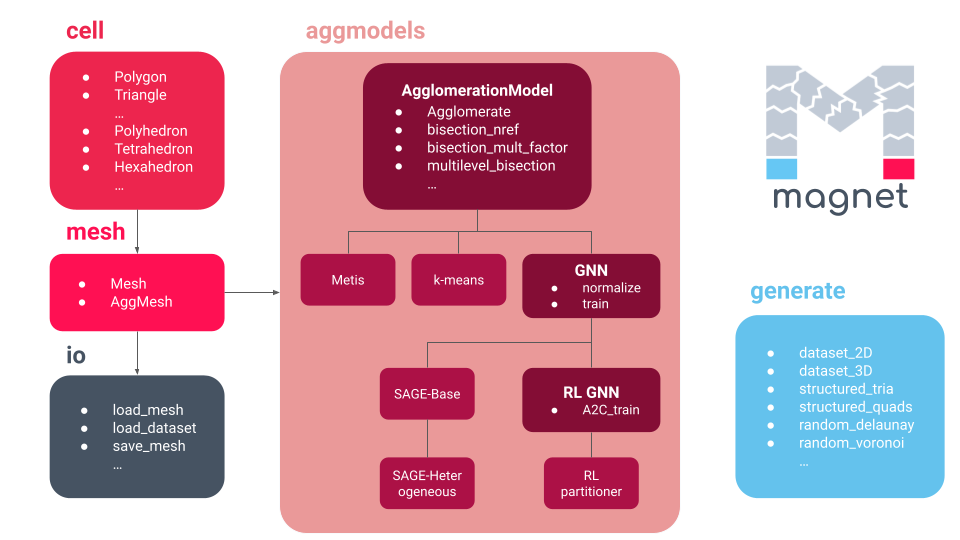}
  \caption{Code structure of \maggnn, highlighting the agglomeration models class hierarchy. Arrows denote import relations.}
  \label{fig:package_structure_diagram}
\end{figure}

\maggnn is organized as reported in Figure~\ref{fig:package_structure_diagram}, namely it is composed by the following sub-modules:
\begin{itemize}
  \item \py{mesh}, \py{cell}: defines the mesh and cell data structures;
  \item \py{aggmodels}: contains the agglomeration algorithms, in particular, the recursive bisection algorithms and GNN architectures definition;
  \item \py{io}: functions needed to load and write meshes together with the graph extraction process, plus additional post-processing functions for visualization.
  \item \py{generate}: utility functions for generating meshes for training GNNs and testing.
\end{itemize}
The package additionally includes the \py{lymphcomm.py} script, which can be called from the command line to agglomerate a single mesh, and is also used for communication with the \texttt{lymph} library (see Section~\ref{sec:lymph_interface}). 
Finally, the \py{models} folder contains some pre-trained GNN models according to the specifications detailed in the previous sections.
\subsubsection{The Agglomeration Models}
The core of \maggnn is the \py{aggmodels} Sub-module. It contains the
\py{AgglomerationModel} abstract class, from which all the implemented models inherit, containing the definition of the recursive bisection algorithms. In particular, here, the SAGE-Base, SAGE-Heterogeneous, RL Partitioner, and RL Refiner Neural network architectures and their training algorithms are defined. METIS and k-means are also included as state-of-the-art methods to which compare them to. 
The GNNs are all implemented using Pytorch Geometric \cite{fey2019PytorchGeometric} and defaulting to a GPU backend (if possible) to speed up computations.
  
\subsubsection{The IO Sub-module}
This module contains functions for loading meshes and datasets, extracting their graph, and saving agglomerated meshes to a file.
The package relies on \py{meshio} \cite{meshio2024} for input, so most of its wide assortment of supported formats can be agglomerated; polygonal, tetrahedral, hexahedral, and pyramidal cells are currently supported. As far as output is concerned, \py{meshio} is once again used in the two-dimensional case for flexibility in output format; on the other hand, in the three-dimensional case, we rely on the \py{vtk} library \cite{vtkBook}, since \py{meshio} has some limitations when dealing with polyhedra.
The graph extraction process consists of the computation of the adjacency matrix describing the mesh, together with areas/volumes and centroid coordinates for each mesh element. 
The adjacency matrix is computed starting from the mesh connectivity data using Algorithm \ref{alg:adj_extraction}. The algorithm achieves linear complexity w.r.t.\ the number of cells in the average case by the use of a hash map. 

\begin{algorithm}[t]
\caption{Adjacency matrix computation algorithm}
\label{alg:adj_extraction}
\textbf{Input:} $\texttt{cells}$: List of $N$ cells, where each cell has a list of faces.\\
\textbf{Output:} $\mathbf A$: Adjacency matrix of the mesh in sparse format.
\begin{algorithmic}[1]
    \State $\texttt{face\_to\_cell} \gets \{\}$ \Comment{Initialize to empty dictionary, with default element list}
    \For{$\texttt{cell\_id} \gets 0 \text{ to } N - 1$}
    \For{\texttt{face\_id} \textbf{ in } \texttt{cells[cells\_id].faces}}
        \State \texttt{face\_to\_cell}[\texttt{face\_id}].\texttt{append(cell\_id)}
    \EndFor
    \EndFor
    \State $\mathbf A_{ij} \gets 0, \, \forall \, i, j = 1, \dots, N$ \Comment{Initialize adjacency matrix as sparse empty matrix}   
    \For{$i \gets 0 \text{ to } N - 1$}  \Comment{Loop over the id of the cell}
    \For{\texttt{face\_id} \textbf{ in } \texttt{cells[cells\_id].faces}}
        \For{$j$ \textbf{ in } \texttt{face\_to\_cell[face\_id]}}
        \If{$j \neq i$} \Comment{Exclude the cell itself}
            \State $\mathbf A_{ij} \gets 1$
        \EndIf
        \EndFor
    \EndFor
    \EndFor    
    \State \textbf{return} $\mathbf A$   
\end{algorithmic}
\end{algorithm}

\subsubsection{The Generate Sub-module}
\begin{figure}[t]
  \centering
  \subfloat[Quads]{\includegraphics[width=0.16\textwidth]{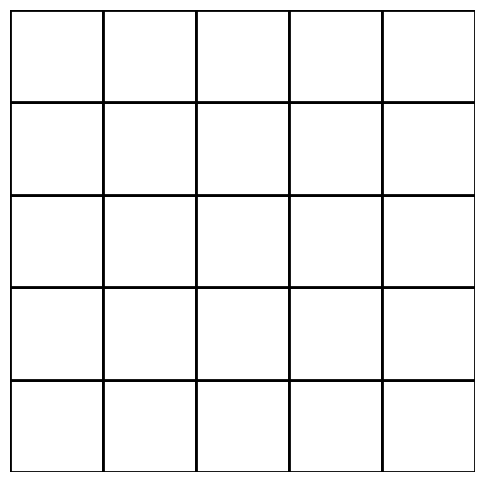}} \quad
  \subfloat[Triangles]{\includegraphics[width=0.16\textwidth]{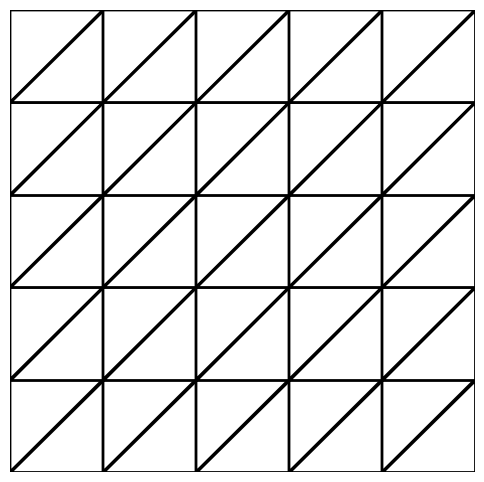}}\quad
  \subfloat[Delaunay]{\includegraphics[width=0.16\textwidth]{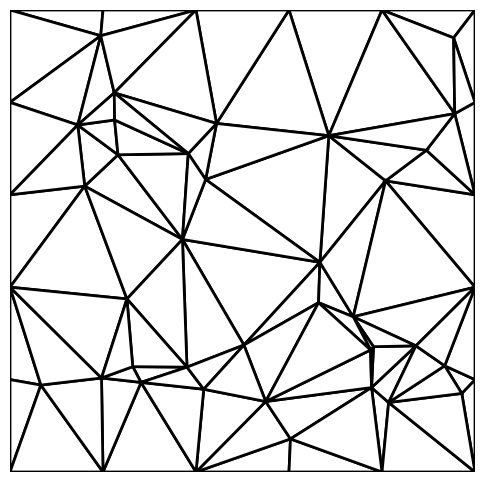}}\quad
  \subfloat[Voronoi]{\includegraphics[width=0.16\textwidth]{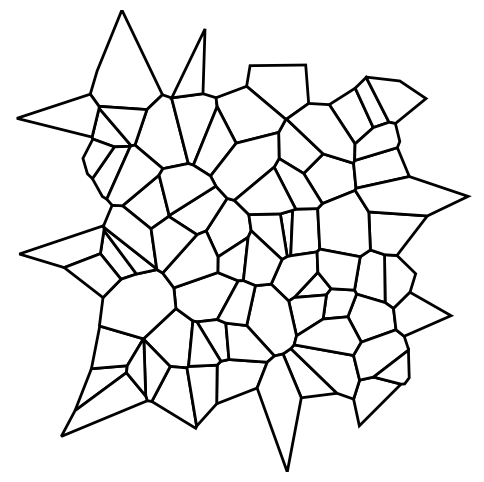}} \quad
  \subfloat[Holes]{\includegraphics[width=0.16\textwidth]{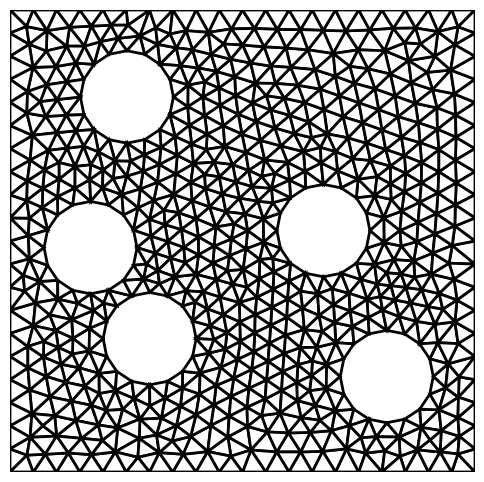}}\\
  \subfloat[Cube]{\includegraphics[trim=300 240 500 300, clip, width=0.25\textwidth]{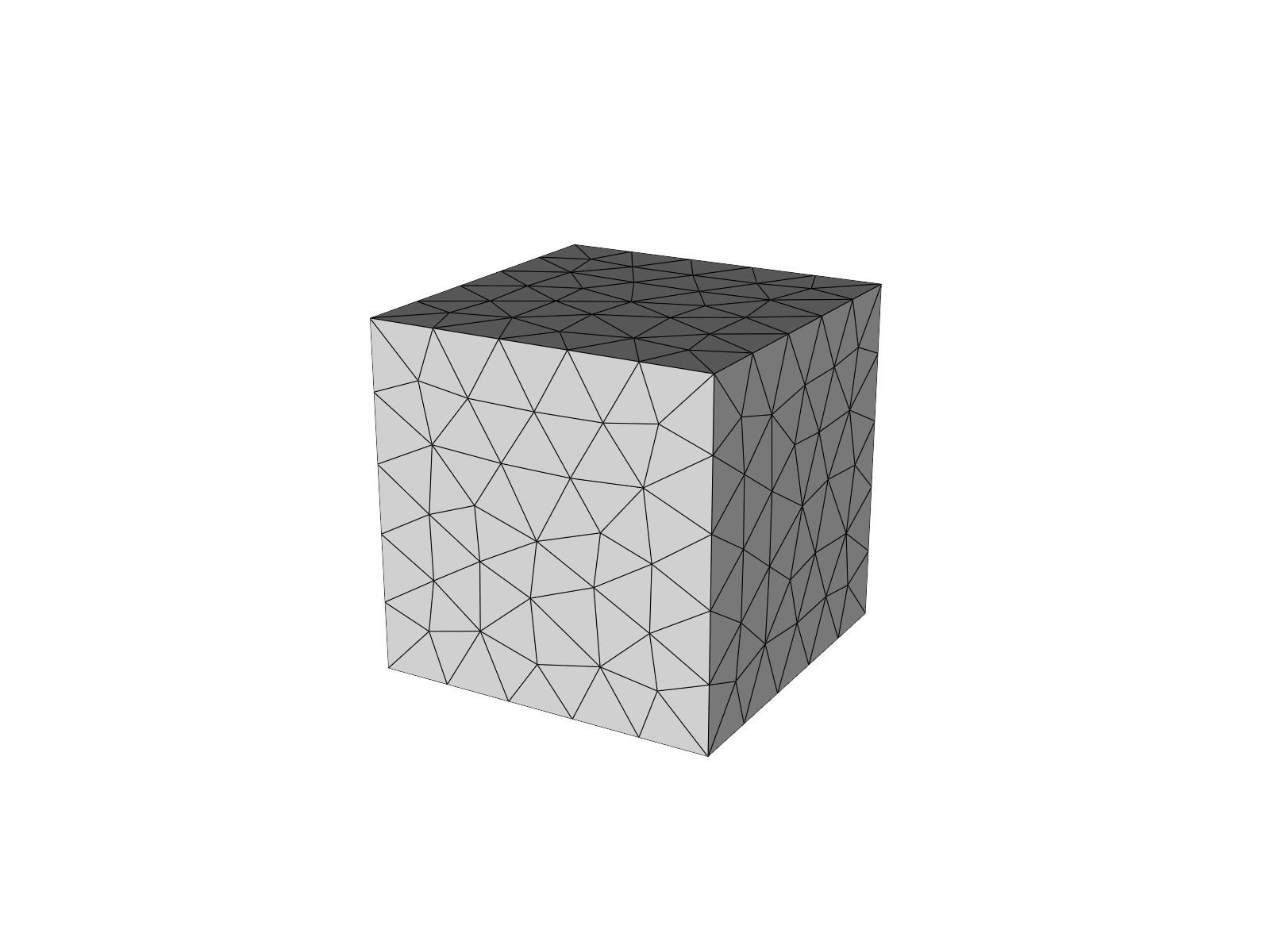}} \quad
  \subfloat[Portion \label{subfig:cube_portion}]{\includegraphics[trim=100 50 30 10, clip, width=0.25\textwidth]{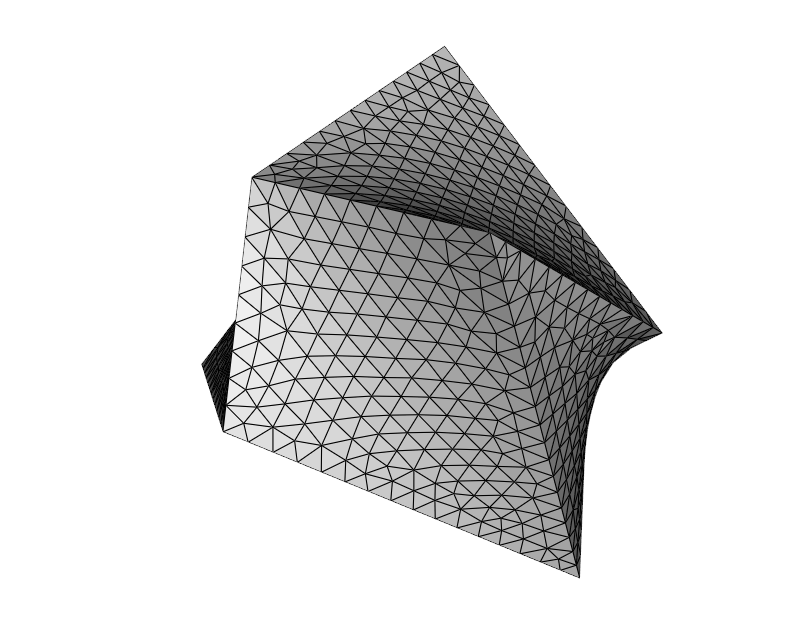}}
  \caption{Examples of the meshes that can be generated using the \py{generate} Sub-module of  \maggnn.}
  \label{fig:generable_meshes}
\end{figure}

The \py{generate} Sub-module enables the generation of simple meshes of the unit square (including structured triangles, structured quadrilaterals, random Delaunay triangulations, random Voronoi tessellations, random circular holes and inclusions) and of the unit cube or portions of it. Some examples of meshes that can be generated through this module are shown in Figure~\ref{fig:generable_meshes}. For mesh generation, \maggnn relies on the open source software \py{Gmsh} \cite{geuzaine2009gmsh}. The main purpose of this module is to supply a convenient way of randomly generating large datasets through the \py{dataset_2D} and \py{dataset_3D} functions for training Neural Networks and testing different agglomeration approaches.
  
\subsection{Agglomeration Algorithms}
\label{sec:core}
At the core of \maggnn is the interplay between agglomeration models and agglomeration modes. Agglomeration models, like SAGE-Base or METIS, are subclasses of the \py{AgglomerationModel} class and define how the mesh is bisected and what data is needed to perform it. Agglomeration modes, on the other hand, dictate how the bisection is applied in practice and act on an algorithmic level. \maggnn provides a few different agglomeration modes:
\begin{itemize}
  \item \emph{direct k-way}: the mesh graph is partitioned in $k$ parts in one shot. This is not available for the GNNs methods, as they rely on recursive bisection.
  \item \emph{number of refinements}: the mesh is bisected recursively $N_{ref} \in \mathbb{N}$ times, leading to a total of $2^{N_{ref}}$ agglomerated elements.
  \item \emph{target size}: the mesh is bisected recursively until all elements have a diameter that is smaller than the given input quantity.
  \item \emph{multiplicative factor}: similar to the previous one, but the target dimension is computed as a fraction of the diameter of the entire mesh, expressed by the input parameter \py{mult_factor}.
  \item \emph{segregated}: exclusively for heterogeneous meshes; the heterogeneous parts of the mesh are separated and their graphs partitioned independently, relying on the \emph{multiplicative factor} algorithm; the node classification is then used to perform the actual agglomeration of the cells in one shot.
  \item \emph{coarsen}: allows to agglomerate (according to one of the previous modes) only part of the mesh by passing the indices of the cells corresponding to it.
  \item  \emph{multilevel recursive bisection}: the mesh graph is initially coarsened using a greedy heavy edge matching algorithm \cite{dhillon2007graclusandmultilevel}, then bisected; the coarse partition is successively projected unto the finer graph and refined (using, for example, the RL Refiner of Section~\ref{sec:rlrefiner}); this process is recursively repeated to obtain the desired partition. See
  \texttt{\href{https://github.com/lymphlib/magnet/blob/main/magnet/_absaggmodels.py}{\_absaggmodels.py}}
  line 395 for details.
\end{itemize}
While methods like METIS and k-means generally perform better when partitioning the mesh directly in the desired number of parts with \emph{k-way} mode instead of applying them recursively, it is still possible to do so using one of the other modes. This is important because recursive bisection algorithms naturally generate a hierarchy of nested coarser grids that can be used in multigrid methods to accelerate numerical solvers.
We note that for really small \py{mult_factor}, the necessary number of bisections and corresponding recursive calls in the stack grows rapidly; this can quickly use up a lot of memory, so the implementation is sequential rather than recursive. A schematic of the \emph{multiplicative factor} bisection algorithm is reported in Algorithm \ref{alg:multfactor}.
\begin{algorithm}[!t]
  \caption{Sequential \emph{multiplicative factor} algorithm}
  \label{alg:multfactor}
  \textbf{Input:} $\mathcal{G}, \texttt{mult\_factor}$: graph to be partitioned, multiplicative factor\\
  \textbf{Output:} Elements classified based on their agglomerated element
  \begin{algorithmic}[1]
    \State $\hat{h}$ $\gets$ $\text{diam}(\mathcal{G}) \, \texttt{mult\_factor}$ \Comment{ $\text{diam}(\mathcal{G})$ is the largest pair distance of vertexes of $\tau_h$.}
    \State $\texttt{parts} \gets [\mathcal{G}]$
    \State $\texttt{output} \gets []$
    \While{$\texttt{parts}\ \textbf{is not}\ \emptyset$}
      \State{$\texttt{new\_parts} \gets []$}
      \For{$\mathcal{S} \ \textbf{in}\ \texttt{parts}$}
        \If{diam{$(\mathcal{S})$} $>$ $\hat{h}$} 
          \State{$\mathcal{S}_1, \mathcal{S}_2 \gets$ \Call{bisect}{$\mathcal{S}$}} \Comment{If the element is too ``large", bisect it}
          \State{$\texttt{new\_parts} \gets \texttt{new\_parts} + [\mathcal{S}_1, \mathcal{S}_2]$}
        \Else
          \State{\texttt{output.append}$(\mathcal{S})$} \Comment{Otherwise, insert it in the final output}
        \EndIf
      \EndFor
      \State{$\texttt{parts} \gets \texttt{new\_parts}$}
    \EndWhile
    \State \Return{$\texttt{output}$}
  \end{algorithmic}
\end{algorithm}
We remark that some agglomeration approaches, especially k-means and SAGE-Base, tend to agglomerate parts of the mesh that are not connected, in particular when holes are present or the underlying geometry is very complex; as such, at the end of the graph partitioning process an additional check on the connectedness of each subgraph is performed, and if more than one connected component is found they are separated (an efficient algorithm for this task is described in \cite{Pearce2005scipyconectcomponents}).
Finally, since we always frame mesh agglomeration as a graph partitioning problem, the actual geometric merging of the cells can be performed later after the nodes have been classified; this modularity allows for great flexibility in the classification process while sharing the same implementation of the merging step across all agglomeration models.

\subsection{Extensibility}
The agglomeration model class hierarchy has been designed with extensibility, making adding new bisection strategies easy. When defining a new agglomeration model inheriting from the abstract base class, only two new methods strictly need to be defined: the \py{get_graph} method that specifies the graph data needed for bisection, and the \py{bisect_graph} method that performs bisection on it; GNN models will additionally need an initialization method defining their architecture and a \py{forward} method. Of course, further customization is always possible by overriding existing methods or implementing new ones.
  
\subsection{Heterogeneous Mesh Agglomeration}
\label{sec:heterogeneous_mesh}
Traditional mesh agglomeration strategies, like METIS, cannot deal with meshes describing heterogeneous domains, i.e.\ meshes that have regions with different physical parameters. The only way to do this using these methods is to agglomerate different regions independently and then merge them; however, the merging process can be rather expensive. Also, this strategy is not well suited for domains with micro-structures, where independent agglomeration of the structures leads to poorer quality agglomerated elements. \maggnn can agglomerate heterogeneous meshes in two different ways:
\begin{enumerate}
  \item By using a GNN (SAGE-Heterogeneous) that takes as input an additional node feature, the physical group, that describes the heterogeneity of the mesh, by using the strategy illustrated in Section~\ref{sec:sagehetero}.
  \item By first separating the heterogeneous parts of the domain and partitioning each corresponding graph separately, using \emph{segregated} mode, described in Section~\ref{sec:core}.
\end{enumerate}
The first approach is less expensive, but can only handle a certain class of physical groups and could wrongly agglomerate elements with different physical parameters. The second approach always separates correctly physical groups and can exploit any of the supported agglomeration strategies, but is more expensive, especially in the case of many small inclusions in the domain.
\subsection{Quality Metrics}
\label{sec:quality_metrics}
The notion of mesh quality for triangular, quadrilateral, and tetrahedral meshes is deeply explored in the literature \cite{shewchuk2002triangulartightqulity}. However, there is no shared consensus on the characterization of good quality for polygonal and polyhedral meshes, which is still an object of ongoing research \cite{sorgente2021polyhedralmeshqualityindicator, zunic2004newconvexitymeasure}.  
For example, a standard assumption on shape regularity for FEMs on polytopal meshes is that each element of the mesh is star-shaped with respect to a point; however, numerical experiments have shown that these methods can reliably solve differential problems also on much more irregular meshes \cite{attene2019benchmarkpolygonqualitymetrics}, so the practical requirements are much weaker. 
\maggnn is flexible in this regard and allows for easy implementation of new quality metrics and integrates them with the rest of the framework. The following quality metrics described in \cite{antonietti2022refinementpolygridsCNN, attene2019benchmarkpolygonqualitymetrics}, which are defined element-wise and take values in [0,1], are provided in \maggnn:
\begin{itemize}
  \item Circle Ratio (CR): the ratio between the radius of the smallest sphere inscribed in the element $P$ and the biggest sphere containing it:
  \begin{equation}
    \label{eq:CR}
    \text{CR}(P)=\frac{\displaystyle\max_{B(r) \subset P} r}{\displaystyle\min_{B(r) \supset P} r}
  \end{equation}
  where $B(r)$ is the ball of radius $r$. It is a measure of the roundness of the elements: the closer the circle ratio is to 1, the closer the element is to a sphere. 
  \item Area to Perimeter Ratio (APR): the ratio between the area of the polygon $P$ and its perimeter.
  \begin{equation}
    \label{eq:APR}
    \text{APR}(P)=\frac{4\pi |P|}{|\partial P|^2}
  \end{equation}
  It is also known as the iso-perimetric quotient because it is the ratio between the area of $P$ and that of a circle with the same perimeter $P$. This metric is a measure of the compactness of the polygon: it has a maximum value of $1$ for the circle and decreases for less compact polygons (i.e.\ the shape is more spread out).
  An equivalent metric in three-dimensions is sphericity:
  \begin{equation}
    \label{eq:sphericity}
    \Psi = \frac{\sqrt[3]{36 \pi |P|^2}}{|\partial P|}
  \end{equation}

  \item Uniformity Factor (UF): the ratio between the diameter of the element and the mesh size:
  \begin{equation}
    \label{eq:UF}
    \text{UF}(P)=\frac{\text{diam}(P)}{h}, \quad h = \max_{P \in \tau_h}\text{diam}(P),
  \end{equation}
  where $\text{diam}(P)$ indicates the maximum among the distances of any pair of vertices of $P$. Values closer to one denote that mesh elements have similar diameters.
  
  \item Volumes Difference (VD): the relative difference between the volume of the cell $P$ and the volume that each cell should have if they all had the same volume $\hat{V}$:
  \begin{equation}
    \text{VD}(P)=\frac{|\text{Vol}(P) - \hat{V}|}{\hat{V}}, \quad \hat{V}=\frac{\sum_P \text{Vol}(P)}{N}
  \end{equation}
  This quantity can take any positive value, so to have a metric between zero and one we consider instead:
  \begin{equation}
    \label{eq:VD}
    \widetilde{\text{VD}}(P)=\frac{1}{1+\text{VD}(P)}
  \end{equation}
  A value of one denotes that all elements have the same volume, while lower values correspond to a less homogeneous volume distribution across the cells.

  \item Heterogeneity Preservation (HP): the relative presence of different physical tags $p_i = 0, 1$ in an agglomerated cell $P$ of $n_P$ fine cells.
  \begin{equation}
    \label{eq:HP}
    \text{HP}(P)= \max\{\bar{p}, 1 - \bar{p}\}, \quad \bar{p} = \frac{1}{n_P}\sum_i^{n_P} p_i
  \end{equation}
  A value of one denotes that all elements have the same physical tag, while lower values correspond to a less homogeneous distribution across the cells.
  
\end{itemize}
  
  \section{Numerical Results}
  \label{sec:examples}
  In this section, we provide a brief step-by-step guide, showing how to train a GNN and how to agglomerate meshes with \maggnn, starting from the two-dimensional homogeneous case and moving on to more complex applications. We also use these examples to comment on the performance and main features of the considered methods, comparing the GNN strategies with METIS and k-means using the metrics of Section~\ref{sec:quality_metrics}.
  
  \subsection{GNN Training}
  \label{sec:training_guide}
\begin{figure}[!t]
    \centering
    \subfloat[2D SAGE-Base]{\includegraphics[trim=60 0 60 0, clip, width=0.5\textwidth]{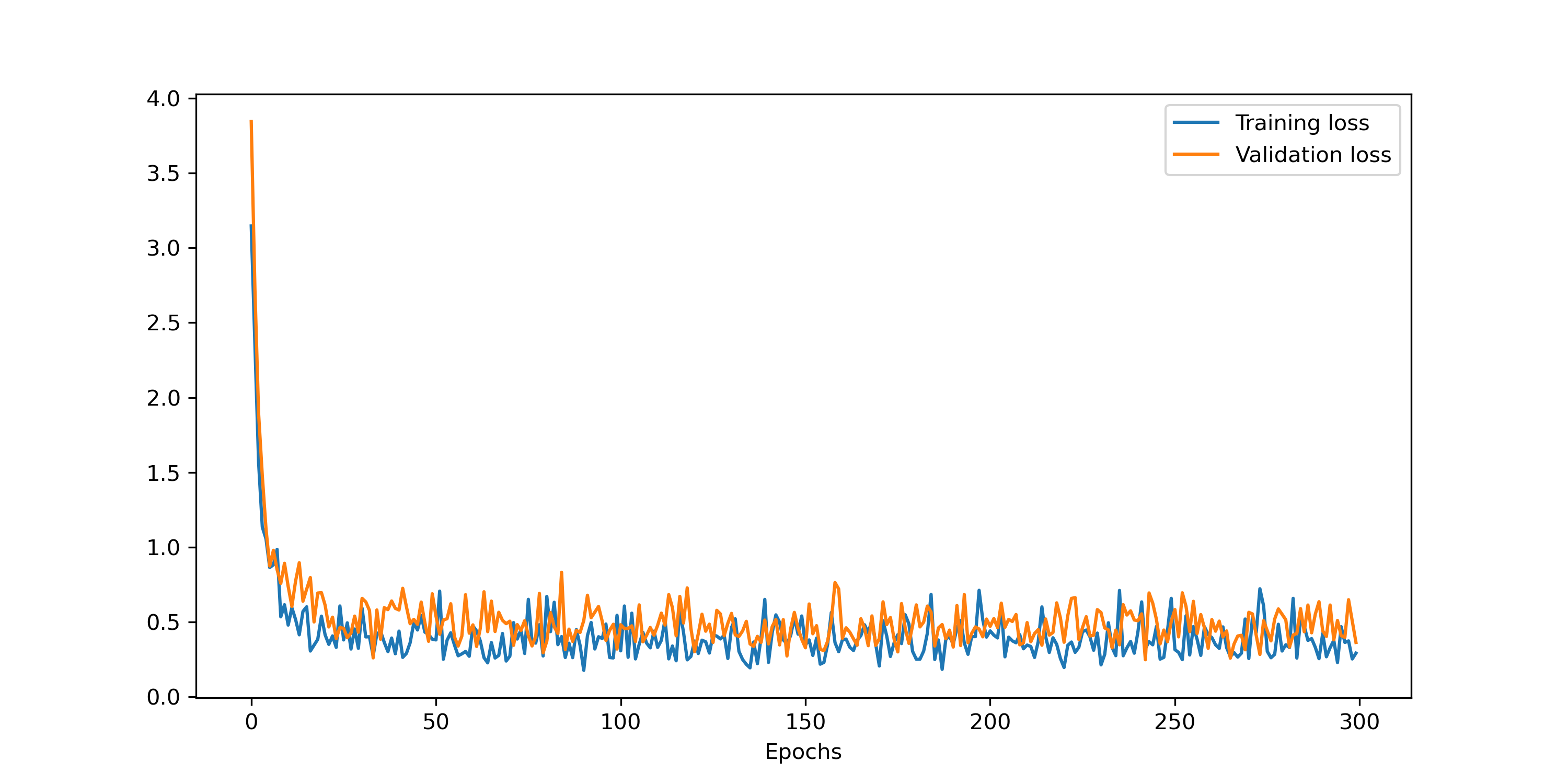}}
    \subfloat[3D SAGE-Base]{\includegraphics[trim=60 0 60 0, clip, width=0.5\textwidth]{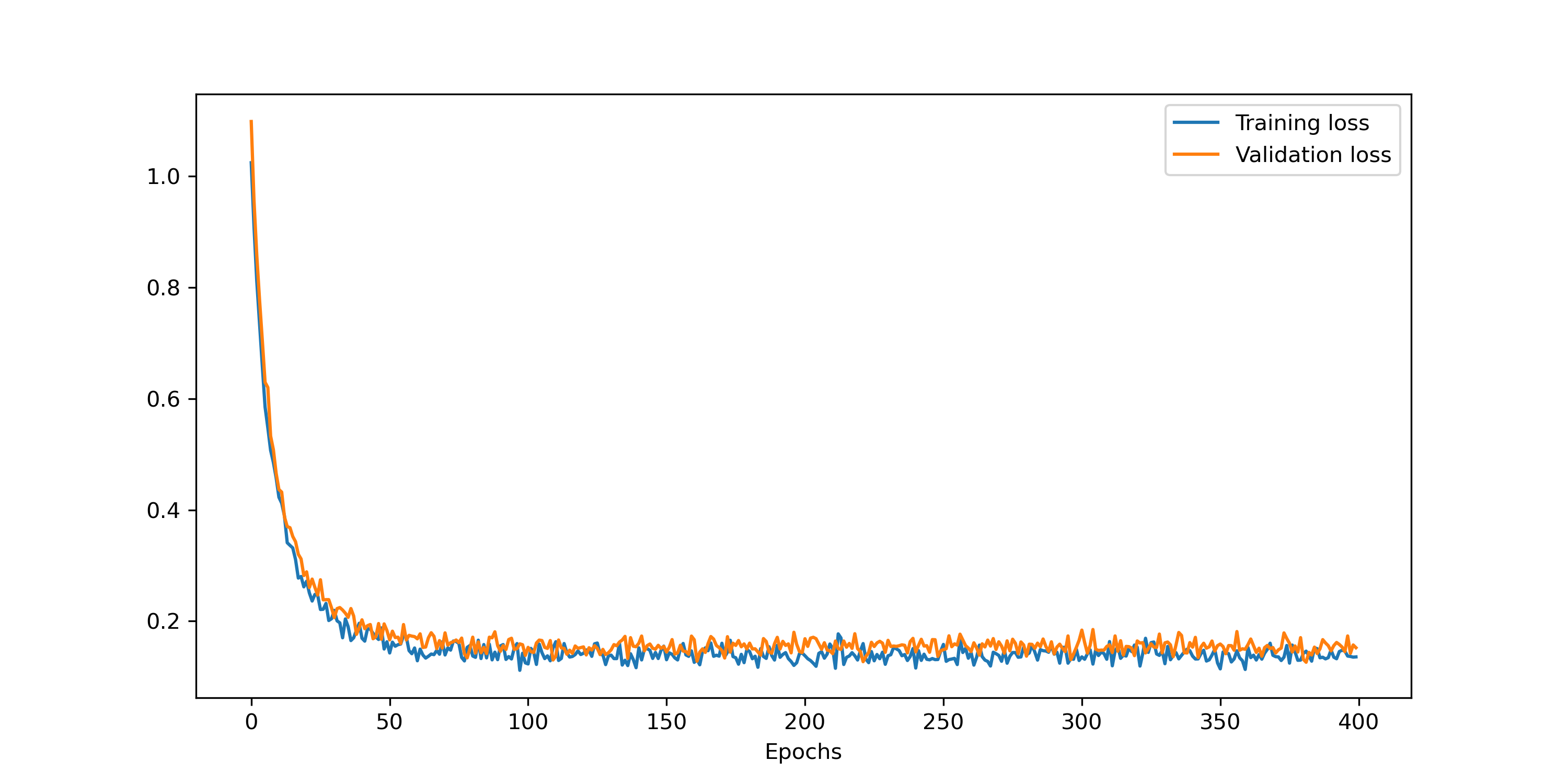}}
    \caption{Training and validation losses history plot of the two-dimensional and three-dimensional versions of the SAGE-Base model of Section~\ref{sec:sagebase}.}
    \label{fig:loss_history}
  \end{figure}

  To train a GNN, we need a training dataset and a validation dataset to evaluate the model's performance during the training. To speed up training, mesh graphs (adjacency matrix, centroid coordinates, volumes) need to be pre-computed and stored. \maggnn provides a \py{generate} Sub-module that generates meshes while computing their graph data; it is also possible to create a dataset from a folder of existing meshes by using \py{create_dataset}. After instantiating the GNN model, training can be initiated with the \py{train_GNN} method; Listing \ref{lst:gnntrain} gives an example of how the whole pipeline looks in \maggnn.
  \begin{lstlisting}[language=Python, caption= Example of the code for GNN training., label={lst:gnntrain}]
    import magnet
    # create a training dataset and load it
    magnet.generate.dataset_2D({'random_delaunauy':200},'datasets/trainig_set')  
    training_dataset = magnet.io.load_dataset('datasets/trainig_set')
    # initialize GNN
    NN = magnet.aggmodels.SageBase(64,32,3,2).to(magnet.aggmodels.DEVICE)  
    NN.train_GNN(training_dataset, epochs=200, batch=4, lr=0.0001)# training loop
    NN.save_model('models/test_training.pt')
  \end{lstlisting}  
  For the reinforcement learning approach, the method \py{A2C_train} is used instead.
  When training is completed, you have the option to save a training summary log file and the loss history; Figure~\ref{fig:loss_history} shows as an example the loss history plots for the SAGE-Base models of Section~\ref{sec:sagebase}. 

  \subsection{Computational Cost}
\begin{figure}[!t]
    \centering
    \includegraphics[width=1.0\textwidth]{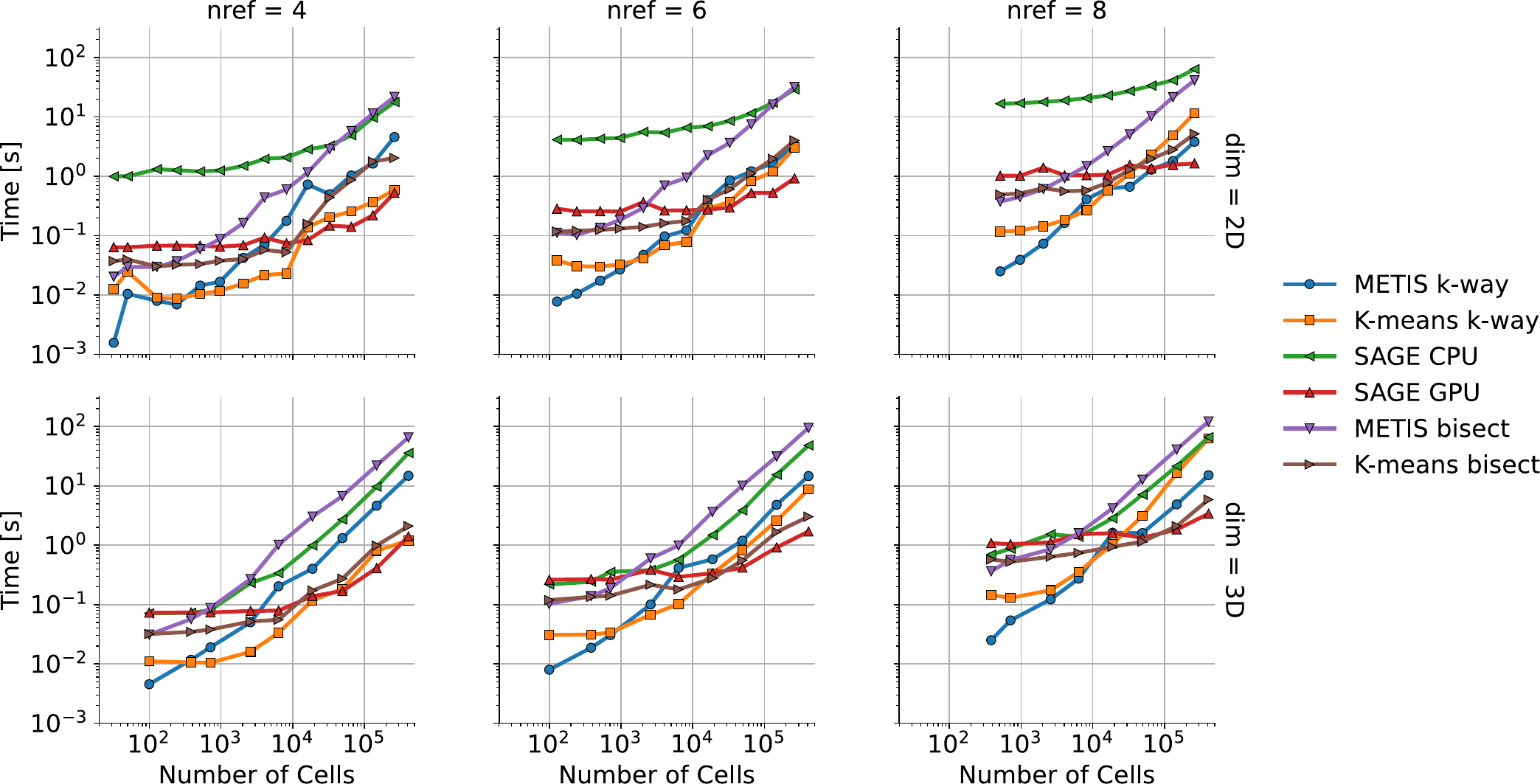} 
    \caption{Comparison of the computational cost of the algorithms present in \maggnn for the agglomeration of 2D meshes (top) and 3D meshes (bottom).}
    \label{fig:comp-cost}
  \end{figure}

In this section, we explore the advantages of the GNN-based agglomeration from a point of view of computational cost. In particular, we provide quantitative measurements of the wall-clock time for METIS, k-means, and GNN-based agglomeration. For the GNN, we compare the wall-clock time for both a CPU-based and a GPU-based run. In order to provide a comparison that is as fair as possible, we use the resources offered for free in Google Colab, namely an Intel(R) Xeon(R) CPU @ 2.20GHz and an Nvidia T4 GPU. In Figure~\ref{fig:comp-cost}, we provide an in-depth analysis of the computational cost depending on the space dimension (2D/3D) and the number of refinements. If we compare the k-way and the bisect variants, we see that the k-way algorithm is much cheaper than the bisect version of METIS. This holds in both 2D and 3D since METIS considers only the dual graph and is thus space independent. On the other hand, k-means scales linearly with respect to the space dimension, the number of cells, and the number of clusters. Hence, the k-way version of k-means performs worse than the bisect version as we refine the mesh (namely, the number of cells and the number of refinements increase). This is particularly evident in the case \py{nref=8} where we can see that the k-way k-means has a much steeper slope compared to the other methods. When running the GNN method on a CPU, the cost is asymptotically comparable to the bisect version of METIS. On the other hand, the GPU version, while having a larger fixed cost until around $10^4$ cells, is the cheaper method for larger meshes.
These results strongly suggest that the GPU hardware is essential to fully unlock the scalability benefits of GNN-based agglomeration, making it a cheaper alternative to traditional methods for medium to large-scale problems.
 
  
  \subsection{Test Case 1: Two-Dimensional Brain Slice}
  \label{sec:brain_slices}
\begin{figure}[!htbp]
    \centering    
    \subfloat[Original mesh]{\includegraphics[width=0.3\textwidth]{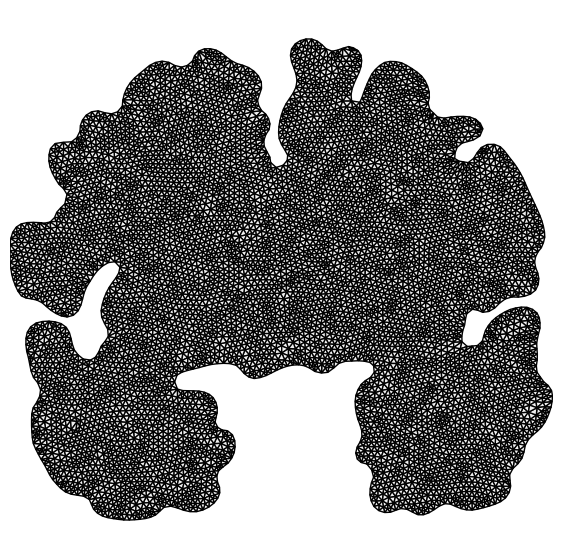}} \quad
    \subfloat[METIS]{\includegraphics[width=0.3\textwidth]{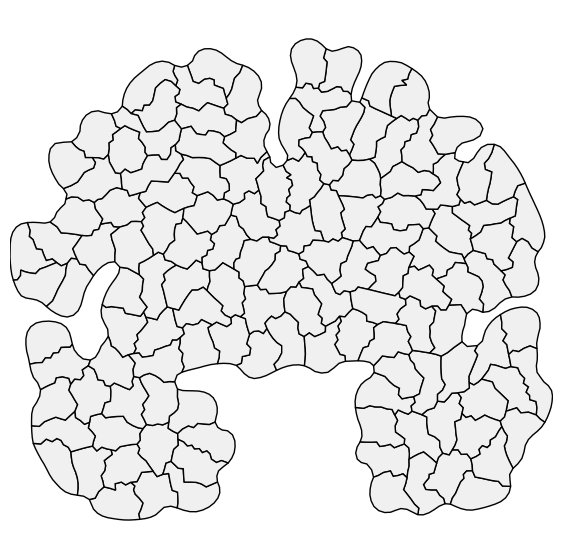}} \quad
    \subfloat[k-means]{\includegraphics[width=0.3\textwidth]{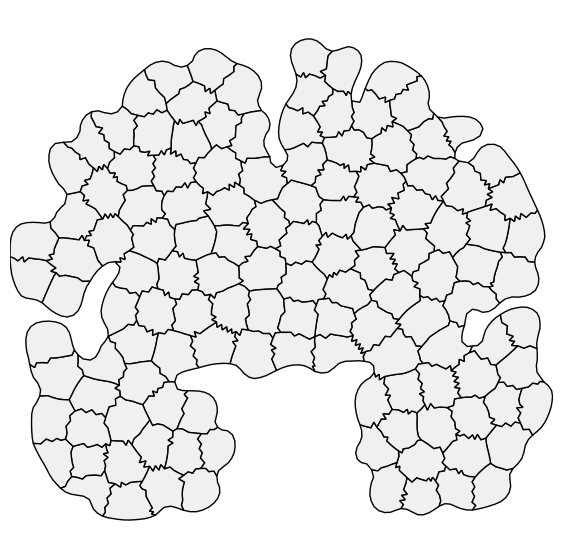}} \\
    \subfloat[SAGE-Base]{\includegraphics[width=0.3\textwidth]{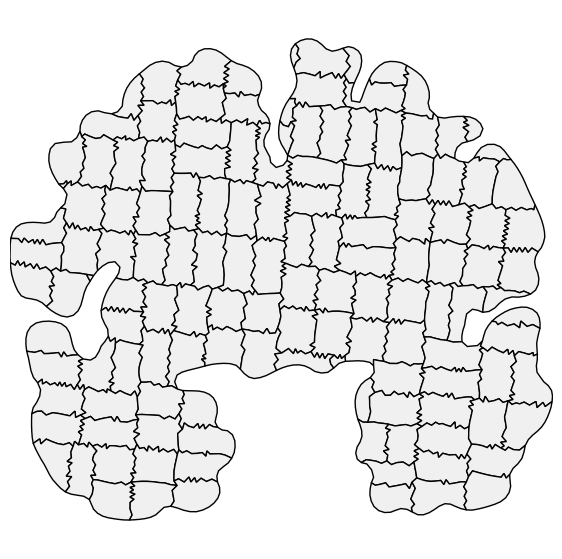}} \quad
    \subfloat[SAGE-Base + RL Refiner]{\includegraphics[width=0.3\textwidth]{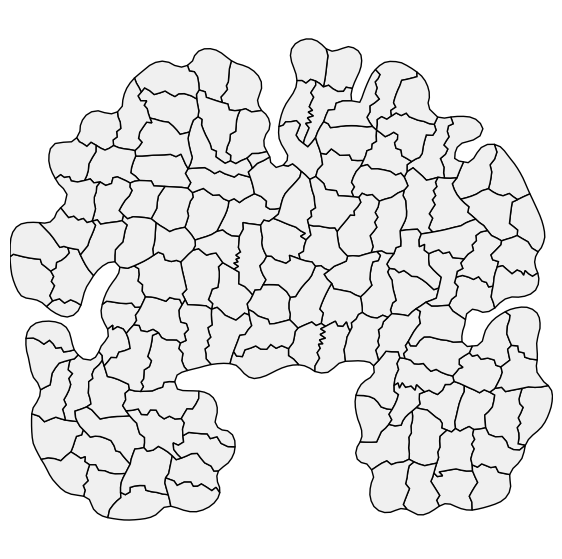}} \quad
    \subfloat[RL Partitioner + RL Refiner]{\includegraphics[width=0.3\textwidth]{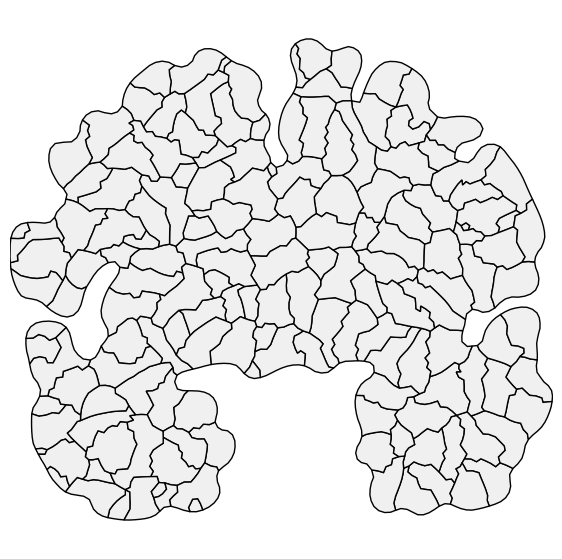}}\\
    \subfloat[Quality metrics box plots \label{subfig:brainslice_qual}]{\includegraphics[width=1\textwidth]{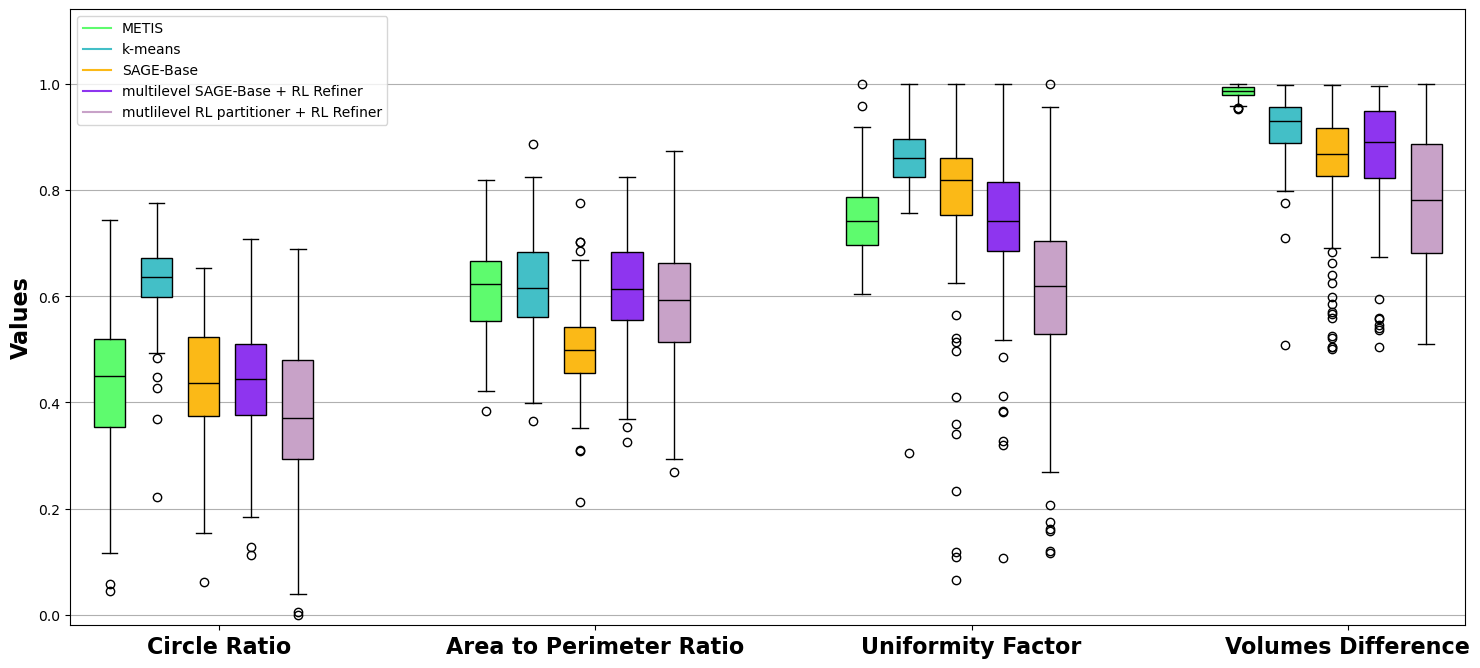}}
    \caption
    {Test Case 1.a: unstructured mesh of a human brain section consisting of 14372 triangles, agglomerated using different methods (METIS, k-means, SAGE-Base, SAGE-Base and RL coarse partitioner in multilevel framework with RL Refiner), together with the box plots of the computed quality metrics (CR, APR, UF and $\widetilde{\text{VD}}$, defined as in Eq.~\eqref{eq:CR}-\eqref{eq:VD}).}
    \label{fig:brain_slice_agg}
  \end{figure}

\begin{figure}[!htbp]
    \centering  
    \subfloat[Original mesh]{\includegraphics[width=0.3\textwidth]{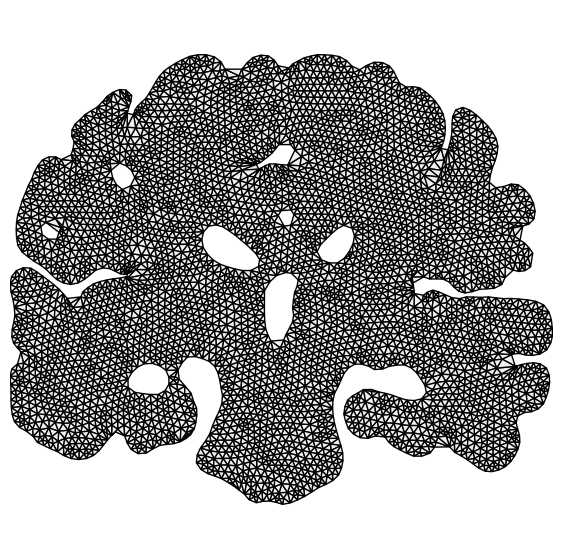}} \quad
    \subfloat[METIS]{\includegraphics[width=0.3\textwidth]{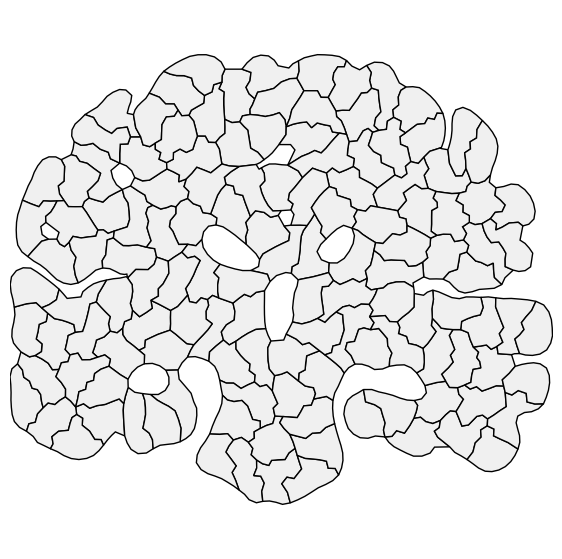}} \quad
    \subfloat[k-means]{\includegraphics[width=0.3\textwidth]{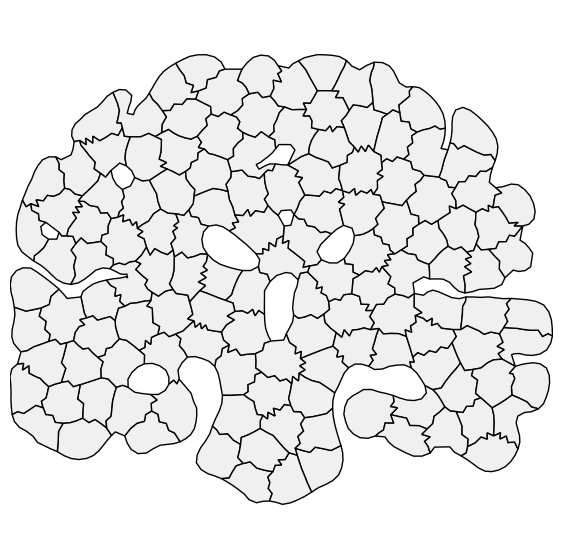}} \\
    \subfloat[SAGE-Base]{\includegraphics[width=0.3\textwidth]{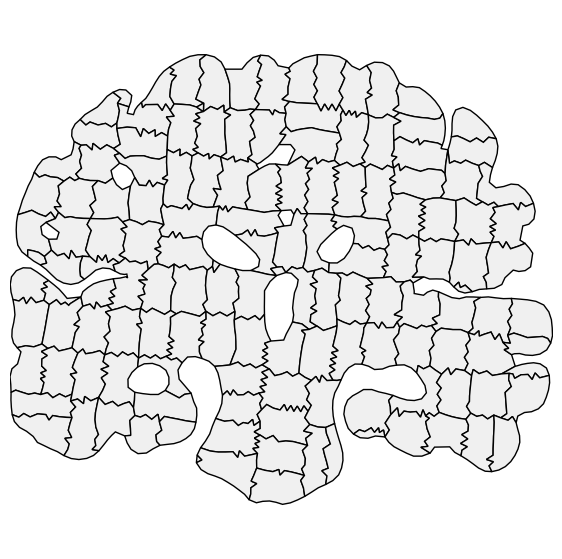}} \quad
    \subfloat[SAGE-Base + RL Refiner]{\includegraphics[width=0.3\textwidth]{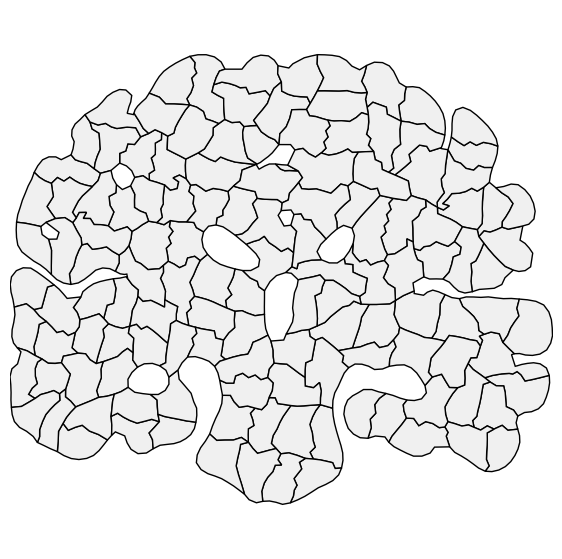}} \quad
    \subfloat[RL Partitioner + RL Refiner]{\includegraphics[width=0.3\textwidth]{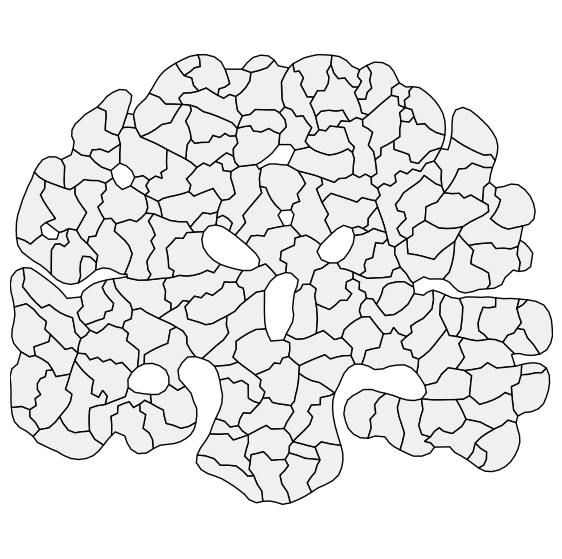}}\\  
    \subfloat[Quality metrics box plots \label{subfig:holesbrain_qual}]{\includegraphics[width=1\textwidth]{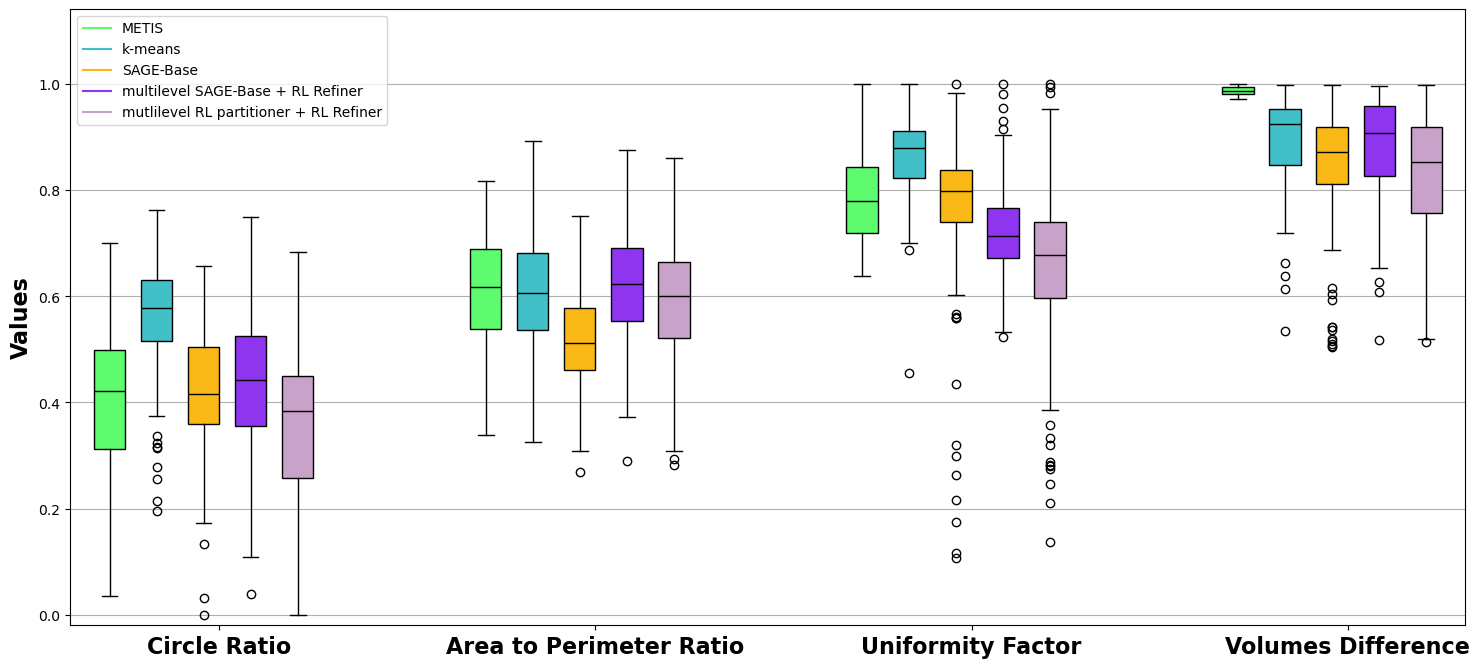}}   
    \caption{Test Case 1.b: unstructured mesh of a human brain section consisting of 8597 triangles, agglomerated using different methods (METIS, k-means, SAGE-Base, SAGE-Base and RL coarse partitioner in multilevel framework with RL Refiner), together with the box plots of the computed quality metrics (CR, APR, UF and $\widetilde{\text{VD}}$, defined as in Eq.~\eqref{eq:CR}-\eqref{eq:VD}).}
    \label{fig:holesbrain_agg}
  \end{figure}

  To showcase the agglomeration of two-dimensional homogeneous meshes and the generalization capabilities of the presented GNN models, we consider as test cases two human brain sections coming from Magnetic Resonance Imaging (MRI); these domains are very complex, have narrow sections, and also include holes in the second case. 
  The code for Test Case 1 is reported in Listing \ref{lst:hom_agg_2d}.
  First, we load the mesh by using \py{io.load_mesh}, which also extracts the mesh graph needed for agglomeration. Then, we initialize our GNN agglomeration model and load it from a state dictionary. We can then proceed to agglomerate the mesh using one of the modes introduced in Section~\ref{sec:core} and plot it. Finally, using the \py{get_quality_metrics} method, the quality metrics described in Section~\ref{sec:quality_metrics} are computed and plotted. 
  \begin{lstlisting}[language=Python, caption=Example of the code for agglomerating one mesh., label={lst:hom_agg_2d}]
    from magnet import io, aggmodels
    brain_mesh = io.load_mesh('datasets/BrainCoronal.vtu')
    NN = aggmodels.SageBase(64,32,3,2).to(aggmodels.DEVICE)  # initialize GNN
    agglomerated_brain = NN.agglomerate(brain_mesh, mode='Nref', nref=7)
    agglomerated_brain.view(colors='grey')
    agglomerated_brain.get_quality_metrics(boxplot=True)
  \end{lstlisting}
We agglomerate the two meshes using five different models (METIS, k-means, SAGE-Base,  multilevel approach with RL Refiner, and either SAGE-Base or RL Partitioner as coarse partitioner) using the template of Listing \ref{lst:hom_agg_2d}. For the GNN models, \emph{number of refinements} mode was used with \py{nref=7}, creating a total of 128 elements, while METIS and k-means were used with \emph{direct k-way} mode and \py{k=128} to have a fair comparison. We report the results in Figure~\ref{fig:brain_slice_agg} and Figure~\ref{fig:holesbrain_agg}. In particular, Figure~\ref{subfig:brainslice_qual} and Figure~\ref{subfig:holesbrain_qual} show the quality metrics computed for each agglomerated mesh, denoting meshes obtained with different agglomeration models with distinct colors; since these metrics are defined element-wise, each box plot shows their distribution in a single mesh.
  We observe that k-means is by far the best one in terms of both circle ratio and uniformity factor, which is not surprising considering that it is a clustering algorithm that exploits the geometric information of the mesh. METIS performs extremely well in terms of volume uniformity due to its strict algorithmic constraint on partition balance and the fact that we are using cell volumes themselves as node weights. The other models perform roughly the same in terms of the other metrics, except the RL coarse partitioner model, which is slightly worse across the board. Finally, we notice that SAGE-Base tends to bisect the mesh along straight lines, leading to the formation of elements with squared corners and a lower area-to-perimeter ratio as a consequence.
  
  \subsection{Test Case 2: Two-Dimensional Domain with Inclusions}
\begin{figure}[!t]
    \centering
    \includegraphics[width=1\textwidth]{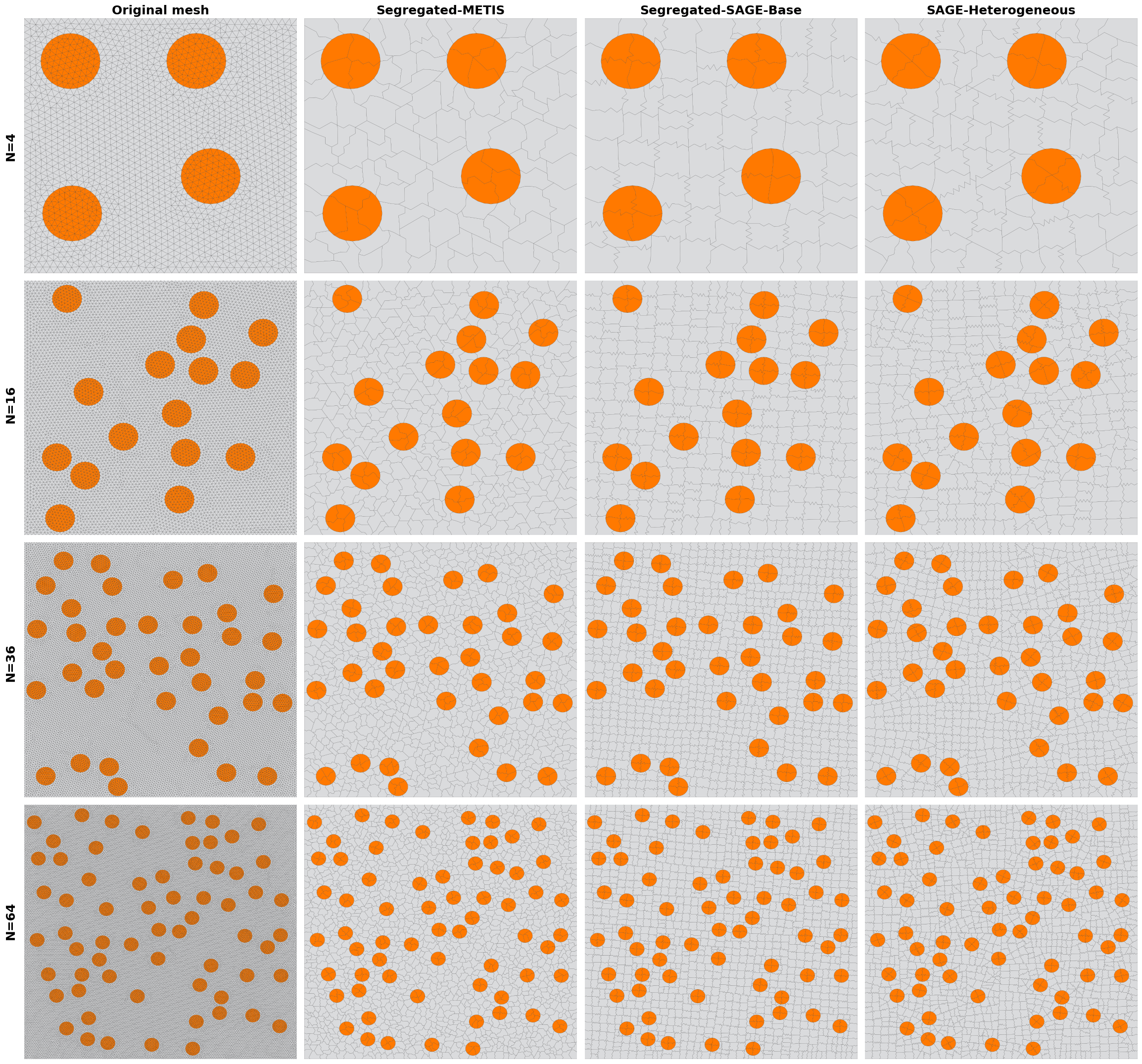}
    \caption{Test Case 2: four meshes of the unit square containing 4, 16, 36, 64 circular inclusions representing microstructures in the domain, consisting of 3796, 15374, 35132, 62766 triangles respectively. The meshes have been agglomerated with METIS and SAGE-Base in \emph{segregated} mode and with SAGE-Heterogeneous,  using a target size equal to 3/4 of the circle's diameter.}
    \label{fig:inclusions_grid_plot}
  \end{figure}

  We consider a set of four triangular meshes of the unit square with an increasing number of circular inclusions that represent microstructures in the underlying domain, with respectively 4, 16, 36, and 64 inclusions each. The radii of the circles have been chosen so that they always cover 15\% of the total area. We agglomerate them both with the SAGE-Heterogeneous model and by using \emph{segregated} mode with both METIS and SAGE-Base using a target size equal to 3/4 of the circle diameter, so that the agglomerated elements have a size comparable to that of the inclusions. The results are reported in Figure~\ref{fig:inclusions_grid_plot}. In this case, the SAGE-Heterogeneous model correctly separates the inclusions; however, since from the second bisection onward only one of the two physical groups appears, k-means is called instead.

\begin{figure}[!t]
    \centering
    \includegraphics[width=1\textwidth]{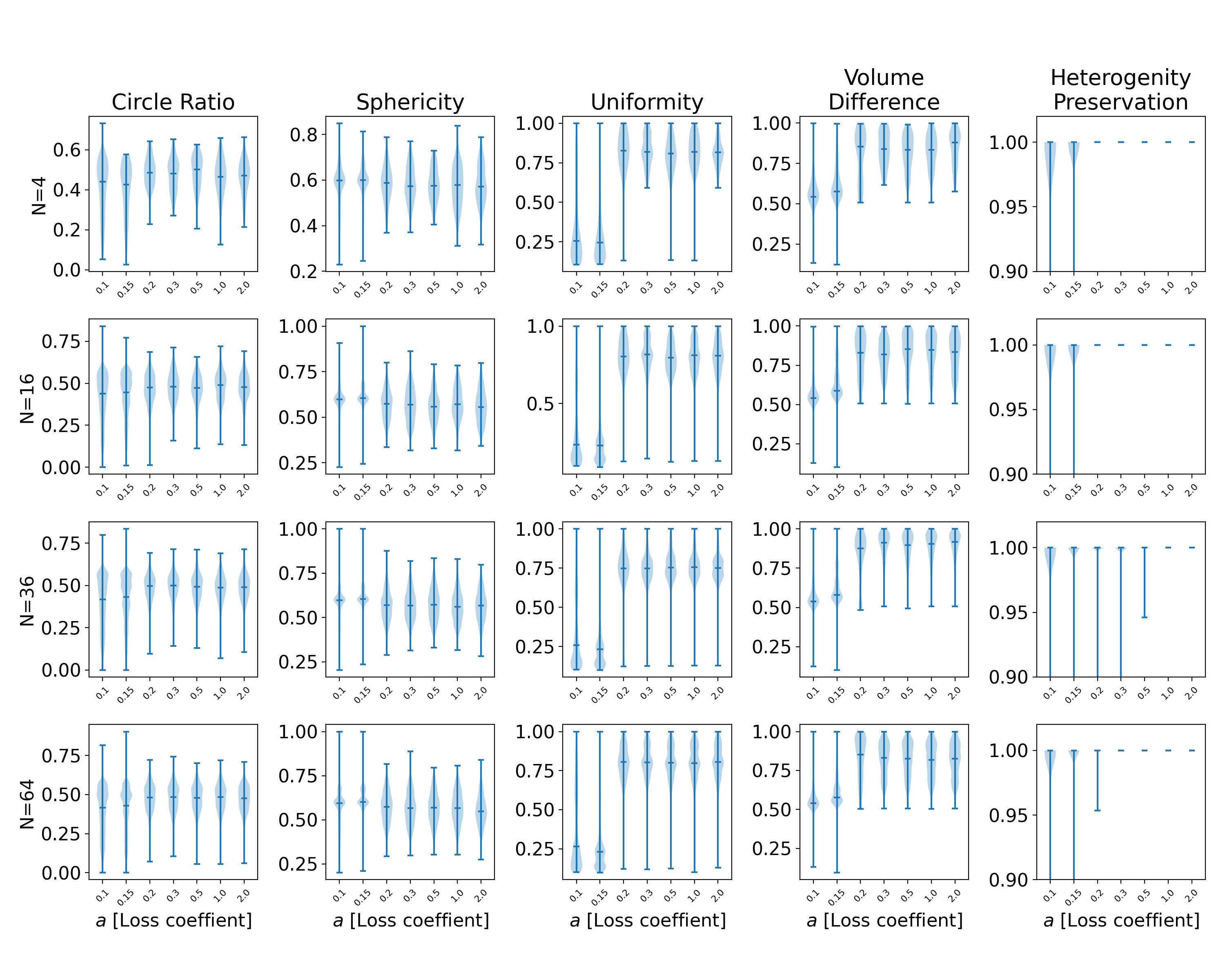}
    \caption{Test Case 2: sensitivity analysis for the parameter $a$ of the loss function Eq.~\eqref{eq:sagehetero-loss}. On each row, the model is tested on a different mesh containing 4, 16, 36, 64 circular inclusions. The four test meshes are the ones considered also in Figure~\ref{fig:inclusions_grid_plot}.}
    \label{fig:sensitivity_analysis_a}
  \end{figure}

  A sensitivity analysis on the parameter $a$ of Eq.~(\ref{eq:sagehetero-loss}) is provided in Figure~\ref{fig:sensitivity_analysis_a}. It shows that most metrics (Sphericity, Uniformity) remain largely unaffected, while Circle Ratio and Volume Difference display higher variability. Heterogeneity Preservation is highly sensitive, with values below 0.95 indicating poor performance. In particular, a value of $a > 0.5$ seems to be the minimum to guarantee consistent results in heterogeneous agglomeration.
  \begin{figure}[t]
    \centering
    \subfloat[Original mesh]{\includegraphics[trim=0 300 0 300, clip, width=0.5\textwidth]{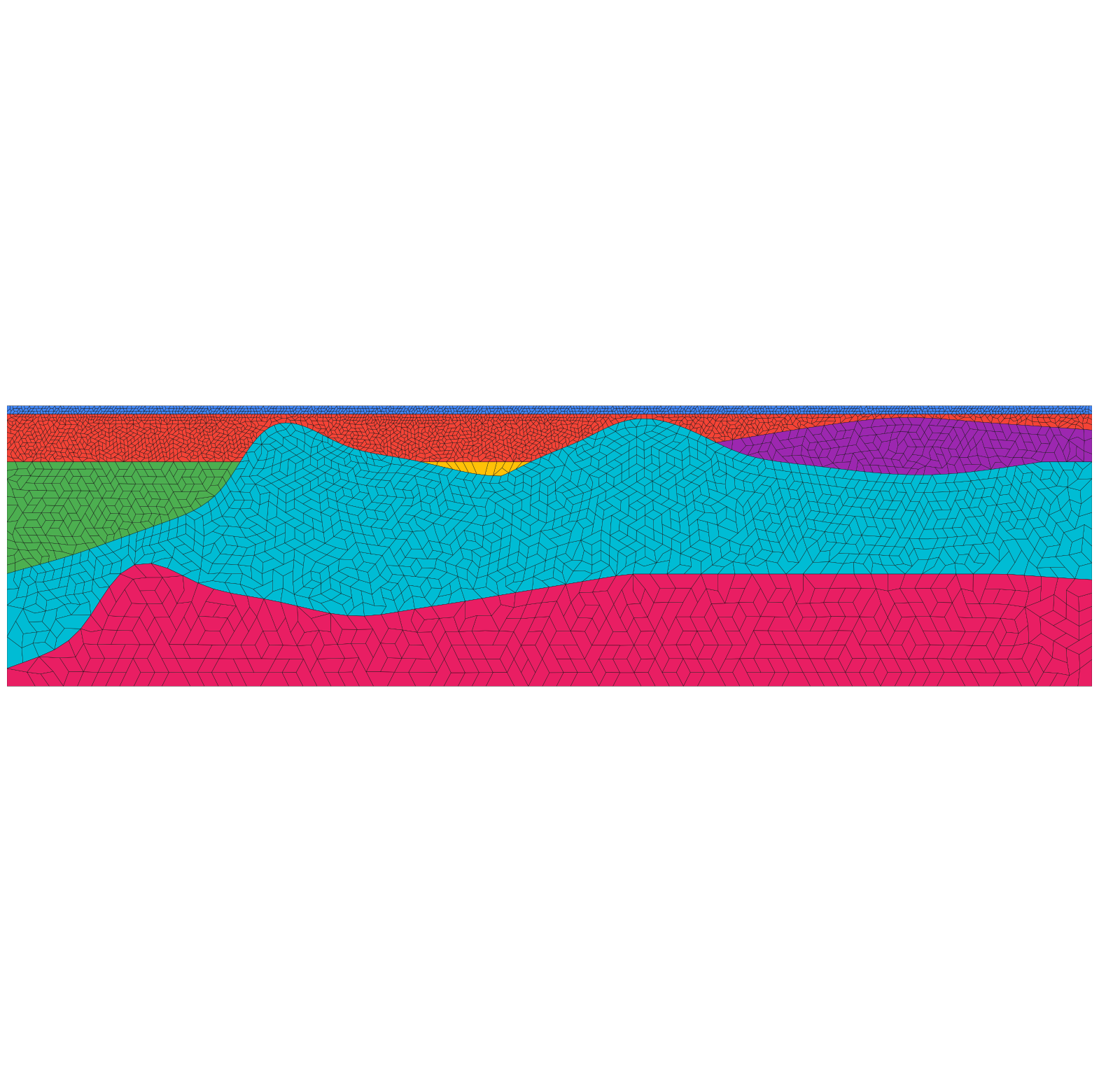}}
    \subfloat[Agglomerated mesh]{\includegraphics[trim=0 300 0 300, clip, width=0.5\textwidth]{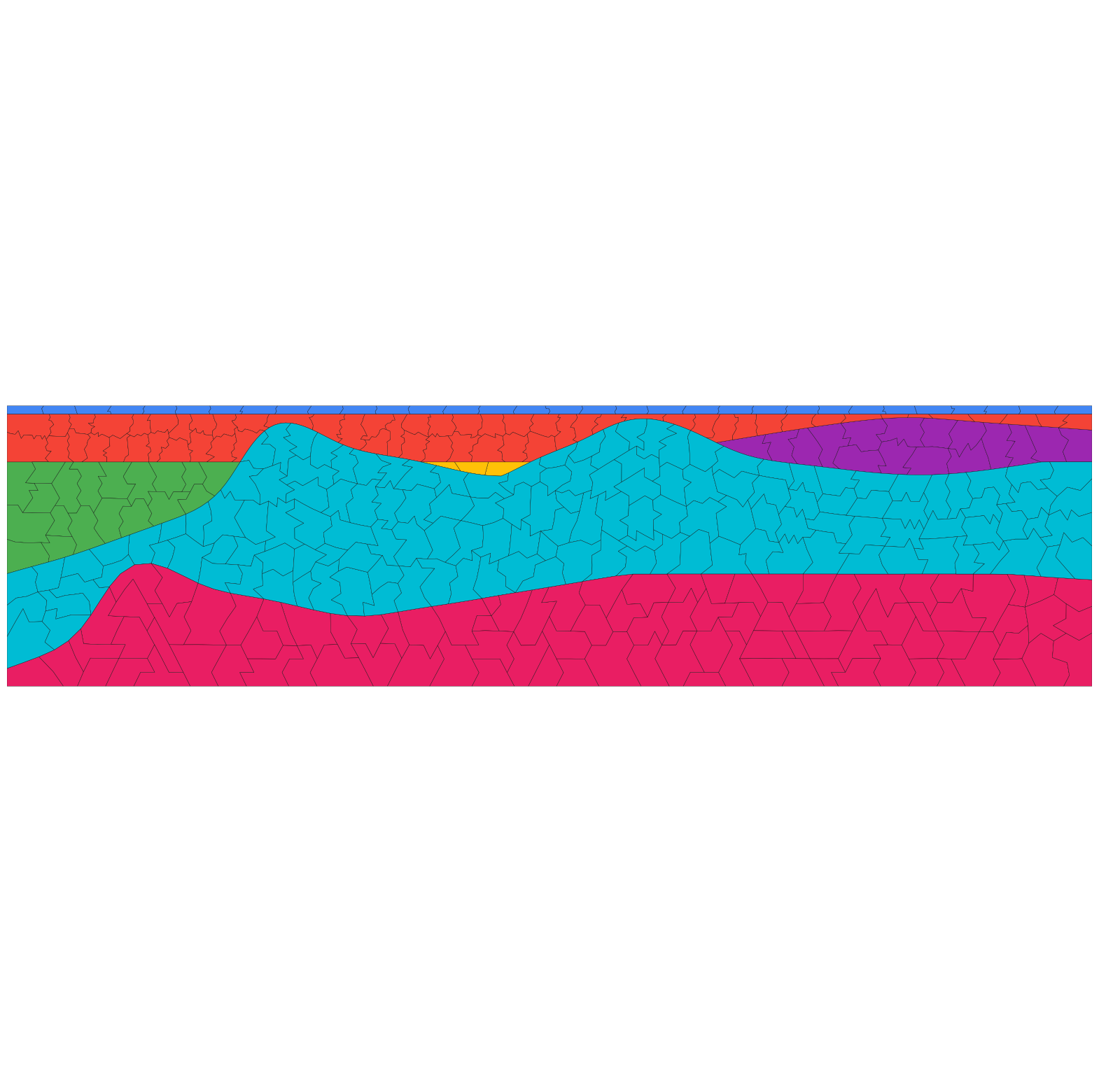}}   
    \caption{Test Case 3: an unstructured polygonal grid of a layered medium; materials with different physical properties are denoted with different colors. On the left, the original mesh with 4870 elements; on the right, the corresponding agglomerated mesh obtained with the SAGE-Base model using \emph{segregated} mode and multiplicative factor of 0.04, resulting in 363 elements.}
    \label{fig:mesh_emilia}
  \end{figure}
  
  \subsection{Test Case 3: Two-Dimensional Layered Medium}
  In this test case, we showcase the ability of \maggnn to handle physical domains associated with multiple physical quantities.
  As explained in Section~\ref{sec:sagehetero}, SAGE-Heterogeneous can only handle two distinct physical parts at most. If we consider a mesh with more than two physical parameters, the only option is to use the \emph{segregated} mode. We demonstrate this by agglomerating a mesh of a layered medium, including seven different physical parameters, using the SAGE-Base model. The mesh represens an idealized bidimensional Earth's cross-section $\Omega = (0, 38.4) \text{km} \times (0, 10) \text{km}$ used for the simulation of seismic wave propagation \cite{antonietti2021PolyDGgeophysicalsimulations}. The results are reported in Figure~\ref{fig:mesh_emilia}.
  
  \subsection{Test Case 4: Hybrid Mesh of a Three-Dimensional Domain}
\begin{figure}[!t]
    \centering
    \subfloat[Section of the original mesh]{\includegraphics[trim=450 230 450 270, clip, width=0.3\textwidth]{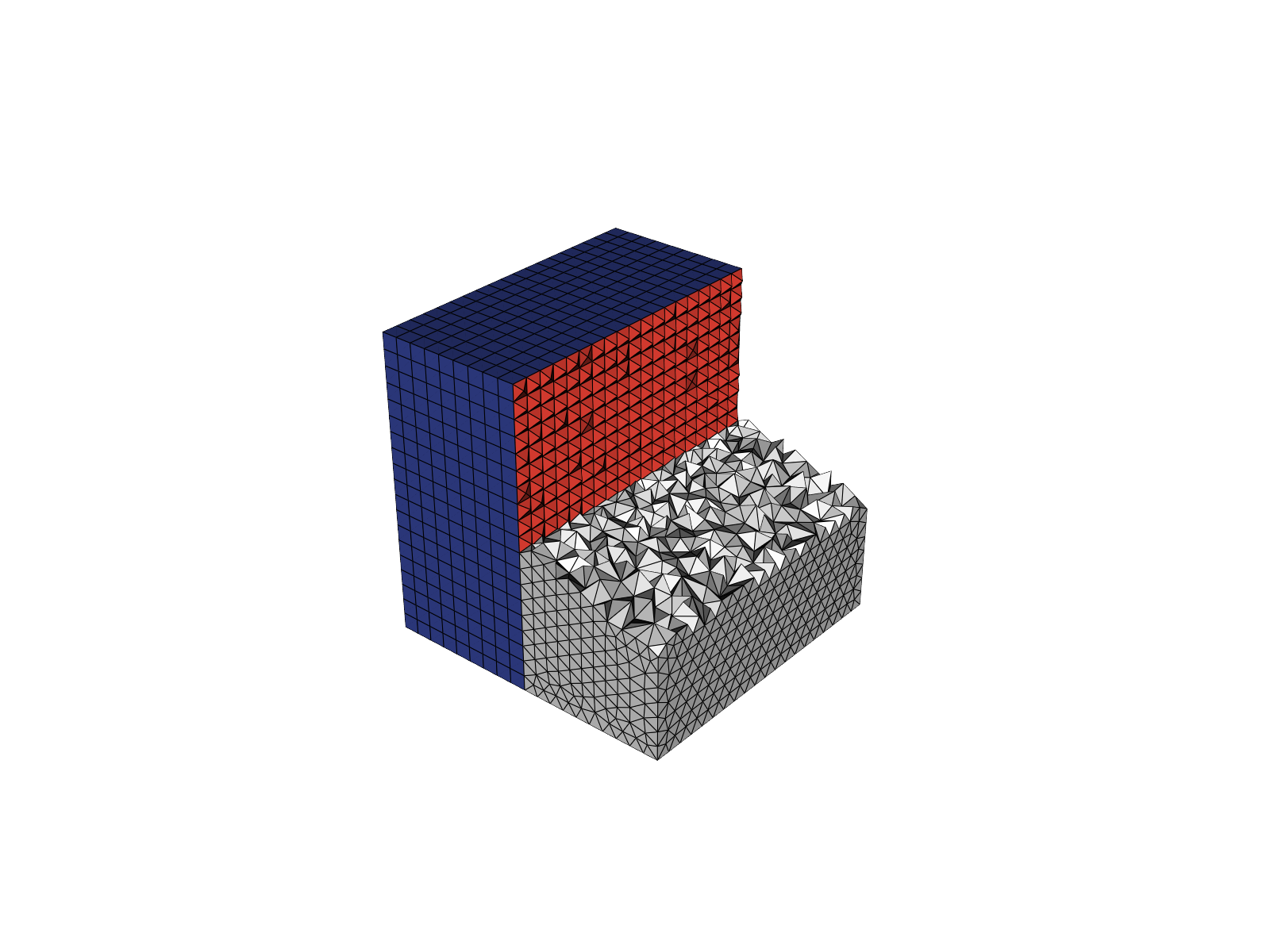}} \quad
    \subfloat[Agglomerated mesh]{\includegraphics[trim=450 230 450 270, clip, width=0.3\textwidth]{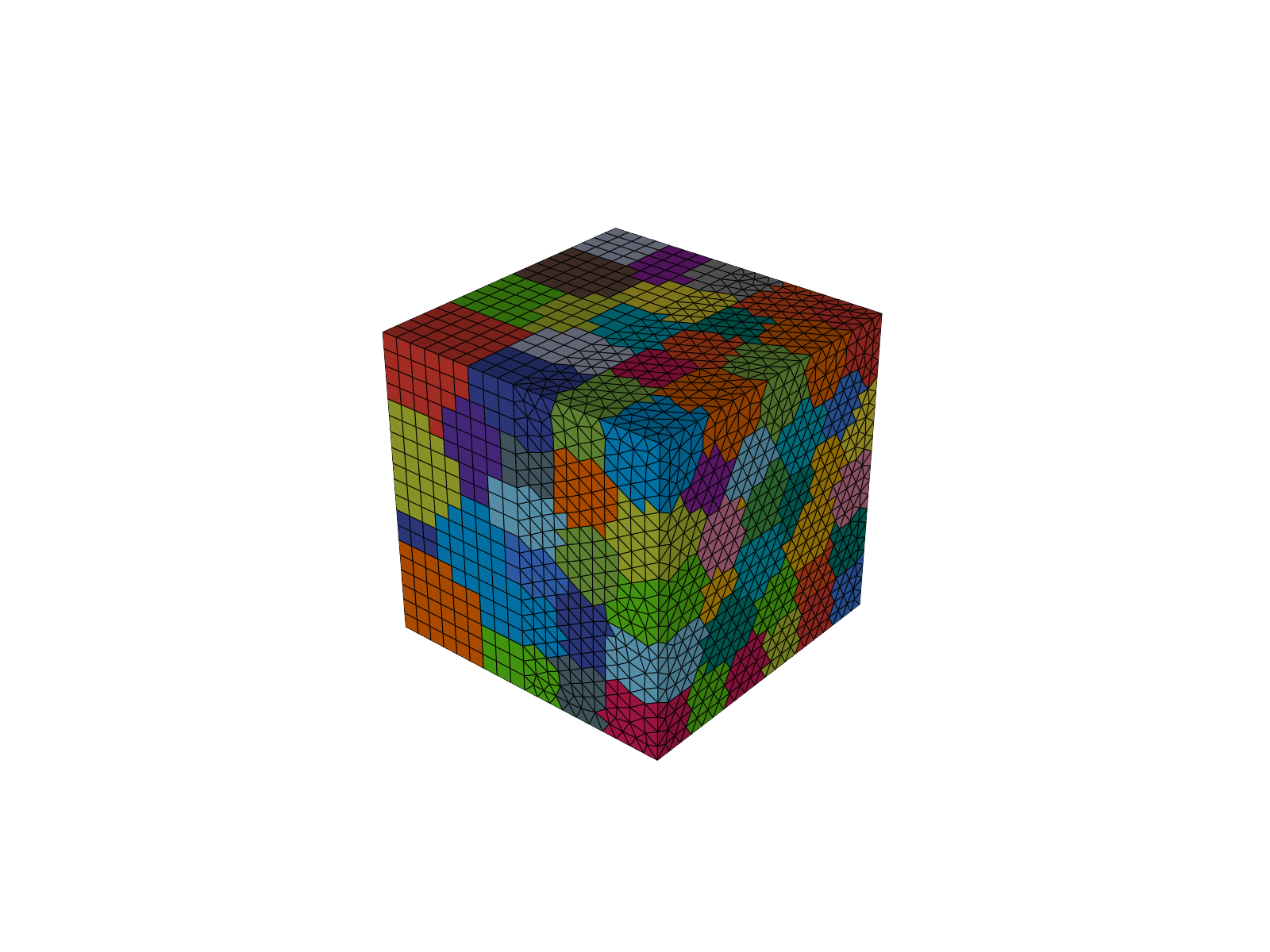}} \quad
    \subfloat[Exploded view]{\includegraphics[trim=450 230 450 270, clip, width=0.3\textwidth]{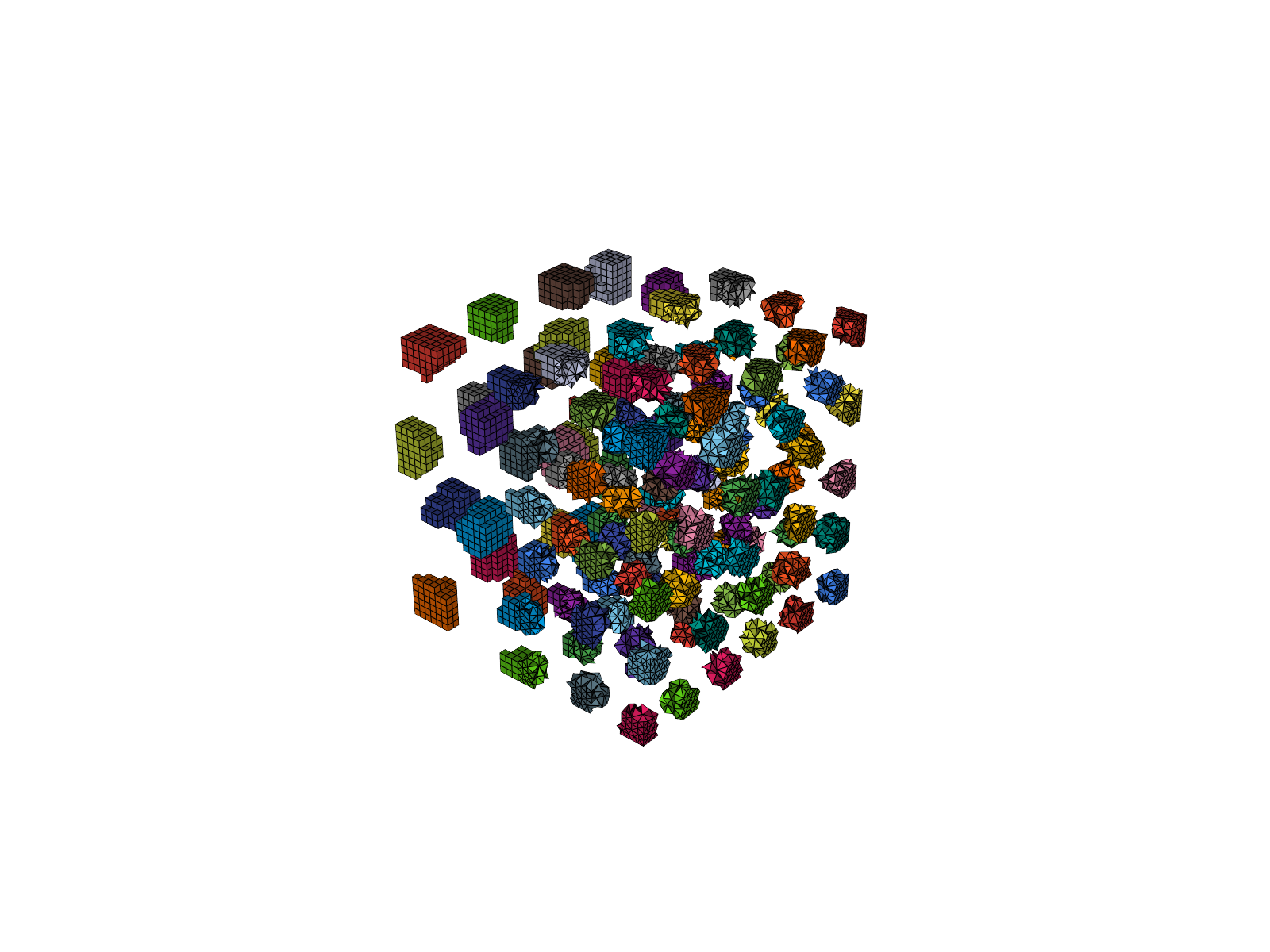}}  
    \caption{Test Case 4: mesh formed by 27484 cells, comprised of 23874 tetrahedra, 3249 hexahedra, and 361 pyramids. On the left, a cut of the original mesh highlighting in red the pyramids for the transition between the white tetrahedra and blue hexahedra. On the center and right, the mesh agglomerated by k-means with \py{k=128}, and its exploded view.}
    \label{fig:mixed_mesh}
  \end{figure}

  In the following test case, we want to highlight that \maggnn can correctly handle three-dimensional meshes that include more than one type of cell. To this end, we consider a mesh of the unit cube formed by tetrahedra and hexahedra, with pyramids for the transition between the two, that has been generated using \py{Gmsh}. The mesh was agglomerated using k-means with $k=128$; the result is shown in Figure~\ref{fig:mixed_mesh}.
    
  \subsection{Test Case 5: Three-Dimensional Brain Geometry}

\begin{figure}[!t]
  \centering  
  \subfloat[View from above]{\includegraphics[trim=400 200 400 200, clip, width=0.3\textwidth]{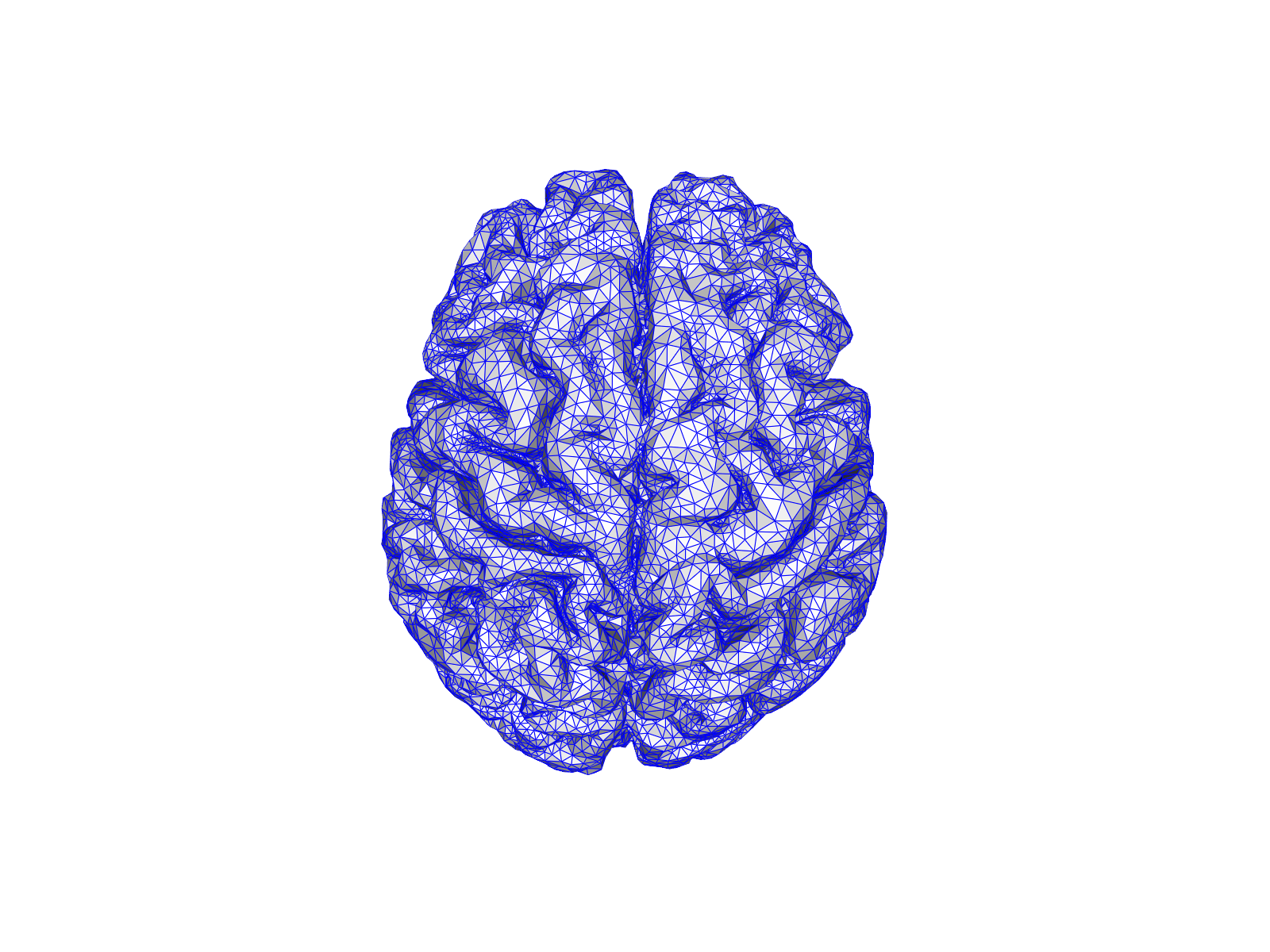}} \qquad
  \subfloat[Lateral view]{\includegraphics[trim=400 200 400 200, clip, width=0.3\textwidth]{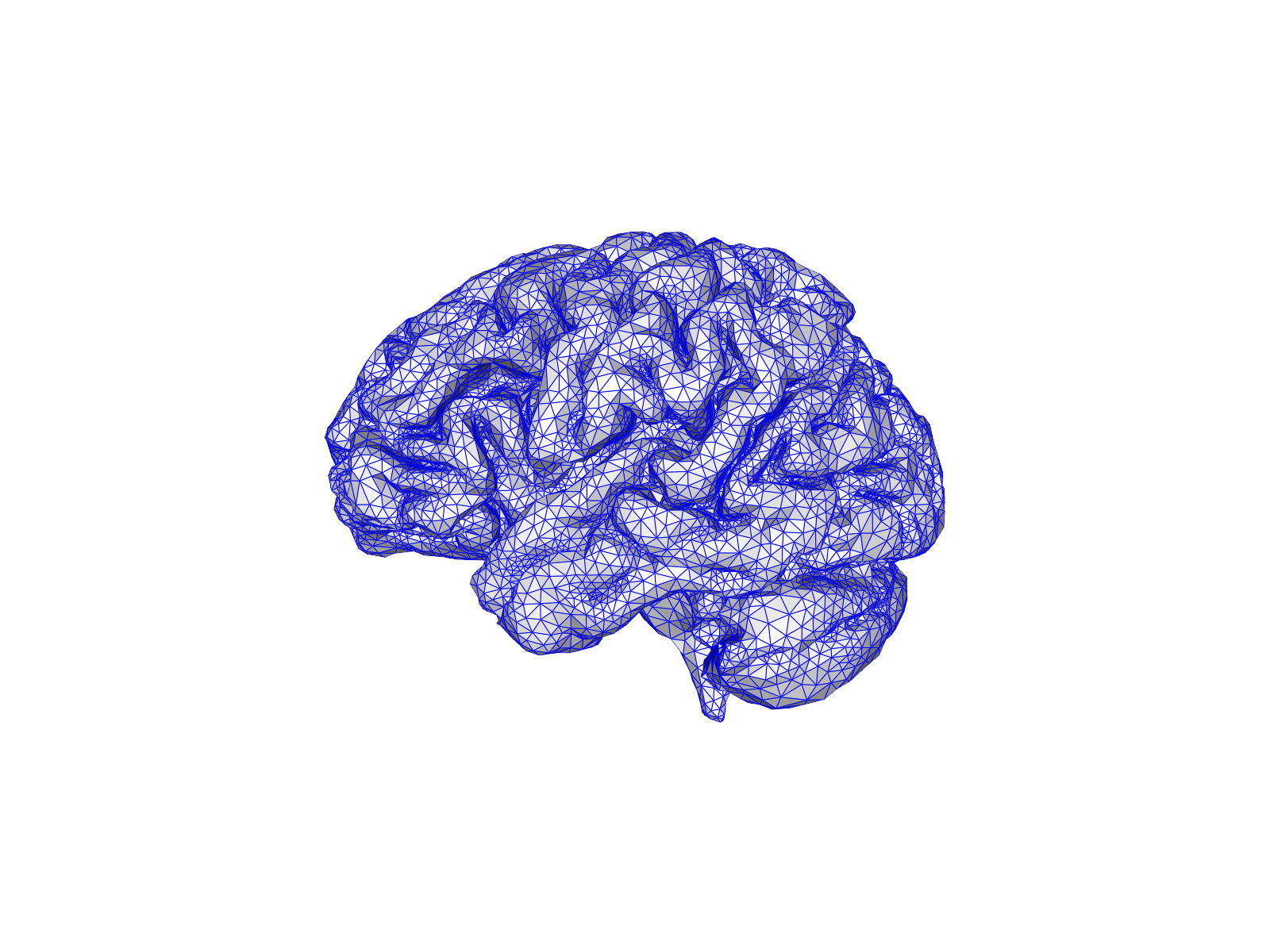}}    
  \caption{Test Case 5: Original mesh of the whole brain coming from Magnetic Resonance Imaging.}
  \label{fig:whole_brain_original}
\end{figure}
\begin{figure}[!t]
  \centering   
  \subfloat{\includegraphics[trim=450 150 450 150, clip, width=0.22\textwidth]{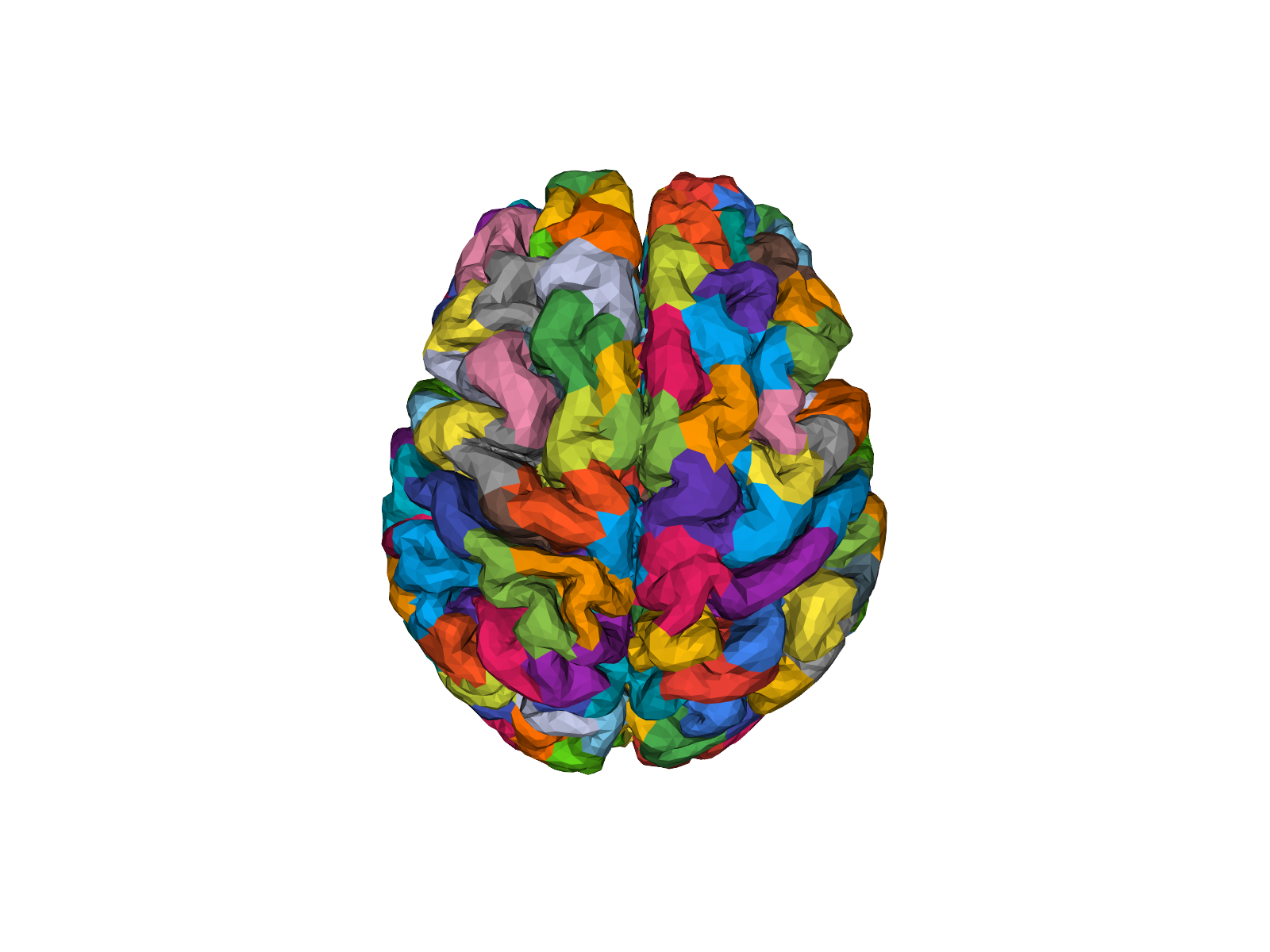}} \quad
  \subfloat{\includegraphics[trim=450 150 450 150, clip, width=0.22\textwidth]{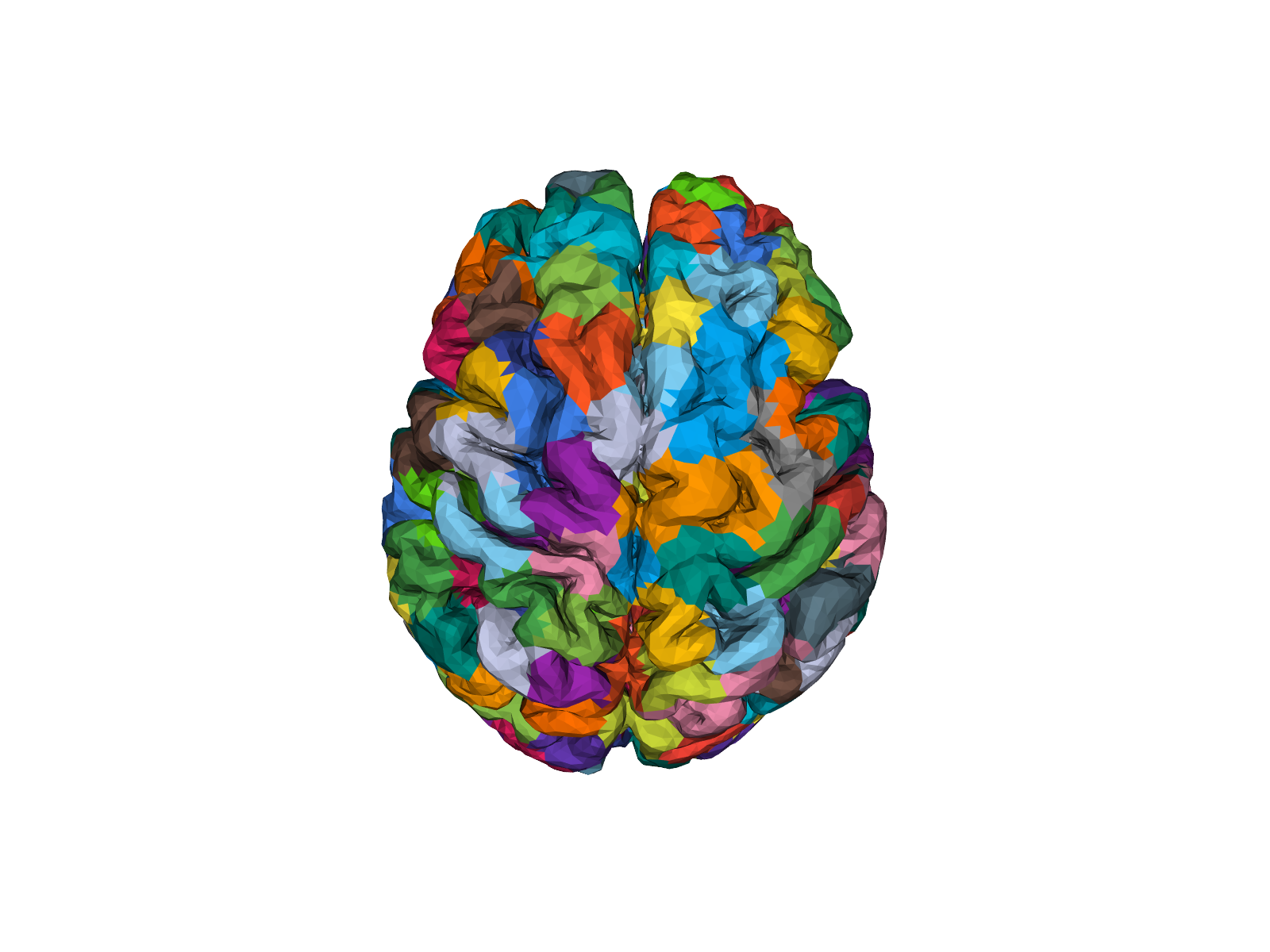}} \quad
  \subfloat{\includegraphics[trim=450 150 450 150, clip, width=0.22\textwidth]{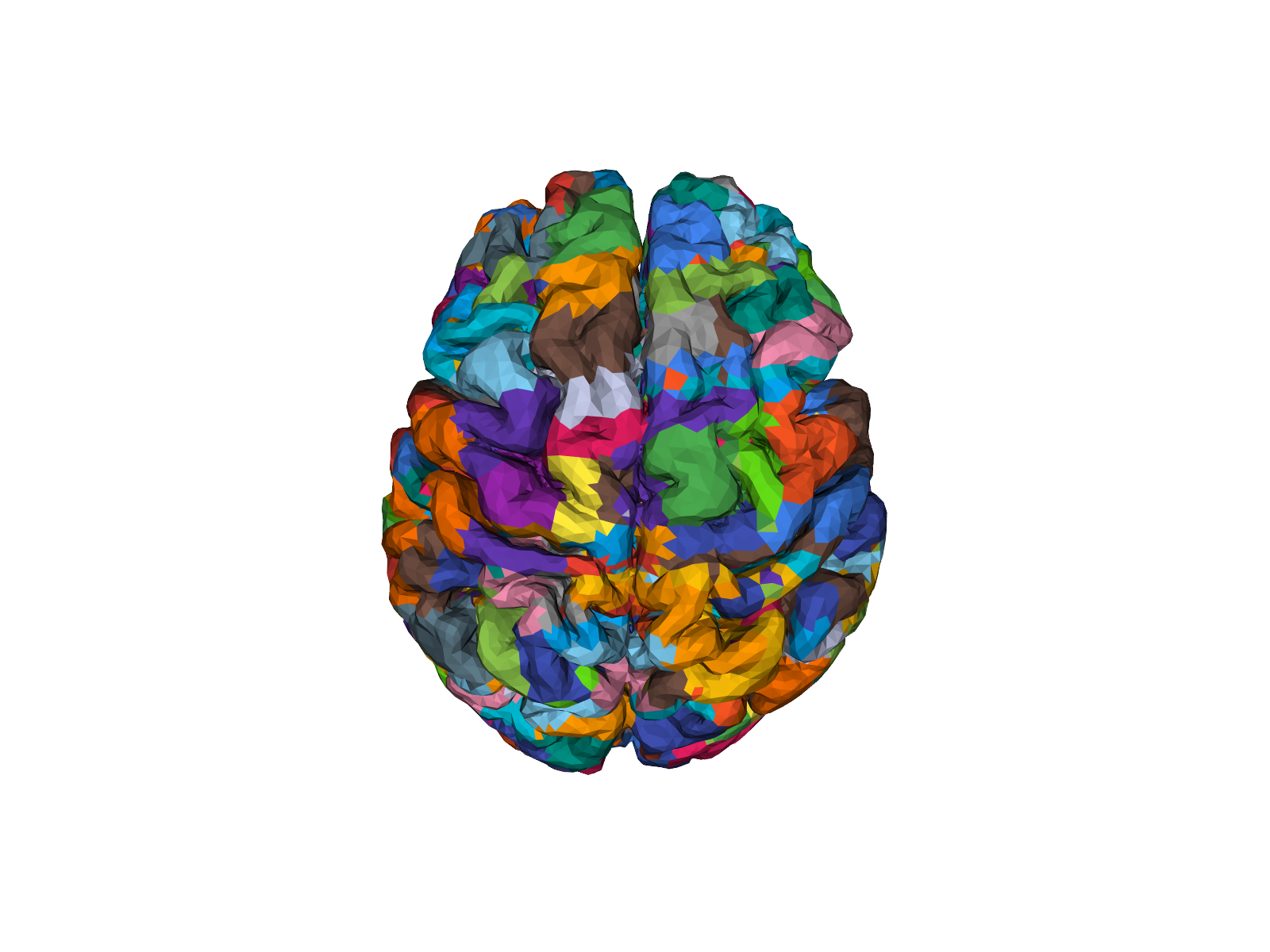}} \quad
  \subfloat{\includegraphics[trim=450 150 450 150, clip, width=0.22\textwidth]{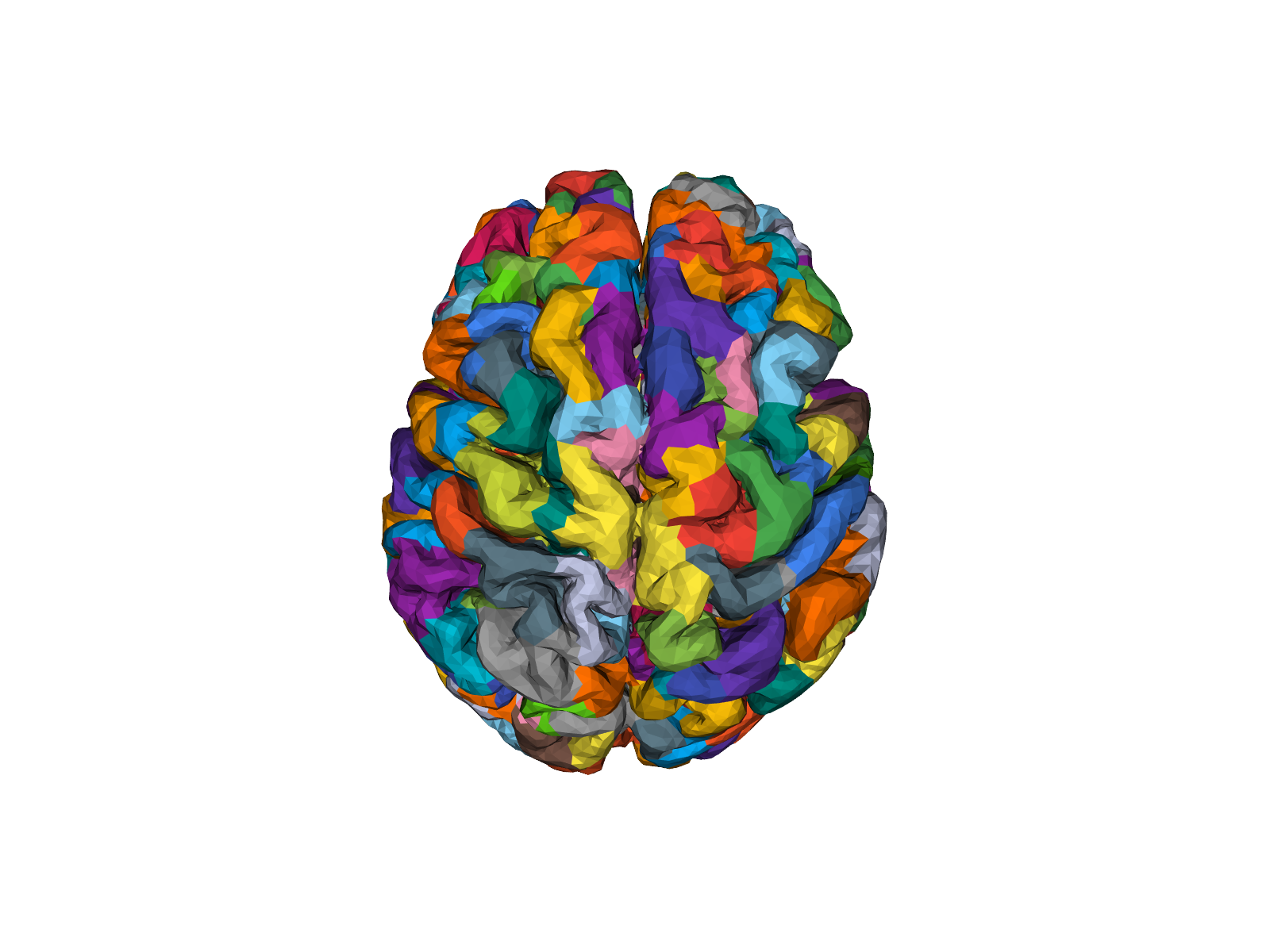}} \\
  \vspace{0.2cm}
  \subfloat[METIS]{\includegraphics[trim=450 150 450 150, clip, width=0.22\textwidth]{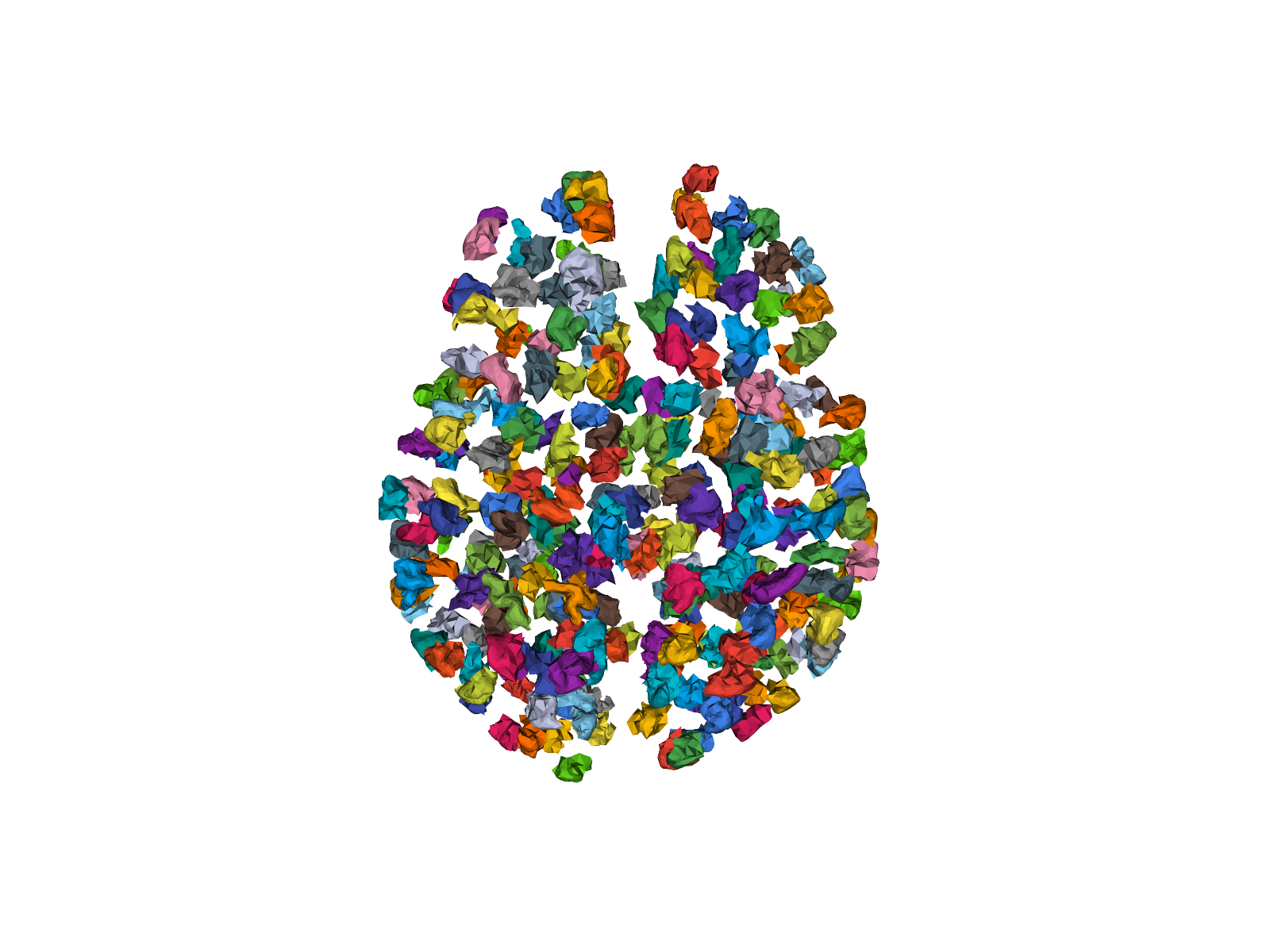}} \quad
  \subfloat[k-means]{\includegraphics[trim=450 150 450 150, clip, width=0.22\textwidth]{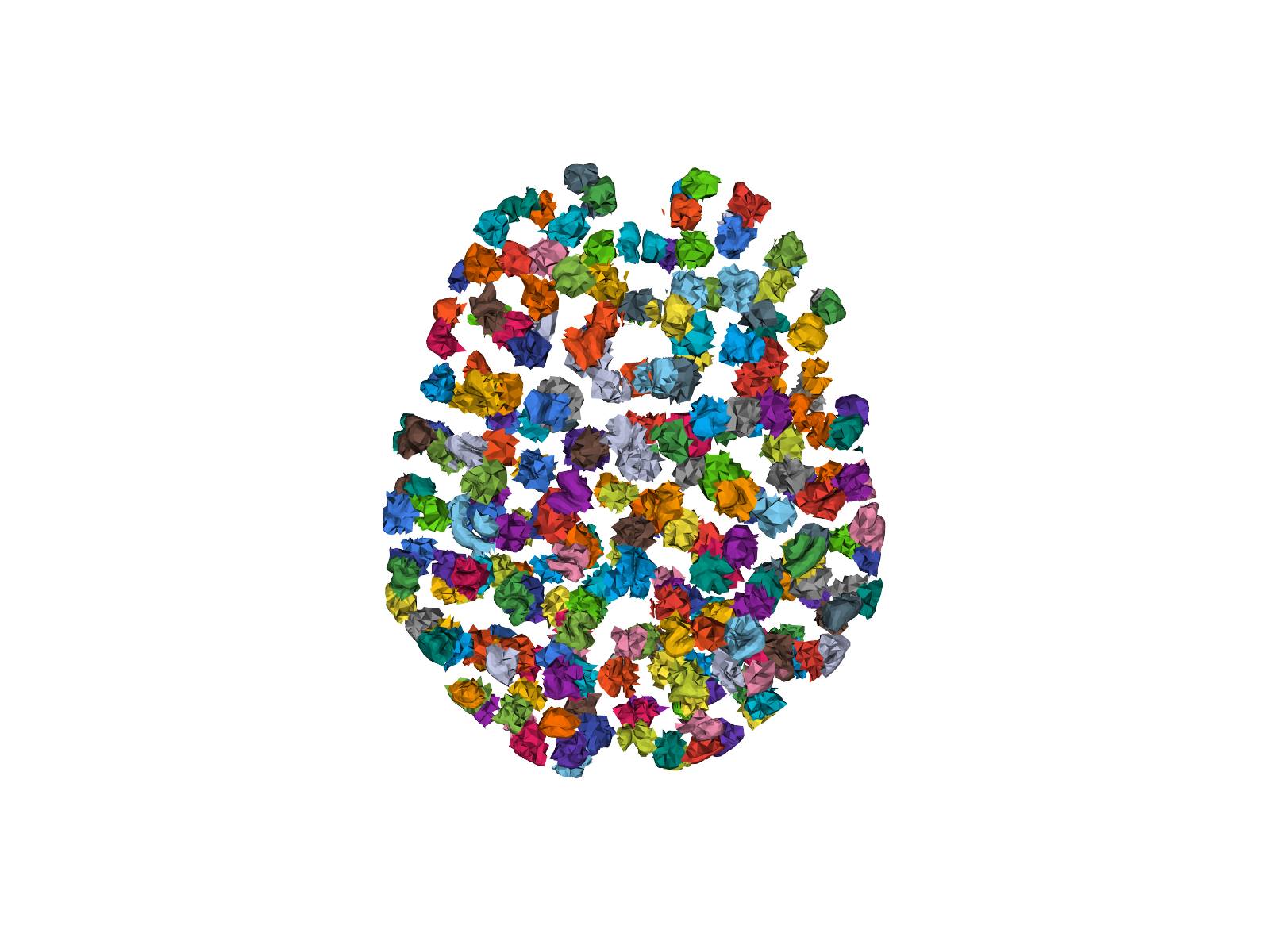}} \quad
  \subfloat[SAGE-Base]{\includegraphics[trim=450 150 450 150, clip, width=0.22\textwidth]{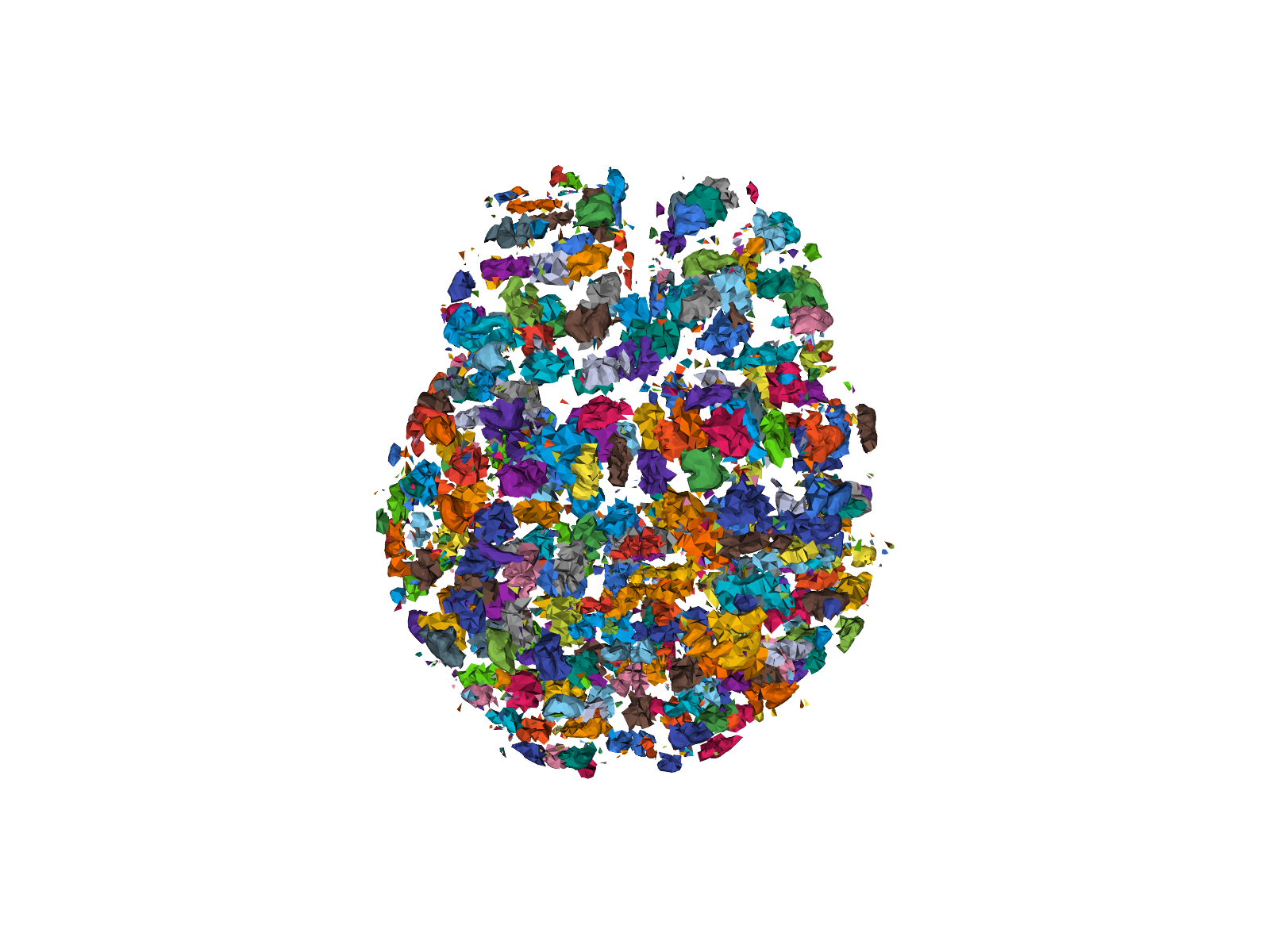}} \quad
  \subfloat[SAGE-Base + RL Refiner]{\includegraphics[trim=450 150 450 150, clip, width=0.22\textwidth]{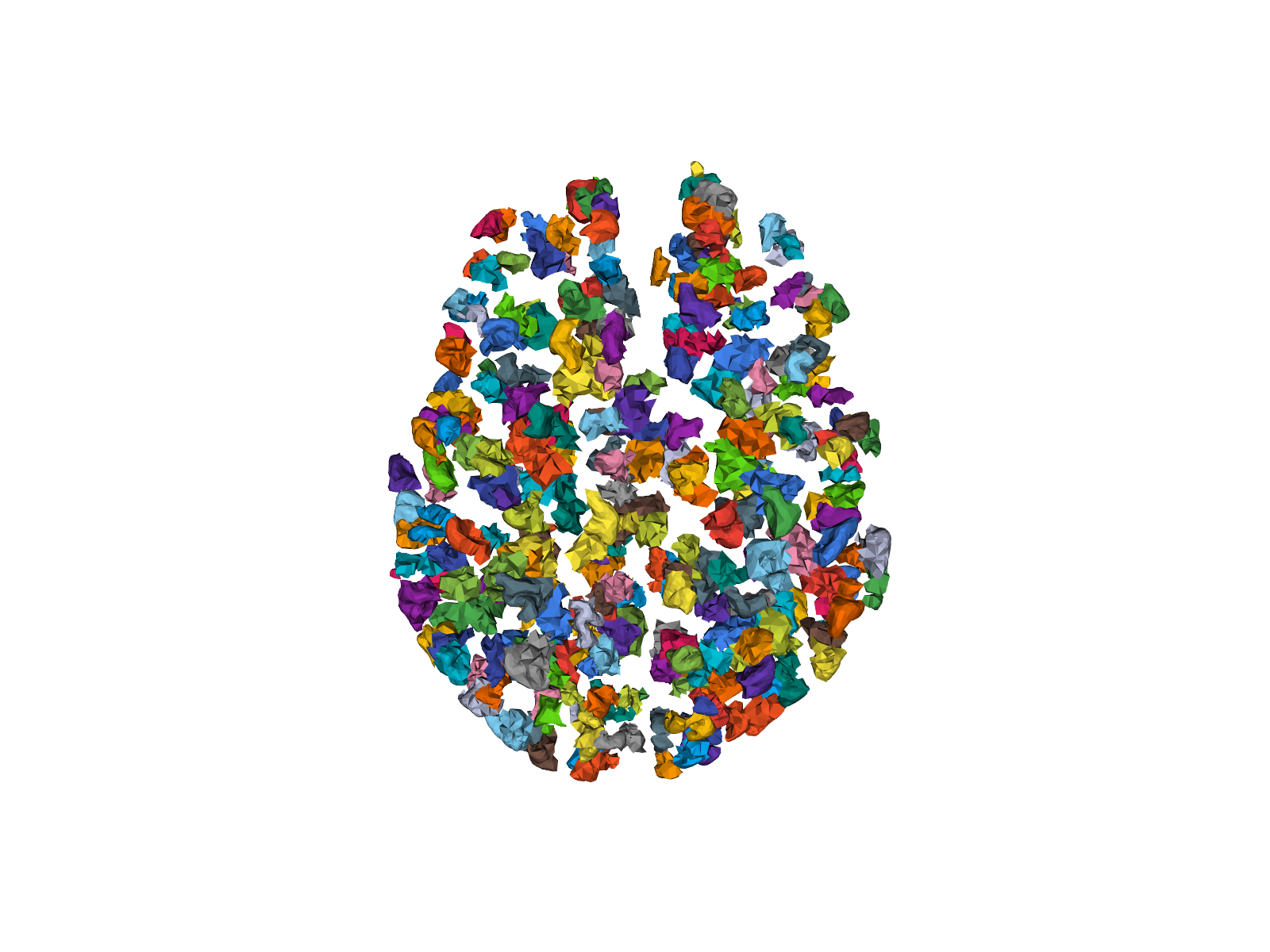}} \\
  \vspace{0.2cm}
  \subfloat[Quality metrics box plots]{\includegraphics[width=1\textwidth]{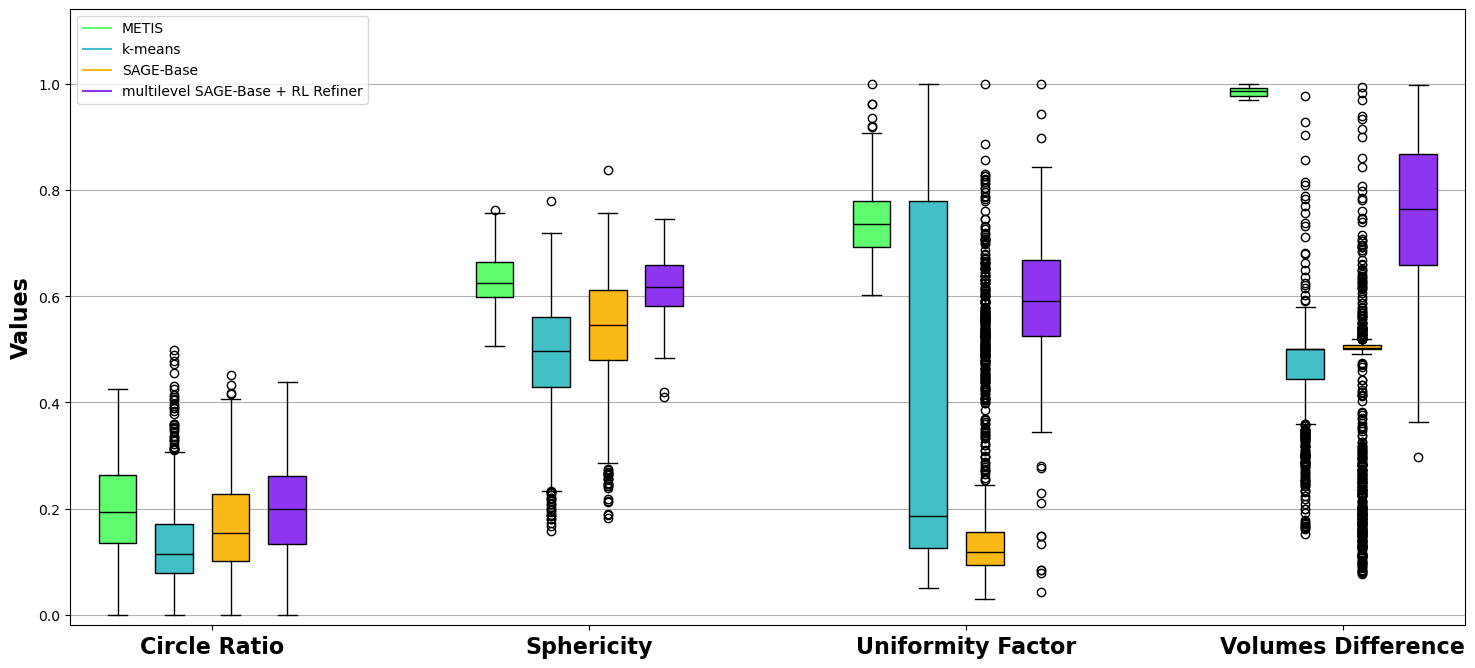}}
  \vspace{0.2cm}
  \caption{Test Case 5: an unstructured mesh of a human brain consisting of 123383 tetrahedra, agglomerated using different methods (METIS, k-means, SAGE-Base, SAGE-Base in multilevel framework with RL Refiner) and their exploded view, together with the box plots of the computed quality metrics (CR, $\Psi$, UF and $\widetilde{\text{VD}}$, defined as in Eq.~\eqref{eq:CR}-\eqref{eq:VD}).}
  \label{fig:whole_brain_agg}
\end{figure}

  To demonstrate the generalization capabilities of the models in the three-dimensional case, we consider once again a mesh coming from an MRI scan of a human brain, shown in Figure~\ref{fig:whole_brain_original}. We agglomerate the mesh using four different models (METIS, k-means, SAGE-Base, and SAGE-Base as coarse partitioner in a multilevel framework with RL Refiner). For the GNN models, we used \emph{number of refinements} mode with \py{nref=8}, while for the other two \emph{direct k-way} with \py{k=256} was used. In Figure~\ref{fig:whole_brain_agg} we report the agglomerated meshes and their exploded view, together with their quality metrics box plots.
  Most of the same considerations of the two-dimensional case hold, with some differences: due to the highly corrugated surface of the brain, k-means and SAGE-Base struggle to create fully connected elements, often leaving a few tetrahedra separated from the rest. As a consequence, the corresponding elements will be much smaller, leading to a lower uniformity factor. This issue can be alleviated in the k-means case by properly choosing a higher $k$, while not much can be done for SAGE-Base since the problem arises from the very first bisection.
  The occurrence of disconnected elements is primarily inherent to the architecture of GNNs. In particular, the message-passing layers preserve knowledge of the topology only within a limited neighborhood of each node in the dual graph. Extending the receptive field to capture more global connectivity would mitigate this issue, but at the expense of significantly increasing computational cost.
  On the other hand, the multilevel approach with the RL refinement mostly solves the issue, improving the quality of the partition.

  \subsection{Test Case 6: Three-Dimensional Statue of Garuda and Vishnu}
\begin{figure}[!t]
  \centering  
  \subfloat[Frontal view]{\includegraphics[trim=400 0 400 0, clip, width=0.35\textwidth]{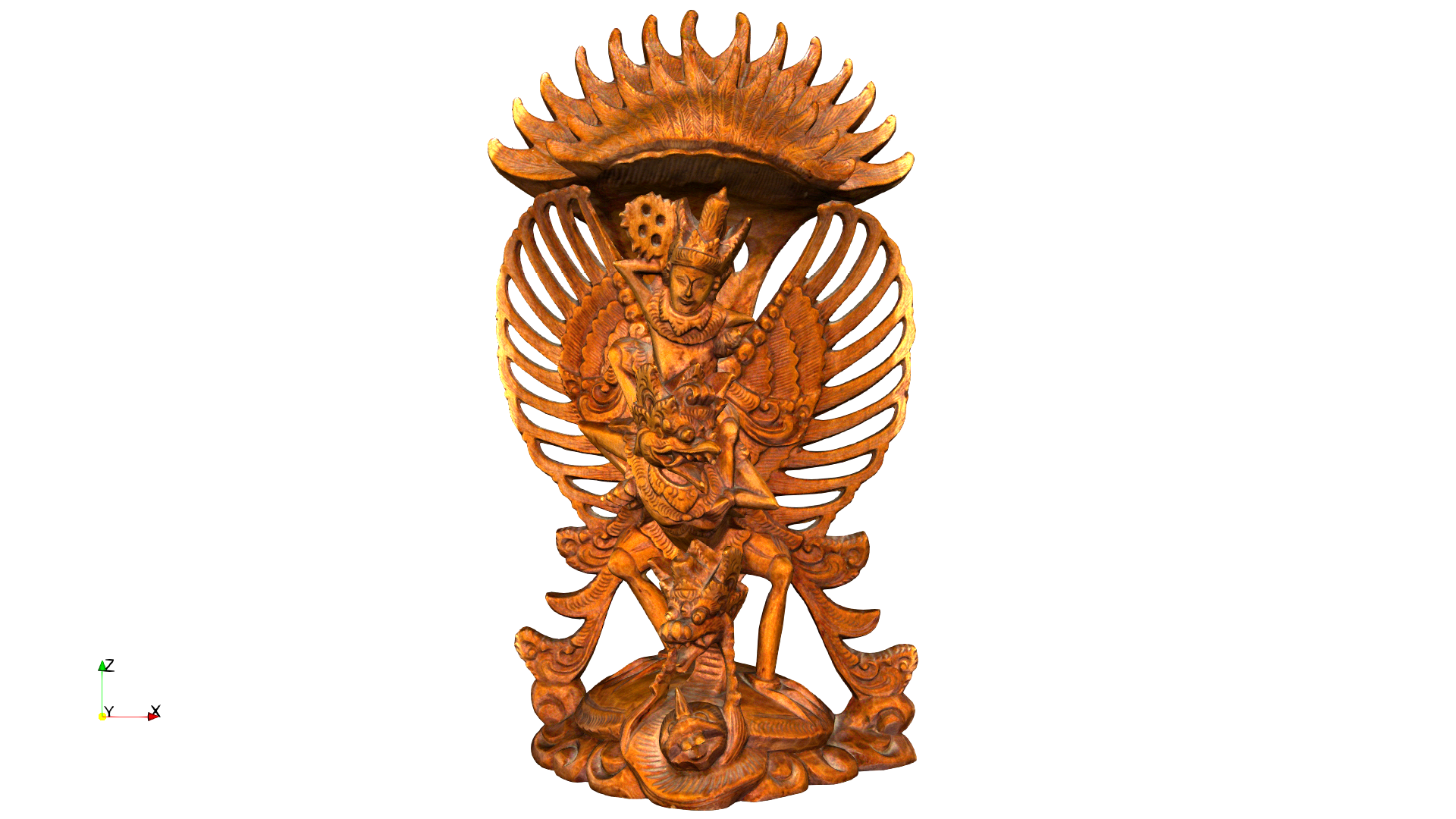}} \qquad
  \subfloat[Lateral view]{\includegraphics[trim=600 40 600 20, clip, width=0.27\textwidth]{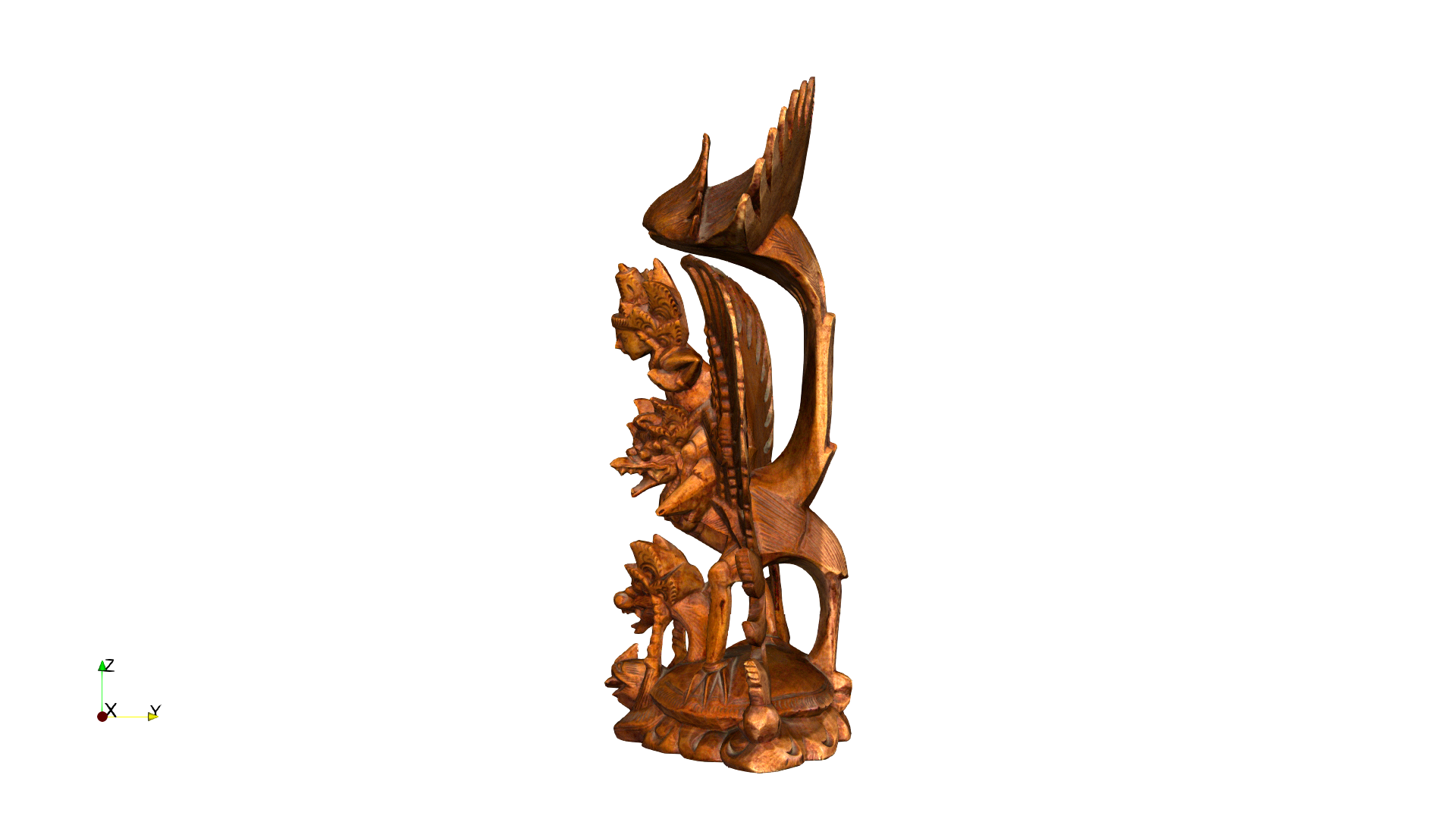}}    
  \caption{Test Case 6: Original mesh of the statue of Garuda and Vishnu coming from a 3D scan.}
  \label{fig:garudavishnu_original} 
\end{figure}

\begin{figure}[!htbp]
  \centering    
  \subfloat{\includegraphics[trim=400 100 400 160, clip, width=0.3\textwidth]{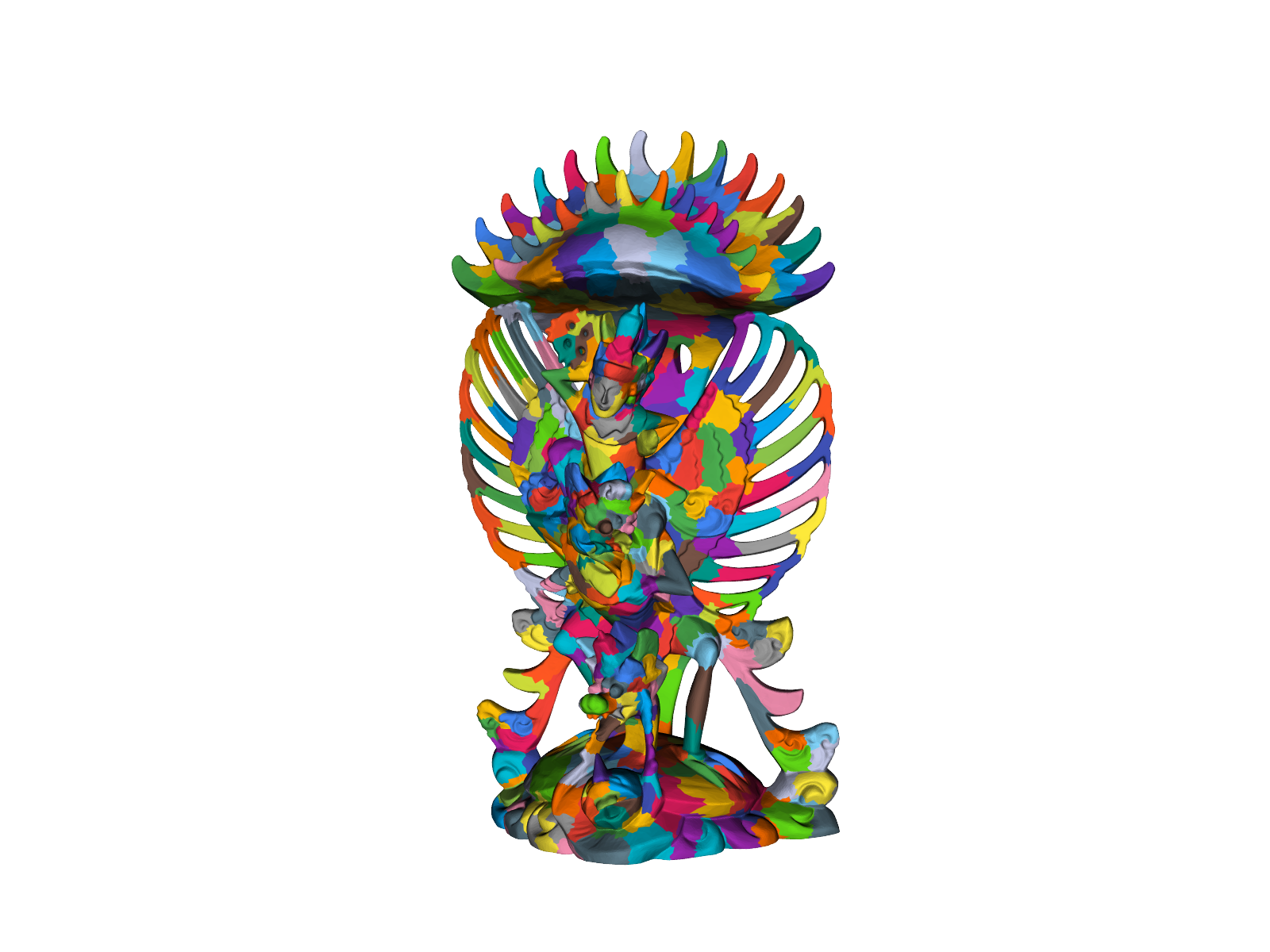}} \quad
  \subfloat{\includegraphics[trim=400 100 400 160, clip, width=0.3\textwidth]{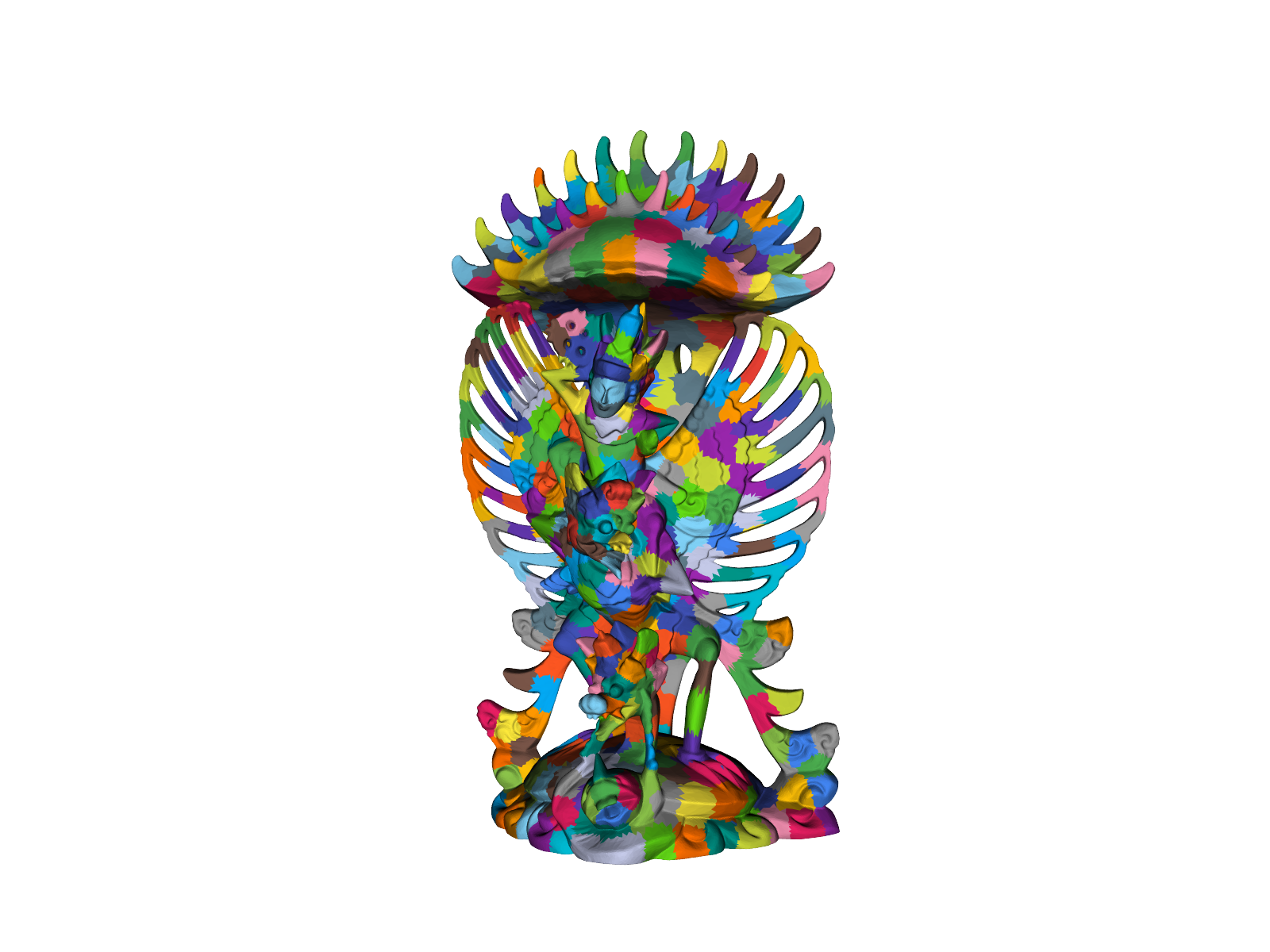}} \quad
  \subfloat{\includegraphics[trim=400 100 400 160, clip, width=0.3\textwidth]{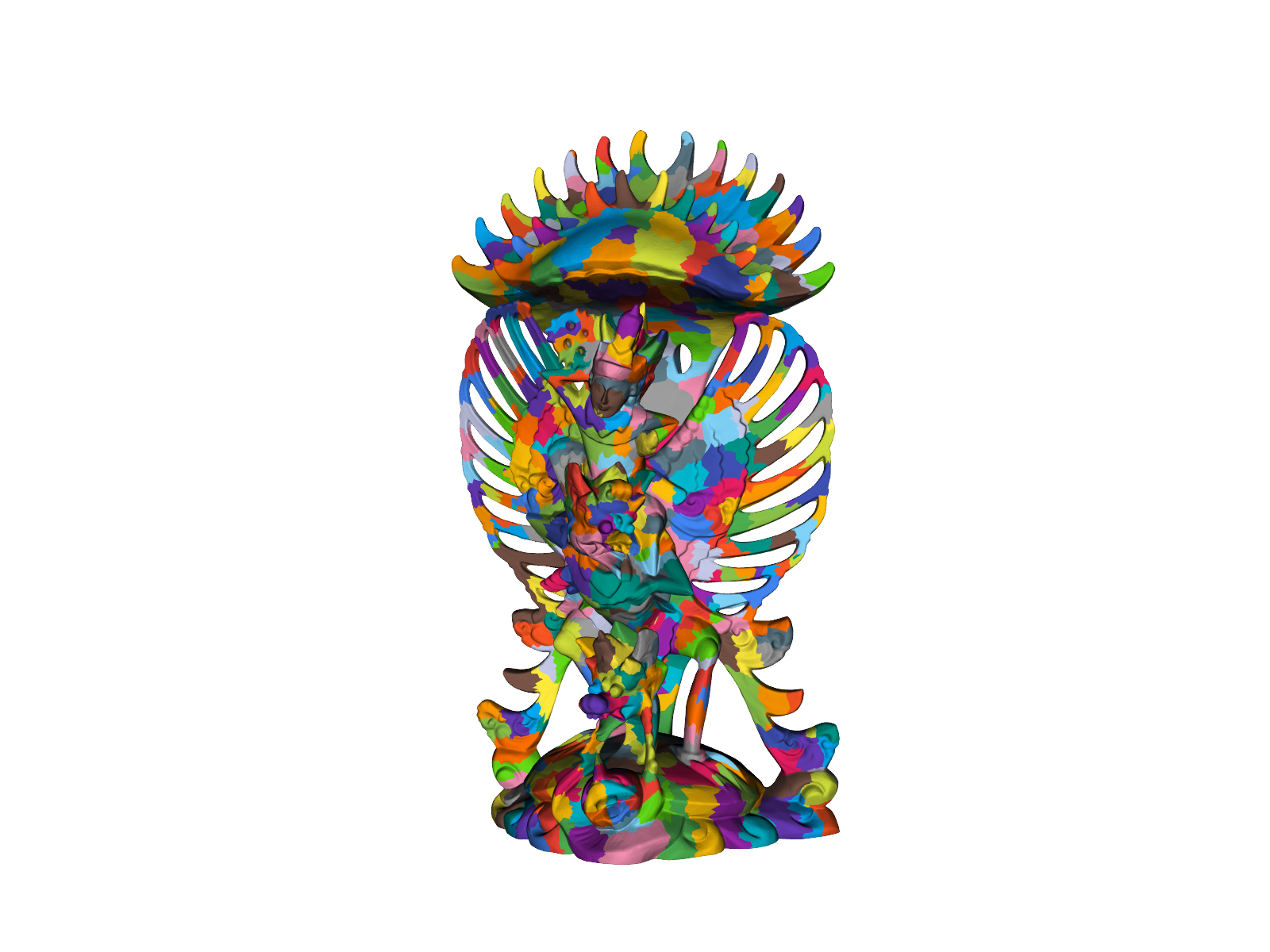}} \\    
  \subfloat[METIS]{\includegraphics[trim=500 100 500 170, clip, width=0.3\textwidth]{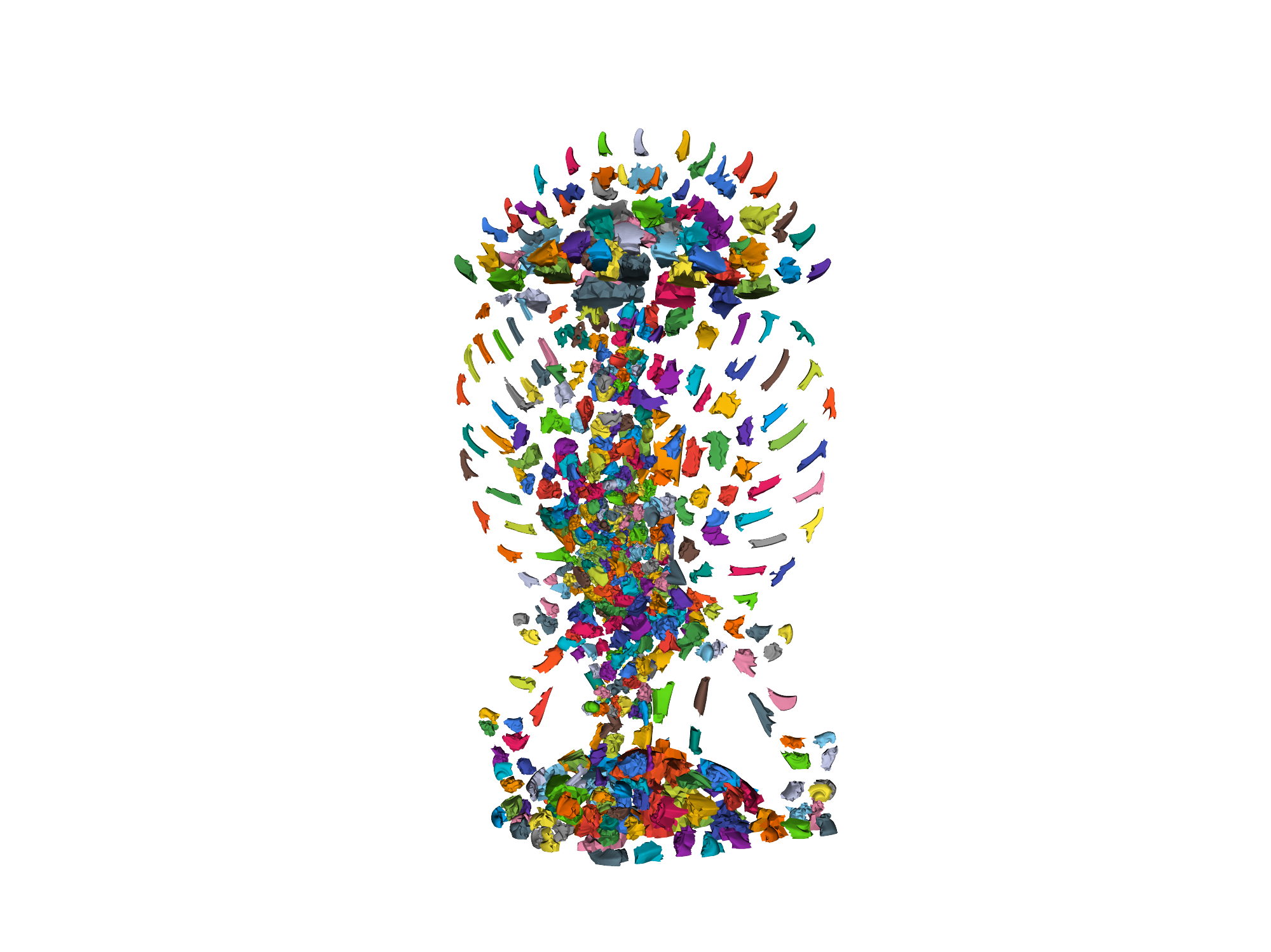}} \quad
  \subfloat[k-means]{\includegraphics[trim=500 100 500 170, clip, width=0.3\textwidth]{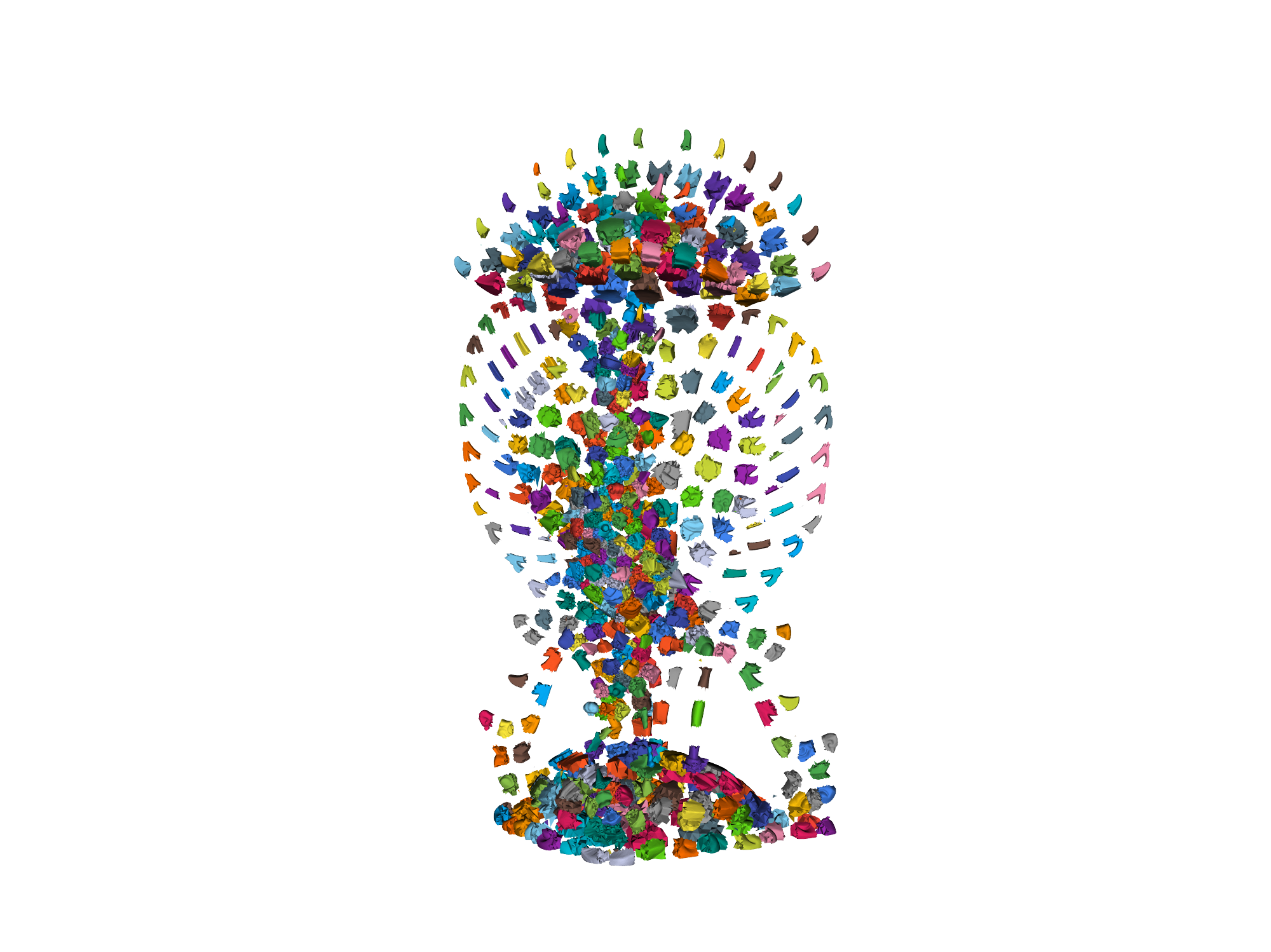}} \quad
  \subfloat[SAGE-Base + RL Refiner]{\includegraphics[trim=500 100 500 170, clip, width=0.3\textwidth]{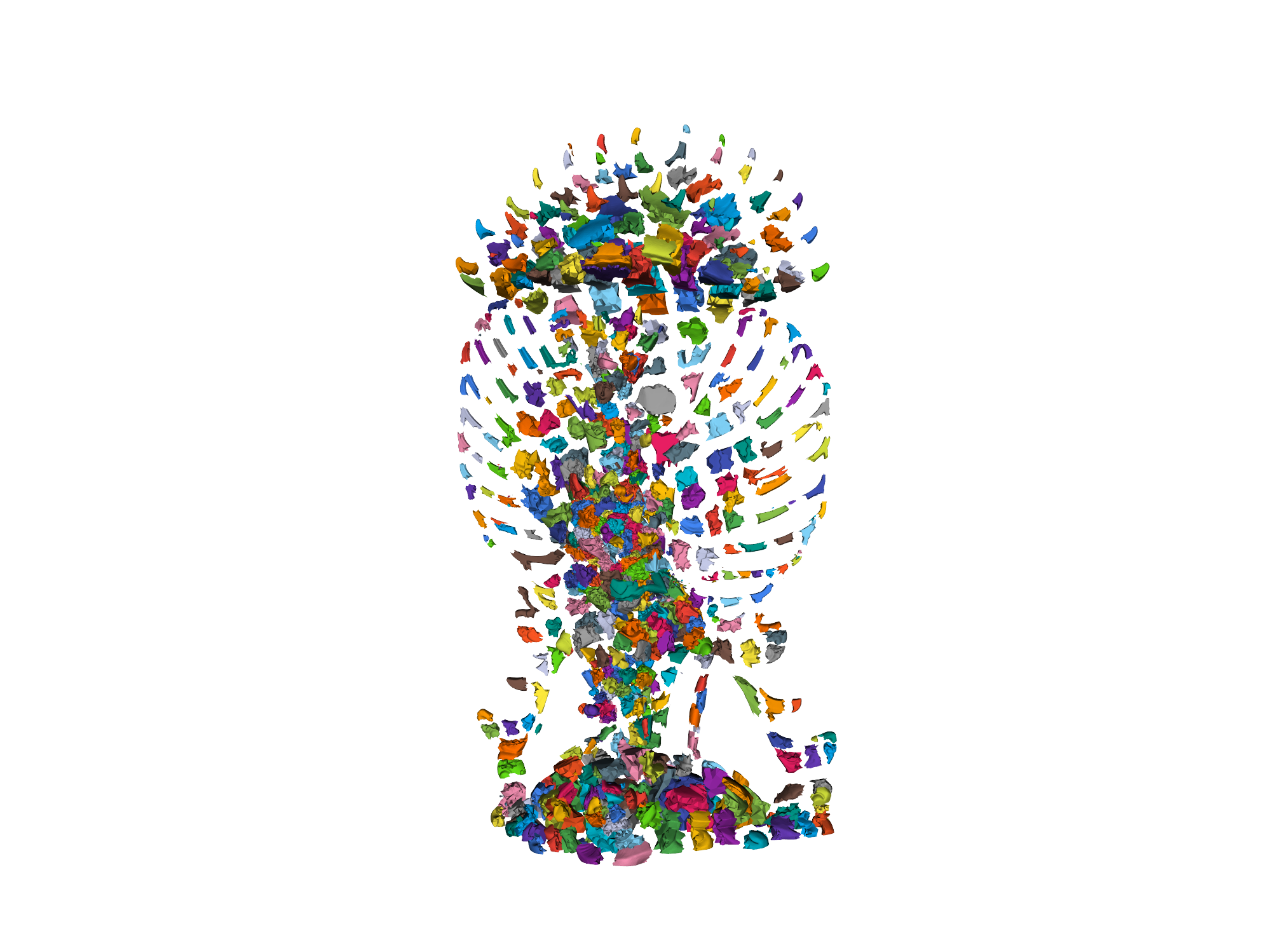}} \\    
  \subfloat[Quality metrics box plots]{\includegraphics[width=1\textwidth]{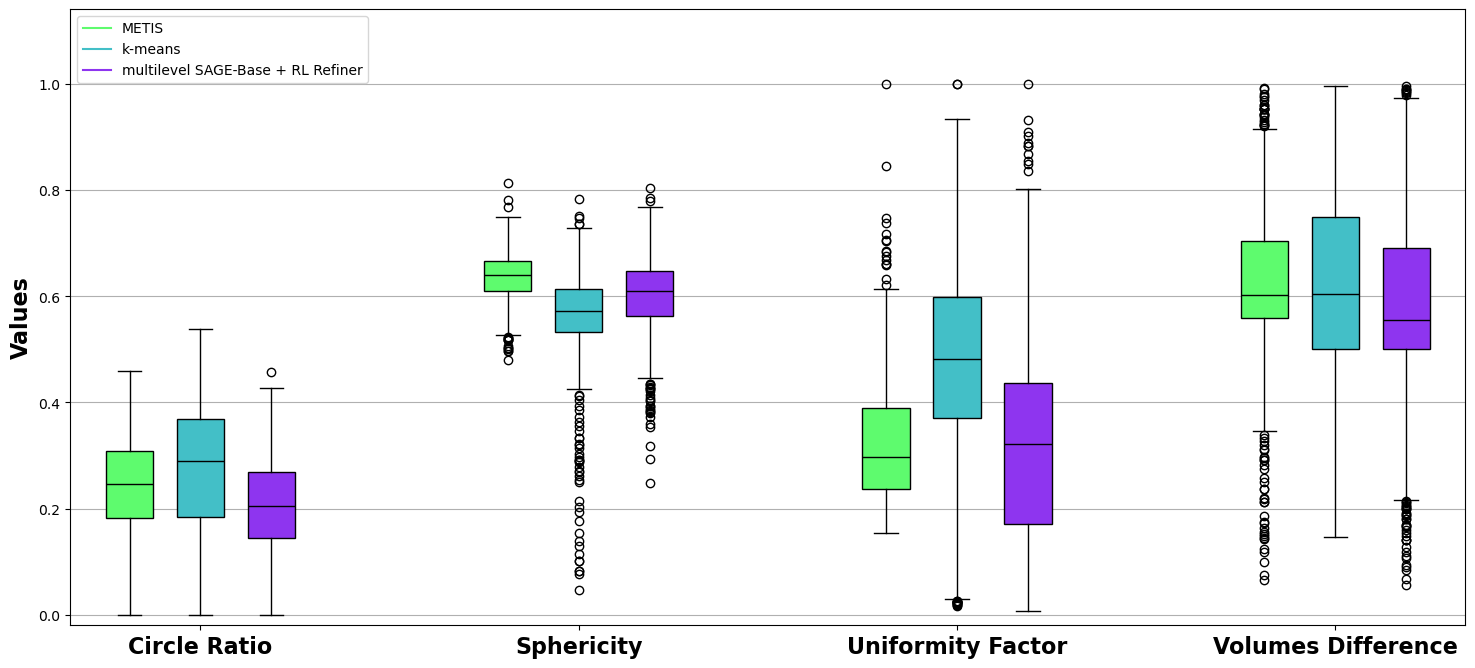}}    
  \caption{Test Case 6: a mesh of a statue of Garuda and Vishnu consisting of 615229 tetrahedra, agglomerated using different methods (METIS, k-means, SAGE-Base in multilevel framework with RL Refiner) and their exploded view, together with the box plots of the computed quality metrics (CR, $\Psi$, UF and $\widetilde{\text{VD}}$, defined as in \eqref{eq:CR}-\eqref{eq:VD}).}
  \label{fig:GarudaVishnu_agg}
\end{figure}

  In this test case, we use a scan of a wooden statue of Garuda and Vishnu, which has very fine geometry and many holes. Thus, it is particularly challenging due to its topology. The tetrahedral mesh has been generated using \py{Gmsh} starting from an STL file coming from a 3D scan of the object (the original STL of Garuda and Vishnu, which is shown in Figure~\ref{fig:garudavishnu_original}, is by \href{http://www.artec3d.com}{Artec3D}).
  We agglomerate it using SAGE-Base in a multilevel framework with RL Refiner and \py{nref=9}, while \emph{direct k-way} with \py{k=512} was used for METIS and k-means; the results are reported in Figure~\ref{fig:GarudaVishnu_agg}.
  Listing \ref{lst:STL_agg} shows the code for the whole agglomeration pipeline of this example when using the multilevel approach. 

 We remark that in this case, METIS fails to partition the graph when passing cell volumes as node weight, so the agglomerated mesh has been created using unitary weights instead. Consequently, the volume uniformity for METIS is much lower compared to previous cases. Other than that, the three methods are fairly comparable across all four metrics.



\begin{lstlisting}[language=Python, caption={Example code showing the whole agglomeration pipeline, from the generation of the mesh from the STL file to the creation of the exploded view plot.}, label=lst:STL_agg]
  from magnet import io, generate, aggmodels
  # Generate tetrahedral mesh from STL file using Gmsh:
  generate.tetrahedra_from_stl('Garuda_Vishnu.stl', remesh=False)
  GV = io.load_mesh('Garuda_Vishnu.vtk')
  sage = aggmodels.SageBase(128,64,4,2).to(aggmodels.DEVICE)  # initialize GNN partitioner
  sage.load_model('models/sage_base.pt')
  lrefiner = aggmodels.RLRefiner(5, 10).to(aggmodels.DEVICE)  # initialize GNN refiner
  lrefiner.load_model('models/rl_refiner.pt')
  agglomerated_GV = sage.agglomerate(GV, mode='multilevel', nref=9, refiner=lrefiner, threshold=500)
  io.exploded_view(agglomerated_GV, scale=0.75, orientation=(30, 50, 0))
\end{lstlisting}
  \section{Integration with the \texttt{lymph} Library}
  \label{sec:lymph_interface}
  \texttt{lymph} is an open-source library for the discretization of multiphysics partial differential equations based on discontinuous Galerkin methods on polytopal grids \cite{antonietti2024lymph}.
  \maggnn comes with \py{Agglomerate.m}, which is a function that allows to agglomerate a mesh by calling the \py{agglomerate.py} script and convert it to \texttt{lymph} format directly from Matlab, making it possible to easily solve differential problems on it. In this section, we give an example usage of this feature, showing that agglomerated meshes created by \maggnn are suitable for PolyDG discretization. 
  
  \subsection{Verification Test 1: Poisson Problem}
  We first provide a brief guide on how to use the \texttt{lymph} API by considering the solution of a Poisson problem in an open and bounded domain $\Omega \subset (0, 1)^2$, cf. Figure~\ref{fig:polydg_poisson_brain},
  \begin{equation}
    \label{eq:Poisson_problem}
    \begin{cases}
      -\Delta u = f(\textbf{x})          & \text{ in } \Omega, \\
      u(\textbf{x}) = g(\textbf{x})      & \text{ on } \partial \Omega.
    \end{cases}
  \end{equation}
  We discretize problem \eqref{eq:Poisson_problem} with a PolyDG method as shown in \cite{antonietti2024lymph}.
  For simplicity, we consider as analytical solution $u(\textbf{x}) = \sin(2\pi \textbf{x}_1) \cos(2\pi \textbf{x}_2)$ and compute the right-hand side $f$ and boundary condition $g$ accordingly.
  We consider two different geometries for $\Omega$. First, we consider the unit square $\Omega=(0, 1)^2$ as a baseline benchmark and then the brain slice we used in Section~\ref{sec:brain_slices} as validation for real-world geometries.
  \subsubsection{Poisson Problem in the Unit Square}
\begin{figure}[!t]
    \centering 
    \includegraphics[width=1\textwidth]{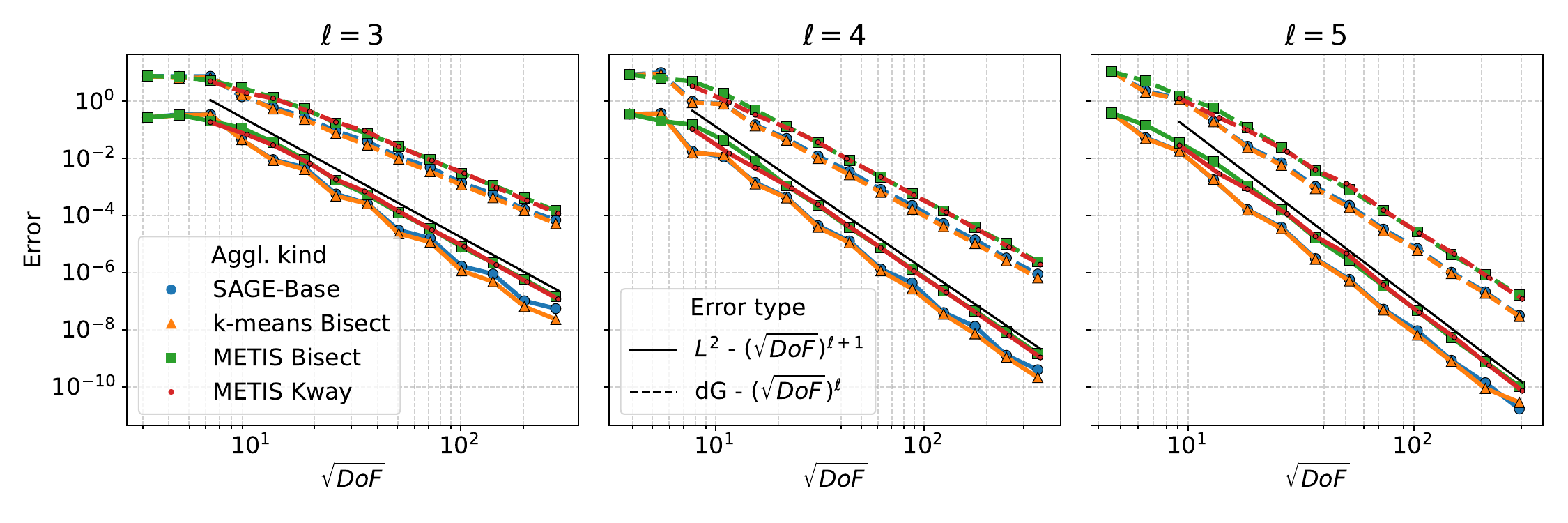}
    \caption{Verification Test 1.1: Computed errors $||u_h-u_{ex}||_{L^2(\Omega)}$ and $||u_h-u_{ex}||_{dG}$ as a function of the square root of degrees of freedom for problem \eqref{eq:Poisson_problem} solved on an agglomerate triangular grid of the unit square of $2\cdot1024^2$ elements with different agglomeration methods. }
    \label{fig:convergence_square}
  \end{figure}

  We consider a structured triangular mesh of the unit square with 1024 elements per side, of the kind shown in Figure~\ref{fig:generable_meshes} (a). We agglomerate the mesh with four different methods (METIS (k-way and bisect), k-means (bisect) and SAGE-Base); we skip the k-way version of k-means due to the large computational cost. In Figure~\ref{fig:convergence_square} we show that we obtain the expected convergence rate for all the methods. However, the higher quality mesh obtained with SAGE and k-means achieve smaller errors with the same number of degrees of freedom, showing that the improved mesh quality directly translates into better accuracy for the same number of degrees of freedom.

  \subsubsection{Poisson Problem in a Brain Slice} 
  We solve this problem on the brain slice we used in Section~\ref{sec:brain_slices} rescaled to fit in  $(0,1)^2$. The mesh is agglomerated from Matlab by calling the \py{Agglomerate.m} function as reported in Listing \ref{lst:poissonbrain}.

\begin{figure}[!t]
    \centering
    \includegraphics[width=1\textwidth]{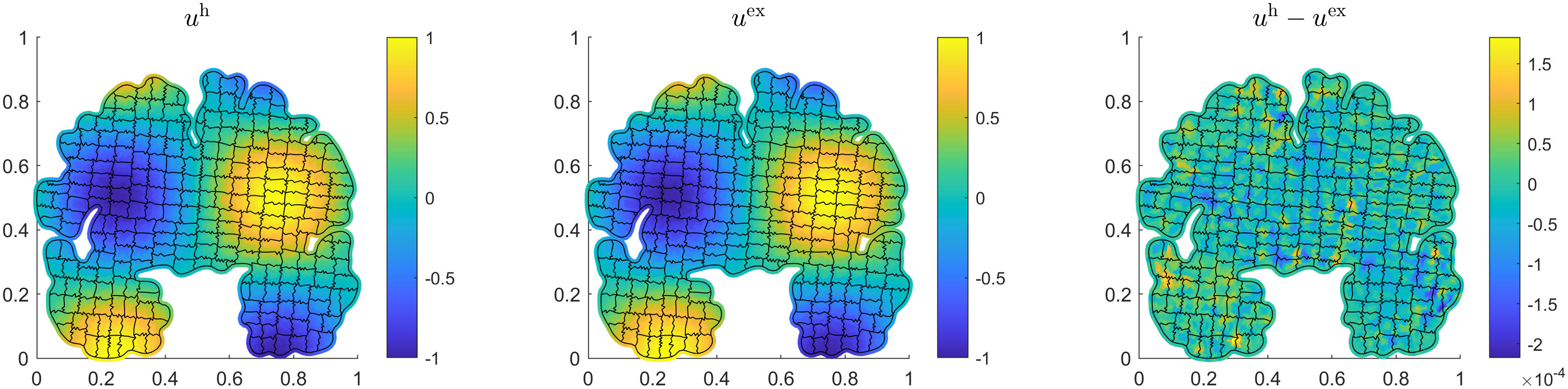}
    \caption{Verification Test 1.2: PolyDG solution of problem \eqref{eq:Poisson_problem} computed on a grid of 256 agglomerated elements (\textit{left}), plot of the exact solution (\textit{center}) and computed numerical error (\textit{right}).}
    \label{fig:polydg_poisson_brain}
  \end{figure}

  \begin{lstlisting}[language=matlab, caption={Example on how to use \maggnn's \texttt{lymph} API.}, label={lst:poissonbrain}]
    run("lymph/Physics/Laplacian/InputData/LapAggTest.m")  % creates problem Data structure
    SimType = 'laplacian';
    mesh_path = 'mesh.vtu';  % Input mesh to be agglomerated
    output_path = fullfile('AggMeshTest.mat'); % Path where the mesh will be saved     
    % agglomeration parameters
    model = 'SageBase2D';  % Agglomeration model
    mode = 'Nref';  % Agglomeration mode
    param = 7;   
    Agglomerate(mesh_path, output_path, ...
    Data, SimType, model, mode, param);
  \end{lstlisting}
  \py{Data} and \py{SimType} are needed to correctly embed the boundary conditions into the FEM region data structure. The full list of accepted parameters can be found at the \href{https://lymphlib.github.io/magnet}{documentation web page}; also, we redirect to it for details on the installation of a Python distribution compatible with Matlab. 
Once the mesh has been agglomerated, the problem can be solved by running the main function while indicating the appropriate mesh to use in the \py{Data} structure. In this case, we agglomerate the mesh using the SAGE-Base model with \py{nref=7} and solve the problem using polynomials of degree $\ell=3$. The results are reported in Figure~\ref{fig:polydg_poisson_brain}.
 
   \begin{figure}[!t]
    \centering   
    \includegraphics[width=0.9\textwidth]{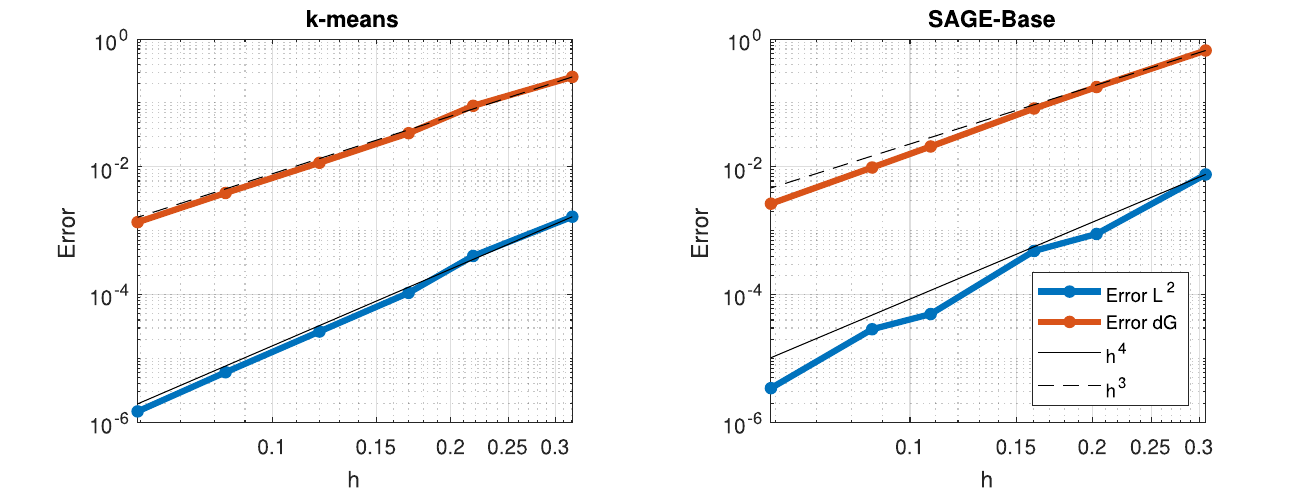}
    \caption{Verification Test 1.2: Computed errors $||u_h-u_{ex}||_{L^2(\Omega)}$ and $||u_h-u_{ex}||_{dG}$ as a function of the mesh size $h$ by fixing the polynomial degree $\ell=3$. 
    On the left, the test is performed with k-means by agglomerating different meshes so that the ratio of original elements to agglomerated elements is the same; on the right, the test is performed by agglomerating the same mesh with SAGE-Base using an increasing number of recursive bisections.}
    \label{fig:convergence_plots}
  \end{figure}
  To ensure that the agglomerated elements created by \maggnn are regular enough to achieve the theoretical order of convergence, we solve problem \eqref{eq:Poisson_problem} with the same data in the unit square $\Omega=(0,1)^2$, using six progressively finer meshes. Namely, we agglomerate the same initial mesh of 23264 elements with an increasing number of recursive bisections using the SAGE-Base model in \emph{number of refinements} mode with $N_{ref}$ going from five to ten. The polynomial degree for the discretization is fixed to $\ell=3$.
  As seen in Figure~\ref{fig:convergence_plots}, we observe the theoretical order of convergence of $\ell+1$ for the error in $L^2$-norm and of $\ell$ for the dG-norm error \cite{cangiani2021hppolyDGpolytopal}. However, we notice two facts: first, the empirical order of convergence is slightly higher than the theoretical one (4.40 and 3.34 versus 4 and 3, respectively). This is probably because in our setup, coarser meshes also have a larger number of edges per element, which is a source of irregularity and numerical errors that become less significant when we move to the finer ones. Second, the $L^2$ error seems to have an oscillatory behavior. Our suggested explanation is that in two-dimensional recursive bisection, if $N_{ref}$ is even, then the agglomerated elements will tend to be square-shaped, while if $N_{ref}$ is odd, the elements will be rectangular, with one side roughly double the other. The lower element quality for $N_{ref}$ odd corresponds to a comparatively higher error for the same $h$, as is reflected by the graph.
  To confirm these suspicions, we perform a second convergence test, this time choosing six meshes with 500 to 16000 elements and agglomerating them with k-means in \emph{direct k-way} mode, selecting $k$ to be a tenth of the original number of elements. In this way, the shape, number of edges, and overall quality of the elements will be roughly the same across all meshes. Indeed, we observe that in this case, both errors are consistently lower than in the first test, and the empirical orders of convergence are much closer to the theoretical ones (4.13 and 3.08).

  \subsection{Verification Test 2: Heat Equation with Discontinuous Boundary Conditions}
  \begin{figure}[!htbp]
    \hspace{0.12\textwidth}\textbf{$\ell=1$}
    \hspace{0.22\textwidth} \textbf{$\ell=3$} 
    \hspace{0.22\textwidth} \textbf{$\ell=5$}\\
    \centering
    \setlength{\fboxsep}{0pt}
    \setlength{\fboxrule}{1pt}
    \rotatebox{90}{\hspace{0.12\textwidth}\textbf{METIS}}
    \subfloat{\fbox{\includegraphics[trim = 235 0 275 0, clip, width=0.27\textwidth]{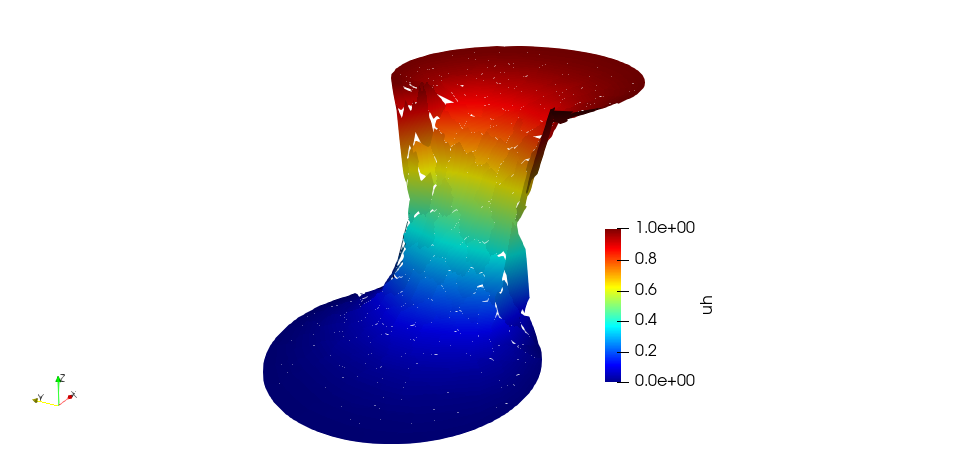}}} \quad
    \subfloat{\fbox{\includegraphics[trim = 235 0 275 0, clip, width=0.27\textwidth]{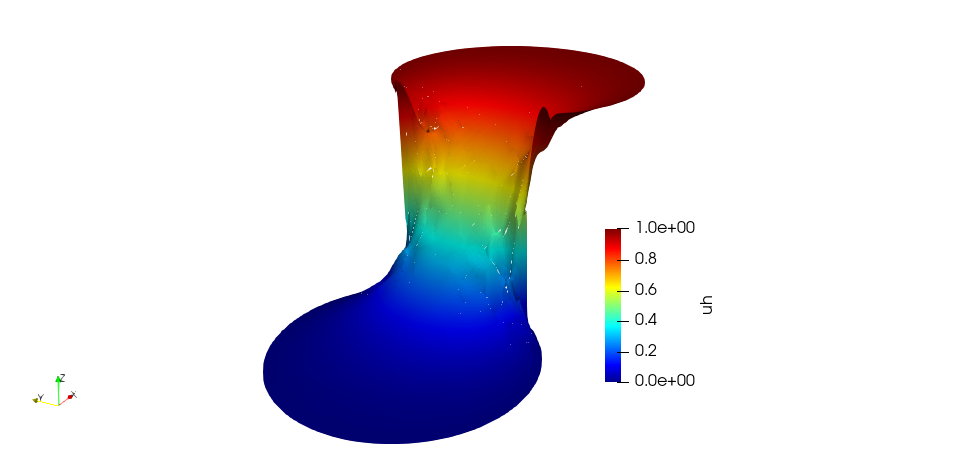}}} \quad
    \subfloat{\fbox{\includegraphics[trim = 235 0 275 0, clip, width=0.27\textwidth]{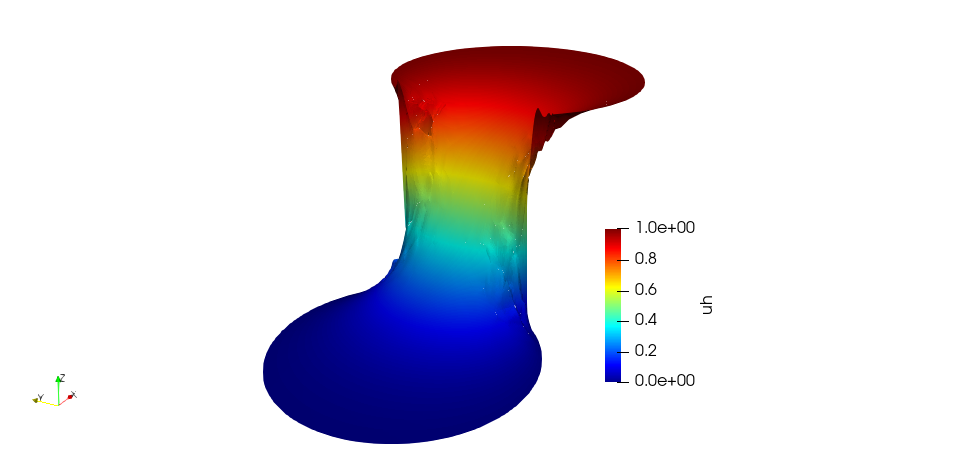}}} \\  
    \rotatebox{90}{\hspace{0.11\textwidth}\textbf{k-means}}
    \subfloat{\fbox{\includegraphics[trim = 235 0 275 0, clip, width=0.27\textwidth]{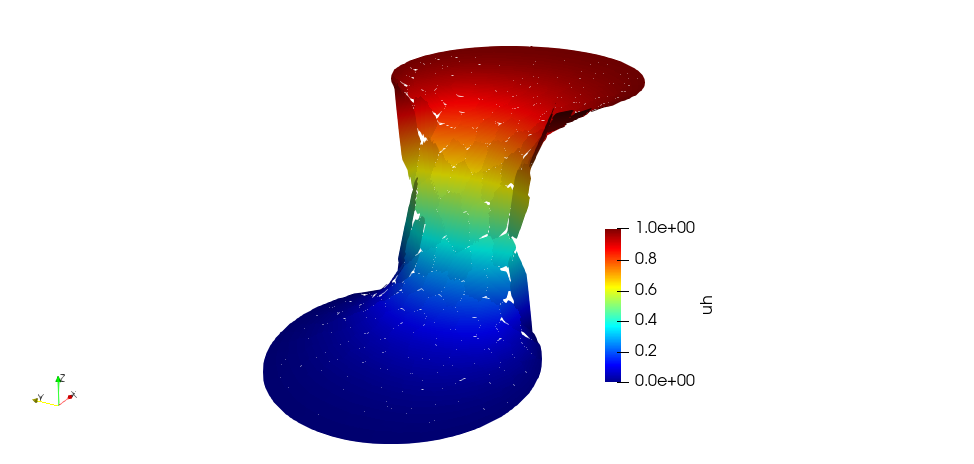}}} \quad
    \subfloat{\fbox{\includegraphics[trim = 235 0 275 0, clip, width=0.27\textwidth]{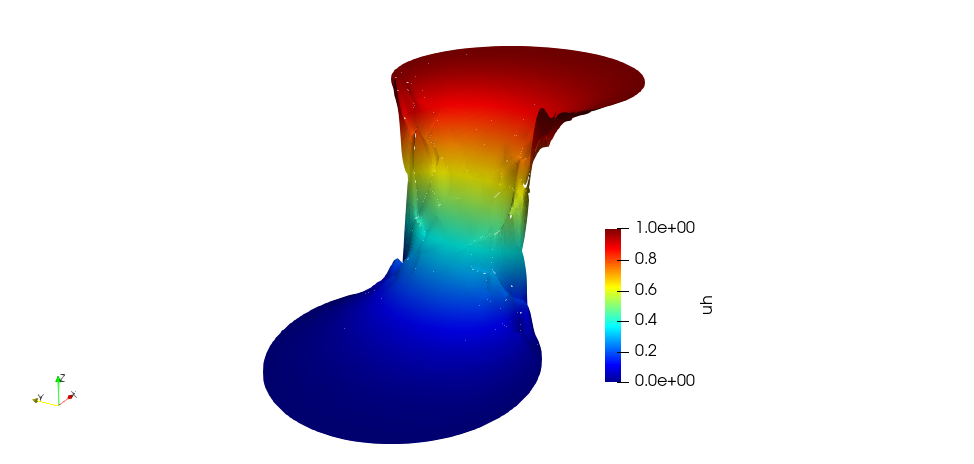}}} \quad
    \subfloat{\fbox{\includegraphics[trim = 235 0 275 0, clip, width=0.27\textwidth]{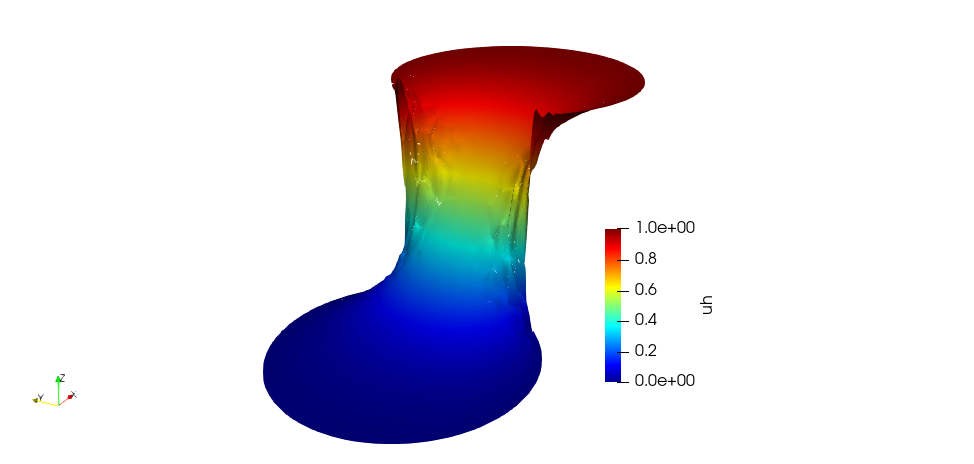}}}\\    
    \rotatebox{90}{\hspace{0.10\textwidth}\textbf{SAGE-Base}}
    \subfloat{\fbox{\includegraphics[trim = 235 0 275 0, clip, width=0.27\textwidth]{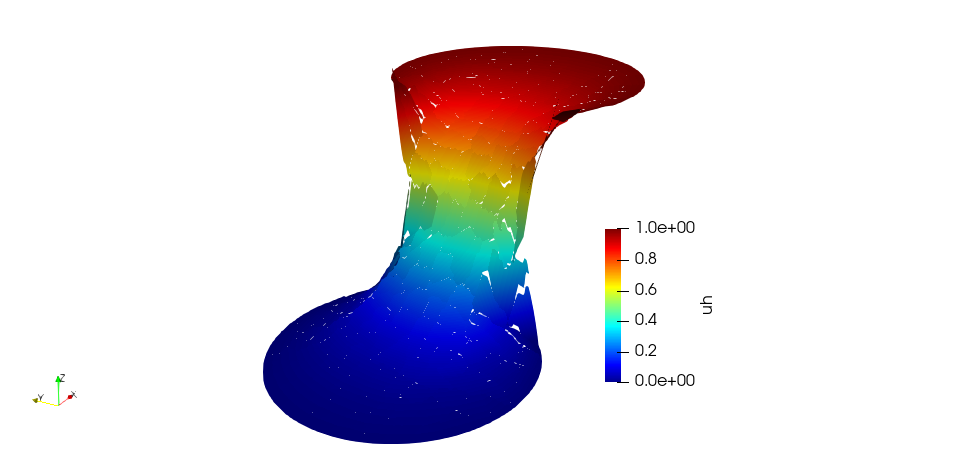}}} \quad
    \subfloat{\fbox{\includegraphics[trim = 235 0 275 0, clip, width=0.27\textwidth]{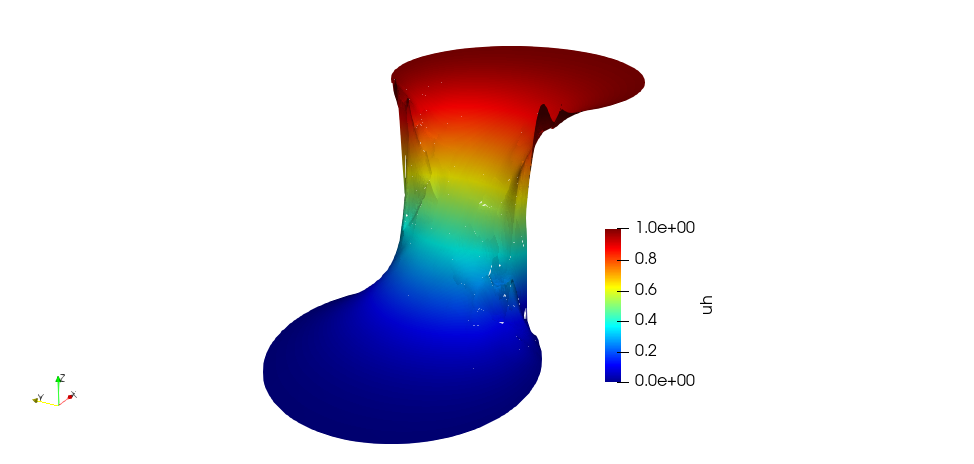}}} \quad
    \subfloat{\fbox{\includegraphics[trim = 235 0 275 0, clip, width=0.27\textwidth]{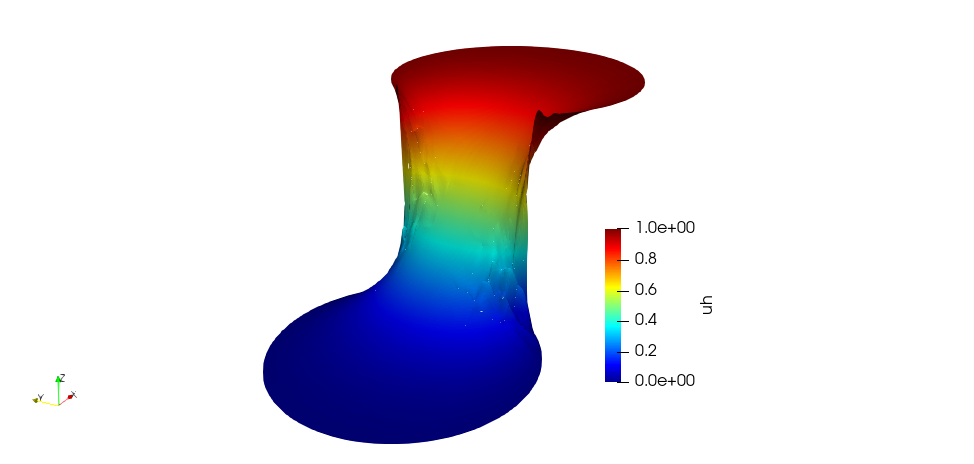}}}\\    
    \rotatebox{90}{\hspace{0.10\textwidth}\textbf{Polymesher}}
    \subfloat{\fbox{\includegraphics[trim = 235 0 275 0, clip, width=0.27\textwidth]{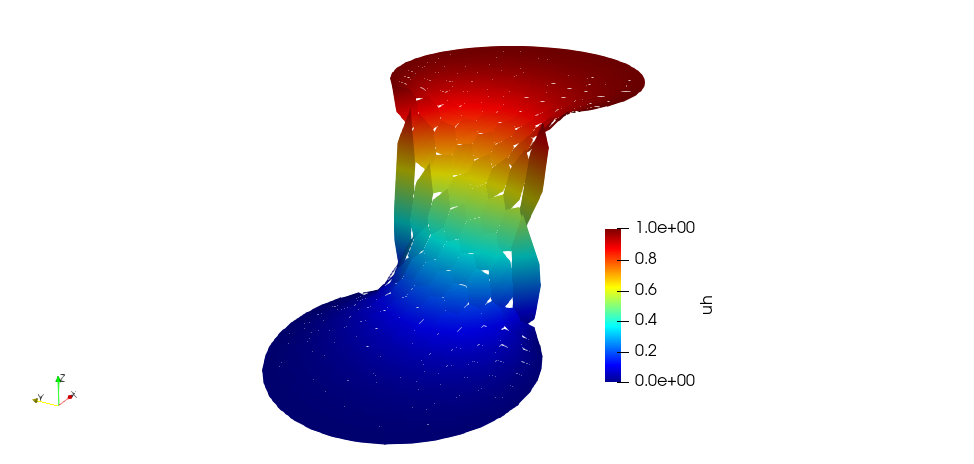}}} \quad
    \subfloat{\fbox{\includegraphics[trim = 235 0 275 0, clip, width=0.27\textwidth]{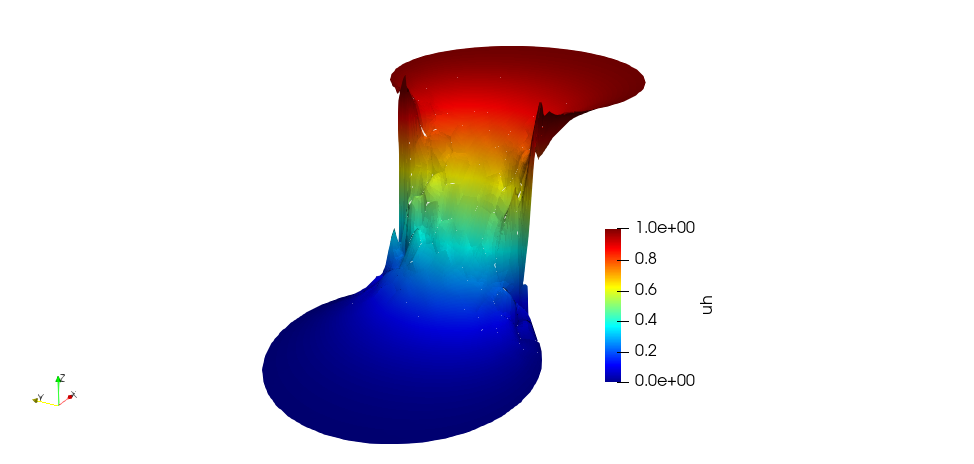}}} \quad
    \subfloat{\fbox{\includegraphics[trim = 235 0 275 0, clip, width=0.27\textwidth]{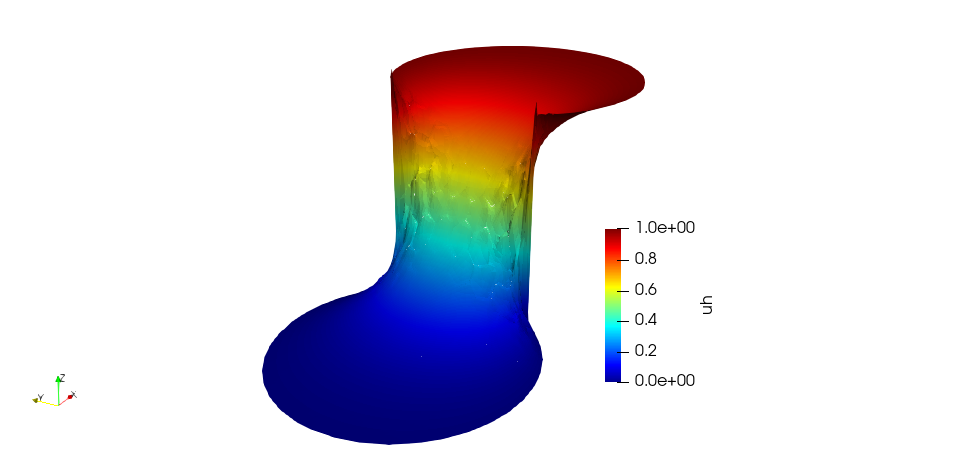}}}\\   
    \caption{Verification Test 2: PolyDG solution of problem \eqref{eq:heat_equation} on three meshes agglomerated by METIS, k-means, and SAGE-Base respectively, with the addition of a polygonal mesh generated with Polymesher, discretized with polynomial degrees $\ell=1,3,5$}
    \label{fig:heat_solutions}
  \end{figure}
  To assess the performance of differently agglomerated meshes in a PolyDG framework, we now consider the heat equation: 
  \begin{equation}
    \label{eq:heat_equation}
    \begin{cases}
      \partial_t u -\nabla\cdot (\mu \nabla u) = 0  & \text{ in } \Omega \times (0,T],\\
      u(\textbf{x},t) = g(\textbf{x},t)  & \text{ on } \partial \Omega \times (0,T],\\
      u(\textbf{x},0) = u_0(\textbf{x})  & \text{ on } \Omega \times \{t = 0\}.
    \end{cases}
  \end{equation}  
  We solve \eqref{eq:heat_equation} through \texttt{lymph},  by combining a PolyDG method with a Crank-Nicolson time integration scheme, \cite{antonietti2024lymph}.
  We consider $\Omega = B_{0.6}(0.5,0) \cup B_{0.6}(-0.5,0)$, i.e., the combination of two balls that intersect. We take $\mu=0.1$, final time $T=1$ and discontinuous Dirichlet boundary conditions $u(\mathbf{x})=1$ for $\mathbf{x}_1\ge0$, $u(\mathbf{x})=0$ for $\mathbf{x}_1<0$, with analogous discontinuous initial condition.
  We solve the problem on three meshes agglomerated by METIS, k-means, and SAGE-Base, respectively, starting from a very fine mesh of 35158 elements, using either \emph{k-way} or \emph{number of refinements} mode so that the resulting mesh has 256 elements. Additionally, we consider for comparison a polygonal mesh with the same number of elements generated using the software Polymesher \cite{talischi2012polymesher}, which tends to create very regular hexagonal cells.
  For each mesh, the problem is discretized using 3 different polynomial degrees $\ell=1,3,5$ and time step $\Delta t = 10^{-3}$. The computed PolyDG solutions at the final time $T=1$  are reported in Figure~\ref{fig:heat_solutions}.
  We observe that all four meshes produce very similar results when considering the same polynomial degree, although the agglomerated meshes have a significantly higher number of edges per element compared to the Polymesher one. The numerical solutions are overall coherent with what is expected from the literature \cite{cangiani2017hp}. For higher polynomial degrees, some numerical oscillations appear in proximity to the boundary conditions discontinuity; this might be due to agglomerated elements not perfectly conforming to the two different sections of the border. 
  Unfortunately, our current architecture does not allow for taking into account information belonging to the frontier of a domain, and it is an active topic of research.
  Indeed, if a single element shares an edge with both sections, there will be great numerical oscillations because a continuous polynomial would need to approximate a discontinuous function. Unfortunately, \maggnn has no way to guarantee this does not happen since the boundary is not taken into account when agglomerating.
  \section{Conclusions}
We have presented \maggnn, an open-source Python library for two and three-dimensional polytopal mesh agglomeration by Graph Neural Networks, and illustrated its core structure and main functionalities. Thanks to its flexibility, extensibility, straightforward integration with other software, and simple interface, we believe that \maggnn can be a useful tool to explore new ML approaches to mesh agglomeration. 
  We have introduced the different GNN methods that are available within it. Namely, the SAGE-Base model can exploit the geometric information of the mesh together with its connectivity and the reinforcement learning partitioner and refiner models, which offer enhanced performance in applications to multilevel frameworks. We tested these agglomeration strategies against state-of-the-art methods like METIS and k-means (also readily available in \maggnn), showing that they produce agglomerated meshes of comparable quality. Finally, thanks to our flexible API, we employed \maggnn in conjunction with the Matlab library \texttt{lymph} for polytopal discontinuous Galerkin approximation. We proved that the agglomerate meshes we produce are suitable for discretizing PDEs.
  We have shown that the explored machine learning-based techniques can match state-of-the-art methods in terms of performance. Moreover, we are confident that the current implementation is still improvable. By exploiting the affinity of ML techniques for modern GPU architectures, we can achieve better speed and scalability and more easily incorporate the geometric and physical features of the problem at hand.
  Future developments include integrating the recursive bisection algorithm with multigrid methods by generating a hierarchy of nested coarser grids. Including edge and face coarsening routines might also be essential to avoid having agglomerated elements with many edges, which are associated with higher computational costs.
  Regarding the GNN models, the training datasets could be extended to improve their generalization capabilities, and different model architectures should be experimented with. In particular, bigger RL partitioner and refiner models should be trained to perform better in the three-dimensional case. Also, the RL Refiner model could be improved by further exploiting the geometric information. 
  Finally, we mention that different combinations of these approaches could be possible to experiment with, e.g., using SageBase as a partitioner and a reinforcement learning model for refinement.

\section*{Declarations}
  \paragraph*{Funding.}
  This work received funding from the European Union (ERC SyG, NEMESIS, project number 101115663). 
  Views and opinions expressed are however those of the authors only and do not necessarily reflect those of the European Union or the European Research Council Executive Agency.
  Neither the European Union nor the granting authority can be held responsible for them. 

  \paragraph*{Declaration of competing interests.} 
  The authors declare that they have no known competing financial interests or personal relationships that could have appeared to influence the work reported in this paper.

  \paragraph*{Acknowledgements.}
  We thank Mattia Corti for segmenting and providing the meshes of the brain.
  The brain MRI images were provided by OASIS-3: Longitudinal Multimodal Neuroimaging: Principal Investigators: T. Benzinger, D. Marcus, J. Morris; NIH P30 AG066444, P50 AG00561, P30 NS09857781, P01 AG026276, P01 AG003991, R01 AG043434, UL1 TR000448, R01 EB009352.
  Paola F. Antonietti, Matteo Caldana and Ilario Mazzieri are members of INdAM-GNCS.
  The present research is part of the activities of “Dipartimento di Eccellenza 2023-2027”, MUR, Italy.

  \paragraph*{Code and data availability.}
  The code and data used in this work are available and can be accessed at \href{https://github.com/lymphlib/magnet}{https://github.com/lymphlib/magnet}. Proper attribution should be given when reusing or distributing the materials.

\paragraph*{Authors' contribution.}
P.F.A.: Conceptualization, Funding acquisition, Methodology, Project administration, Supervision, Writing - review and editing.
M.C.: Conceptualization, Data curation, Formal analysis, Investigation, Methodology, Software, Validation, Visualization, Writing - review and editing.
I.M.: Conceptualization, Methodology, Project administration, Supervision, Writing - review and editing.
A.R.F.: Conceptualization, Data curation, Formal analysis, Investigation, Methodology, Software, Validation, Visualization, Writing - original draft.

\bibliography{bibliography}
  
\end{document}